\RequirePackage{snapshot}
\documentclass[1p]{preprint}
\usepackage{hyperref}
\usepackage{amsmath,amssymb,amsthm}
\usepackage{subfigure}
\usepackage{graphicx}
\usepackage{color}
\usepackage{pdfsync}
\usepackage{mathtools}
\usepackage[mathscr]{eucal}
\usepackage{iga_fsi_bem_shell-private} 
\usepackage[normalem]{ulem}

\usepackage{pgfplots}
\usepackage{pgfplotstable} 
\usepackage{colortbl}
\usepackage{booktabs}
\usepackage[draft,inline,marginclue,mode=multiuser]{fixme}

\usepackage[refpage,noprefix]{nomencl}

\makenomenclature

\journal{Computer Methods in Applied Mechanics and Engineering}

\FXRegisterAuthor{ads}{aads}{ADS}
\FXRegisterAuthor{lh}{alh}{LH}
\FXRegisterAuthor{jk}{ajk}{JK}
\FXRegisterAuthor{ar}{aar}{AR}
\FXRegisterAuthor{all}{aall}{ADS+LH+JK+AR}

\newcommand{\Gm}{\vv R}

\newcommand{\K}{\vv K}
\newcommand{\Fi}{\vv P}
\newcommand{\Fe}{\vv F}

\newcommand{\ud}{\dot{\us}}
\newcommand{\udd}{\ddot{\us}}
\newcommand{\ua}{\us_{k-1+\alpha_f}}
\newcommand{\uda}{\dot{\us}_{k-1+\alpha_f}}
\newcommand{\udda}{\ddot{\us}_{k-1+\alpha_m}}

\newcommand{\OMm}{\mathcal{M}}
\newcommand{\OC}{\mathcal{C}}
\newcommand{\OFi}{\mathcal{P}}
\newcommand{\OFe}{\mathcal{F}}

\usepackage{fixfoot}
\DeclareFixedFootnote{\cofirst}{Co-first authors. These authors contributed equally to the article.} 

\begin{document}

\begin{frontmatter}

  \title{A natural framework for isogeometric fluid-structure interaction based on BEM-shell coupling}
  %\footnote{Preprint SISSA 33-2009-M.}}
  \author[sissa]{Luca Heltai\cofirst} 
  \ead{luca.heltai@sissa.it}
  \author[tubs]{Josef Kiendl\corref{cor}\cofirst} 
  \ead{j.kiendl@tu-braunschweig.de}
  \cortext[cor]{Corresponding author. Tel.: +49 0531 39194360}
  \author[sissa]{Antonio DeSimone} 
  \ead{antonio.desimone@sissa.it}
  \author[unipv,tum]{Alessandro Reali} 
  \ead{alereali@unipv.it}

\address[sissa]{SISSA-International School for Advanced Studies\\
  via Bonomea 265, 34136 Trieste - Italy}
\address[tubs]{Institute for Applied Mechanics, Technische Universit\"at Braunschweig\\
 Bienroder Weg 87, 38106 Braunschweig Germany}
\address[unipv]{Department of Civil Engineering and Architecture, University of Pavia\\
 via Ferrata 3, 27100 Pavia - Italy}
 \address[tum]{Institute for Advanced Study, Technische Universit\"at M\"unchen\\
 Lichtenbergstra{\ss}e 2a, 85748 Garching - Germany}

\begin{abstract}
  The interaction between thin structures and
  incompressible Newtonian fluids is ubiquitous both in nature and in
  industrial
  applications. In this paper we present an isogeometric formulation
  of such problems which exploits a boundary integral formulation of
  Stokes equations to model the surrounding flow, and a non linear
  Kirchhoff-Love
  shell theory to model the elastic behaviour of the structure. We
  propose three different coupling strategies: a monolithic, fully
  implicit coupling, a staggered, elasticity driven coupling, and a
  novel semi-implicit coupling, where the effect of the surrounding
  flow is incorporated in the non-linear terms of the solid solver
  through its damping characteristics. The novel semi-implicit
  approach is then used to demonstrate the power and robustness of our
  method, which fits ideally in the isogeometric paradigm, by
  exploiting only the boundary representation (B-Rep) of the thin
  structure middle surface.
\end{abstract}


\end{frontmatter}
\pagebreak
\section{Introduction}
\label{sec:introduction}

One of the most attractive features of isogeometric analysis (IGA)~\cite{HughesCottrellBazilevs-2005-a,CottrellHughesBazilevs-2009-a} is the ability to bypass mesh generation and to perform direct design-to-analysis simulations, by employing the same class of functions used for geometry parameterization in CAGD packages during the analysis process.

Most modern CAD tools, however, are based on boundary representation
(B-Rep) objects, making the use of volume-based finite element
isogeometric analysis tools (FE-IGA) less attractive, since they
require the extension of the computational domain inside (or outside)
the enclosing (or enclosed) CAGD surface.

For thin structures, isogeometric shell models circumvent this issue
since they only need a surface description of the structure. For fluid
dynamics, isogeometric boundary element methods (IGA-BEM) also
circumvent the issue mentioned above by reformulating the volumetric
flow problem in boundary integral form. Such dimensionality reduction
makes the coupling between boundary integral formulations and shell
theory an ideal combination for a large class of fluid-structure
interaction (FSI) problems, where thin structures interact with
Newtonian incompressible flows, and it fits ideally in the IGA
paradigm, by only requiring surface representations for both fluid and
structural analyses.

FSI problems have been tackled with many different techniques, ranging from interface tracking, based on Arbitrary Lagrangian Eulerian 
(ALE)~\cite{Bazilevs2008,Bazilevs201228,Dettmer20065754,DONEA1982689,Farhat199895,HUGHES1981329} 
or space--time methods~\cite{Hansbo20044195,Hobner20042087,Takizawa2011,TEZDUYAR1992339,Tezduyar1992353,Tezduyar20062002}, to interface capturing~\cite{RotundoKimJiang-2015-a} or immersed boundary (IBM)~\cite{BoffiGastaldiHeltaiPeskin-2008-a,Heltai-2008-a,Heltai2012b,Hesch201251,Zhang20042051}, or immersogeometric methods~\cite{hsu_dynamic_2015,Kamensky2014}. 
Most of the FSI applications in the isogeometric community rely on a FE-IGA approximation of the flow equations.

Some early attempts to model FSI problems using only B-Rep
representations were common in biological applications
\cite{AlougesDeSimoneHeltai-2011-a,ArroyoDeSimoneHeltai-2010-a,Arroyo30102012}
and have been recently extended to a purely B-Rep isogeometric
paradigm for the simulation of inflatable structures
\cite{Opstal2012,VanOpstal2013,VanOpstal2015}.

IGA-BEM and shell techniques have grown separately to mature and
efficient simulation techniques. On the IGA-BEM side, a considerable
effort has been put in the treatment of singular
integration~\cite{Heltai2014,Farutin2014},
adaptivity~\cite{Feischl2015,Feischl2015a}, multipatch and trimmed
surfaces~\cite{Wang2015a,Wang2015b}, efficient
solvers~\cite{TakahashiMatsumoto-2012-a,Marussig2015a}, fracture
simulations~\cite{Peng,Peng2014}, and acoustic wave
problems~\cite{Peake2015}. Flow solvers using IGA-BEM proved to be
very effective in the study of vesicles and
membranes~\cite{ArroyoDeSimoneHeltai-2010-a,Arroyo30102012,Joneidi2015},
ship hydrodynamics~\cite{KimLeeKerwin-2007-a,
  PolitisGinnisKaklis-2009-a,BelibassakisGerostathisKostas-2011-a,
  Belibassakis2013}, and rigid wings and sails~\cite{Manzoni2015},
among many others.

For shell analysis, it can be said that IGA has initiated a
renaissance in rotation-free thin shell models motivated by the high
continuity of the NURBS discretization, which permits a direct
implementation of such models. The first isogeometric formulation for
geometrically nonlinear Kirchhoff--Love shells was introduced in
\cite{kiendl_isogeometric_2009}. A similar high continuity approach
was presented in the pioneering works~\cite{cirak_subdivision_2000,
  cirak_fully_2001}. These formulations have then been
employed for various applications such as wind turbine modeling
\cite{BazHsu11a,BazHsu12b,HsuBaz12a,Korobenko:2013gw}, cloth draping
simulations \cite{lu_dynamic_2014}, explicit finite strain analysis of
membranes \cite{chen_explicit_2014}, PHT-spline shell analysis
\cite{nguyen-thanh_rotation_2011}, and fracture modeling within an
extended IGA approach \cite{nguyen_extended_2015}. Recently, this
formulation was extended to arbitrary hyperelastic materials in
\cite{kiendl_isogeometric_2015}. In the case of multipatch structures,
the lack of rotational degrees of freedom requires additional
treatment at patch connections in order to ensure the necessary
$C^1$-continuity across patch interfaces. Different methods have been
proposed, such as the bending strip method \cite{Kiendl_bending_2010},
penalty formulations as in
\cite{apostolatos_nitsche-type_2013,breitenberger_analysis_2015} or a
Nitsche formulation as in \cite{guo_nitsche_2015}. IGA has created a
lot of interest and developments not only for thin shell models but
also for Reissner-Mindlin shells
\cite{benson_isogeometric_2010,dornisch_isogeometric_2013,dornisch_treatment_2014,uhm_tspline_2009}.
Morevoer, a hierarchic family of shells was presented in
\cite{echter_hierarchic_2013} which includes Kirchhoff-Love,
Reissner-Mindlin, and higher order shells. In
\cite{benson_blended_2013}, a blended shell formulation was presented,
which is a combination of rotation-free shells
\cite{benson_large_2011} in the patch interior and Reissner-Mindlin
shells \cite{benson_isogeometric_2010} at the boundary. Furthermore,
many developments have been done as well on isogeometric solid shells
\cite{bouclier_efficient_2013,caseiro_on_2014,caseiro_assumed_2015,
  hosseini_isogeometric_2013,hosseini_isogeometric_2014}. Thin shell
models, as those presented in
\cite{kiendl_isogeometric_2009,kiendl_isogeometric_2015}, are purely
surface-based, in a sense that the shell is completely defined by its
middle surface and the shell kinematics are completely described by the
middle surface metric and curvature properties. This allows for a direct
integration of IGA into CAD systems
\cite{breitenberger_analysis_2015,schmidt_realization_2010}, it
facilitates the coupling of shell structures and fluids in
fluid-structure interaction (FSI) applications due to the lack of
rotational degrees of freedom
\cite{bazilevs_3d_2011,BazHsu12b,HsuBaz12a}, and it is perfectly
suited for the coupling with an IGA-BEM fluid solver as we propose in
this paper.

The rest of this paper is organized as follows. In
 section~\ref{sec:bound-integr-appr} we present a brief overview of
isogeometric NURBS spaces. The continuous FSI problem we want to
tackle is introduced in section~\ref{sec:fluid-struct-inter}, and an
analysis of the coupling strategies is presented in
section~\ref{sec:algor-analys-fsi}. Sections~\ref{sec:hydr-inter}
and~\ref{sec:KLtheory} present respectively the numerical treatment of
the BEM and Shell parts, while in sections~\ref{sec:numtests} and
\ref{sec:conclusions} we present some numerical examples and draw some
conclusions.

\section{Overview of isogeometric NURBS spaces}
\label{sec:bound-integr-appr}

Given a nondecreasing knot vector
$\Theta = \{k_0, k_1, ..., k_{n+p}\}$, \nomenclature{$\Theta$}{Knot
  span for B-splines} the $n$ \nomenclature{$n$}{Number of univariate
  B-spline basis} B-splines of degree $p$
\nomenclature{$p$}{Degree of B-spline basis (as superscript)} are
defined by the recurrence relation
\begin{equation}
  \label{eq:B-Spline-p0-definition}
    B^{(i,0)}(s) = 
    \begin{cases}
      1, & \mbox{ if } k_i \leq s < k_{i+1}\\
      0, & \mbox{ otherwise},
    \end{cases}
  \end{equation}
for $p=0$, while for $p>0$ we have
\nomenclature{$B^{(i,p)}$}{The $i$-th B-spline basis of order $p$}
\begin{equation}
  \label{eq:B-Spline-definition}
  B^{(i,p)}(s) = \tau^{(i,p)}(s) B^{(i,p-1)}(s) - \tau^{(i+1,p)}(s) B^{(i+1,p-1)}(s),
\end{equation}
for $i = 0, ..., n-1$, where
\begin{equation}
  \label{eq:def-tau}
  \tau^{(i,p)}(s) := 
  \begin{cases}
    \displaystyle
      \frac{s - k^i}{k^{i+p} - k^i} & \mbox{ if } k^{i+p} \neq k^{i}\\
    0, & \mbox{ otherwise}.
  \end{cases}
\end{equation}
\nomenclature{$\tau^{(i,p)}$}{Interpolating function for B-spline
  recursive definition}

The above recurrence relation can be evaluated in a numerically stable
way by the de Boor algorithm (see, for example,
\cite{PieglTiller-1997-a}). Between two distinct knots, a B-spline is
of continuity class $\mathcal C^{\infty}$, at a single knot it is
$\mathcal C^{p-1}$, and, if a knot is repeated $q$ times, the continuity
is reduced to $\mathcal C^{p-q}$. A knot can be repeated at most
$q=p+1$ times resulting in a discontinuity ($\mathcal C^{-1}$) at that
location.

NURBS basis functions are readily obtained from B-Splines by assigning
a positive weight $w_i$ to each basis spline function and defining the
corresponding NURBS basis function as
\begin{equation}
  \label{eq:definition-NURBS}
  N^i(s) := \frac{w^iB^{(i,p)}(s)}{\sum_{j=0}^{n-1} w^j B^{(j,p)}(s)}.
\end{equation}
\nomenclature{$N^i$}{The $i$-th NURBS basis function of order $p$}

Notice that also the NURBS basis have the partition of unity property,
and B-Splines can be considered a special case of NURBS by taking all
weights to be identical. Taking two knot vectors $\Theta_i$, with
$i=0,1$, one can construct the NURBS basis functions for
two-dimensional surfaces embedded in three dimensional space by tensor
products. Indicating with $\vv s := [s_0,s_{1}] $ a point in $\Re^2$
and $\vv i := (i_0,i_1)$ a two dimensional multi-index belonging to
the set \nomenclature{$\vv i, \vv j$}{Multi-indices (bold latin),
  where each
  component $i_k, j_k$ is in the range $[0, n_k)$}
\nomenclature{$n_k$}{Number of B-spline basis functions of the
  component $k$ for multi-variate B-splines}
\nomenclature{$\mathcal J$}{Set of all possible multi-indices for
  scalar basis functions}
\begin{equation}
\label{eq:multi-index-set}
 \mathcal J := \{\vv j =(j_0, j_1), \quad 0\leq j_k < n_k,
 \quad k = 0, 1\},
\end{equation}
the bi-variate B-Splines and NURBS basis functions are given by
\begin{equation}
  \label{eq:definition-splines-multivariate}
  B^{\vv i, \vv p}(\vv s) := B^{(i_0, p_0)}(s_0)  B^{(i_1, p_1)}(s_1), \qquad N^{\vv i}(\vv
  s) := \frac{w{\vv i} B^{(\vv i, \vv p)}(\vv s)}{\sum_{\vv j \in
      \mathcal J} w^{\vv j} B^{(\vv j,\vv p)}(\vv s)},
\end{equation}
where $\vv i$, $\vv p$ and $\vv j$ are all multi-indices. The
multi-index $\vv p = (p_0, p_1)$ is used to keep track of the degrees
of the B-Splines in each direction, while $\vv n = (n_0, n_1)$ is used
to keep track of the number of basis functions in each
direction. Notice that in Equations~\eqref{eq:definition-NURBS}
and~\eqref{eq:definition-splines-multivariate} we dropped the
superscripts $p$ and $\vv p$ from the definition of the NURBS basis
functions $N^{\vv i}$, to ease the notation in the rest of the paper.

As a generalization of the one dimensional case, if we take a
collection of $n := n_0n_1$ control points in $\Re^3$, we can
represent a two dimensional manifold in a three dimensional space as
the image of the map
\begin{equation}
  \label{eq:definition-m-manifold}
  \Re^3 \supset {\vv x}(\vv s) := \sum_{\vv i\in \mathcal J}\vv P^{\vv
    i}
  N^{\vv i}(\vv s) \qquad \vv s \in \Re^2.
\end{equation}
\nomenclature{$B^m$}{Reference domain of m-dimensional patches}
\nomenclature{$\vv s$}{Point in the reference domain $B^m$} 

The set of control points $\vv P^{\vv i}$ with $\vv i \in \mathcal J$
is usually referred to as \emph{control net}.
The domain of the map ${\vv x}(\vv s)$ is the set
\begin{equation}
  \label{eq:domain-of-map}
  B^2 := [k^0_{0}, k^0_{n_0+p_0}] \times
  [k^{1}_{0}, k^{1}_{n_{1}+p_{1}}] \subset \Re^2,
\end{equation}
where $k^i_j$ is the $j$-th knot in the $i$-th knot vector $\Theta_i$.

In what follows, we will use greek indices $\alpha,\beta$ to indicate
components in the two dimensional manifold (i.e., from zero to one)
and latin indices to indicate components in the three dimensional
embedding manifold.

The tangential vectors on a point on the surface are given by the
covariant base vectors $\vv g_{\alpha}$:
\begin{align}
  \vv g_{\alpha}(\vv s) = \frac{\partial \vv x(\vv s)}{\partial s^{\alpha}} =  \vv x_{,\alpha}(\vv s)  \label{eq:base_vec}
\end{align}
\nomenclature{$\alpha, \beta$}{Greek indices, running from zero to one}
\nomenclature{$a, b, i, j, k$}{Latin indices running from one to
  three}
\nomenclature{$\vv g_\alpha$}{Covariant base vectors of the deformed
  configuration}
\nomenclature{$\vv g^\alpha$}{Contravariant base vectors of the
  deformed configuration}
\nomenclature{$\delta^\alpha_\beta, \delta^i_j$}{Kronecker deltas}

Contravariant base vectors $\vv g^{\alpha}(\vv s)$ are obtained through the relation $\vv g^{\alpha}\cdot \vv g_{\beta}=\delta^{\alpha}_{\beta}$, where $\delta^{\alpha}_{\beta}$ is the Kronecker delta.
Furthermore, we introduce the unit normal vector $\vv g_{3}$:
\begin{align}
  \vv g_{3}(\vv s) &= \frac{\vv g_{\alpha}(\vv s) \times \vv g_{\beta}(\vv s)}{|\vv g_{\alpha}(\vv s) \times \vv g_{\beta}(\vv s)|}. \label{eq:normal_vec}
\end{align} 
With the tangential and normal vectors, we can write the first and second fundamental forms of the surface, respectively:
\begin{align}
  g_{\alpha\beta}(\vv s) &= \vv g_{\alpha}(\vv s) \cdot \vv g_{\beta}(\vv s) \label{eq:first_fund2} \\
  b_{\alpha\beta}(\vv s) &= \vv g_{\alpha,\beta}(\vv s) \cdot  \vv g_{3}(\vv s) \label{b_ab}
\end{align}
where $g_{\alpha\beta}$ and $b_{\alpha\beta}$ represent the metric and curvature coefficients of the surface. 
\nomenclature{$\vv g_3$}{Unit normal vector of the deformed configuration}
\nomenclature{$g_{\alpha\beta}$}{First fundamental form of the deformed configuration}
\nomenclature{$b_{\alpha\beta}$}{Second fundamental form of the deformed configuration}

Integrals on the two-dimensional manifold $\vv x(B^2)$ can be pulled
back to the domain $B^2$ using the standard transformation rule
\begin{equation}
  \label{eq:integral-transformation-rule}
  \int_{\vv x(B^2)} f(\vv x) \d A = \int_{B^2} f(\vv x(\vv s))
  J(\vv s) \d \vv s,
\end{equation}
where we indicated with $J(\vv s)$ the square root of the determinant
of the first fundamental form:
\begin{equation}
  \label{eq:det-J}
  J(\vv s) := \sqrt{\det(g_{\alpha\beta}(\vv s))}. 
\end{equation}
\nomenclature{$J$}{Square root of the determinant of $g_{\alpha\beta}$ }

%For the sake of simplicity, we expose the theory and the numerical
%approximation only considering \emph{single patch geometries},
%referring to~\cite{CottrellHughesBazilevs-2009-a} for the treatment of
%multi-patch domains, with different level of regularity across
%different patches, as well as for indications on how to treat
%different refinements on each patch.
%
A standard (scalar) isogeometric finite dimensional space on a
two-dimensional manifold is readily obtained by considering the span
of the functions $\phi^{\vv i} := N^{\vv i}\circ\vv x^{-1}$:
\begin{equation}
  \label{eq:finite-dimensional-iga-basis}
  V_h := \text{span}\{ \phi^{\vv i}(\vv y) \}_{\vv i \in
    \mathcal J}, \qquad \vv y \in \vv x(B^2) \subset \Re^3,
\end{equation}
where $\phi^{\vv i}$ are such that
\begin{equation}
  \label{eq:inverse-phi}
  \phi^{\vv i}(\vv x(\vv s)) = N^{\vv i}(\vv s),\qquad \forall \vv s \in B^2.
\end{equation}
\nomenclature{$\phi^{\vv i}$}{Scalar basis function for the space $V_h$}
\nomenclature{$V_h$}{Scalar isogeometric finite dimensional space}

The dimension of the space $V_h$ is $n = n_0 n_1$ and it is equal to
the number of control points that define the geometry of the problem.
If we introduce the multi-index set $\mathcal J^3$, as done for scalar
functions in equation~\eqref{eq:multi-index-set},
\begin{equation}
\label{eq:multi-index-set-vector}
 \mathcal J^3 := \{\vv j =(j_0, j_1, j_2), \quad 0\leq j_k < 3\times n_k,
 \quad k = 0, 1, 2 \},
\end{equation}
then a finite dimensional space for vector fields of three components
is obtained by considering
\begin{equation}
  \label{eq:finite-dimensional-vector-iga-basis}
  V^3_h := \text{span}\{ \vv \Phi^{\vv I}(\vv y) \}_{\vv I \in
    \mathcal J^3}
  \qquad \vv y \in \vv x(B^2) \subset \Re^3,
\end{equation}
\nomenclature{$V^3_h$}{Vector isogeometric finite dimensional space
  with three components}
where the basis functions $\vv \Phi^{\vv I}$ are such that
\begin{equation}
  \label{eq:vector-basis-functions}
  \vv \Phi^{\vv I} (\vv y) := \vv e_a \phi^{\vv j}(\vv y), \qquad
  \vv I = (3 j_0+a, 3j_{1}+a) = 3 \vv j +a.
\end{equation}
\nomenclature{$\vv I, \vv J, \vv K$}{Multi-indices
  (uppercase bold latin) for vector valued functions of three
  components, where each component $I_k, J_k$ is in the range $[0,
  3*n_k)$} 
\nomenclature{$\mathcal J^3$}{Set of all possible multi-indices for
  vector basis functions}
\nomenclature{$\vv \Phi$}{Vector basis function for the space $V^3_h$}

The multi index $\vv I$ is meant to transform the multi-index
$\vv j \in \mathcal J$ plus the component index $a$ into a unique
global identifier for the $\vv I$-th basis function. In what follows,
we use upper case bold latin indices $\vv I, \vv J$ to indicate the
global numbering of the basis functions defining the space $V^3_h$,
lower case latin indices $i,j,k$ to label spacial coordinates in
$[0,d)$ and greek indices $\alpha,\beta$ to label parameter
coordinates in $[0,1]$. Unless otherwise stated, we use Einstein
summation convention. A vector function of three components
$\vv f(\vv x)$ in the space $V_h^3$ is identified by its coefficient
vector $\vv f$ such that \begin{equation}
  \label{eq:def-finite-f}
  V^3_h \ni \vv f(\vv x) := \vv f^{\vv I} \vv \Phi^{\vv I}(\vv x),
\end{equation}
where, with a slight abuse of notation, we denote the vector of
coefficients $\vv f$ with the same symbol as the function $\vv f(\vv
x)$ but without the argument ``$(\x)$''.

\section{Fluid-Structure Interaction}
\label{sec:fluid-struct-inter}

We are interested in studying the interaction between a thin
deformable elastic body and an incompressible fluid. We consider a
model problem where the inertial terms of the fluid are negligible
when compared with both the fluid viscosity and the inertial terms of
the solid.  We assume that the deformable body occupies at time $t$
the region $\Os(t) = \vv x(\Os_0)\subset \Re^3$ and that the rest of the space is
entirely occupied by an incompressible fluid whose time dependent
domain is $\Omega^\fluid(t) = \Re^3\setminus \Os(t)$.

The fluid and solid domains are coupled through non-slip conditions
and through balance equations across the boundary of the solid domain
$\Gfsi(t):= \partial \Os(t)$. We will describe the fluid equations in
Eulerian form, where the primal variables are the velocity field of
the fluid $\vf$ and its pressure $p$ at fixed points in space, while
we use a Lagrangian description for the the solid, where each
point $\vv X$ represents a fixed material point in $\Os_0$ mapped by
the transformation $\vv x:\Omega_0 \times [0,T] \mapsto \Re^d$ to its
current location $\vv x(\X,t)$ at time $t$.

For convenience, we introduce the deformation field $\us(\X,t)$, such
that $\us(\X,t)=\vv x(\X,t)-\X$.  The transformation map
$\vv x(\cdot,t): \Os_0 \mapsto \Os(t)$ is assumed to be invertible and
bi-lipschitz for each time $t$ in the interval $[0,T]$, i.e., the
determinant $J$ of the deformation gradient
$\deformation := \Grad ~\vv x(\vv X,t) := \nabla_{\X} \vv x(\vv X,t) =
\nabla_{\X} \us + \Id$
is strictly positive and $\deformation$ is bounded.  We denote the
gradient and the divergence with respect to the $\X$ variable with
$\Grad$ and $\Div$.

The equations of motion of the system can be written as:
\begin{subequations}
  \label{eq:incompressible-navier-stokes}
  \begin{alignat}{3}
    \label{eq:conservation-momentum-fluid}
    & -\nabla\cdot\cauchyf := -\eta \Delta \vf +\nabla p = 0\qquad & \text{ in } \Of(t)\\
    \label{eq:incompressibility-fluid}
    & \nabla \cdot \vf = 0 & \text{ in } \Of(t)\\[.3cm]
   \label{eq:conservation-momentum-solid-material}
    & \rhos  \frac{\partial^2 \us}{\partial t^2} - \Div \big(\deformation
    \cdot \secondpiola \big) -\vv b = 0 \qquad & \text{ in } \Os_0\\[.3cm]
    \label{eq:stress-coupling-material}
    & J\cauchyf \cdot\deformation^{-T} \cdot \vv \nu_0 =  
    \deformation\cdot \secondpiola\cdot\vv \nu_0 & \text{ on } \Gfsi_0\\
    \label{eq:non-penetration-condition-material}
    &  \vf(\vv x(\X,t), t) = \frac{\partial \us(\X,t)}{\partial t} =:
    \dot \us(\X,t) 
    \qquad&
    \text{ on } \Gfsi_0\\[.3cm]
    & \us|_{t=0} = \us_0  & \text{ in } \Os_0\\
    \label{eq:initial-conditions-v-solid}
    & \dot \us |_{t=0} = \vs_0  & \text{ in } \Os_0.
  \end{alignat}
\end{subequations}
\nomenclature{$\Of$}{Fluid domain}
\nomenclature{$\Os$}{Solid domain}
\nomenclature{$\Gfsi$}{Boundary of the solid domain}
\nomenclature{$\Os_0$}{Solid reference domain}
\nomenclature{$\vv X$}{Solid material point}
\nomenclature{$\vf$}{Fluid velocity}
\nomenclature{$\eta$}{Fluid viscosity}
\nomenclature{$\rhos$}{Solid density}
\nomenclature{$\us$}{Solid displacement}
\nomenclature{$\dot \us$}{Solid velocity}
\nomenclature{$\ddot \us$}{Solid acceleration}
\nomenclature{$\deformation$}{Solid deformation gradient}
\nomenclature{$\vv S$}{Second Piola-Kirchhoff stress tensor}
\nomenclature{$\cauchyf$}{Cauchy stress tensor for the fluid}
\nomenclature{$\vv \nu_0$}{Outer normal to the reference configuration}
\nomenclature{$p$}{Fluid pressure}
\nomenclature{$\vs_0$}{Initial solid velocity}
\nomenclature{$\us_0$}{Initial solid displacement}

Where $\eta, \cauchyf$ are the fluid viscosity and Cauchy stress
tensor respectively, while $\rhos, \vv S$ are the solid density and
second Piola-Kirchhoff stress tensor, respectively, and $\vv b$ is a
body load acting on the solid, i.e., gravity. The quantities
$\vs_0, \us_0, \vv \nu_0$ are the initial solid velocity, initial solid
displacement, and outer normal to the reference configuration.

We remark here that at low Reynolds numbers time dependency in the
equations of motion of the fluid can only occur due to boundary
conditions and through the time dependent changes in the shape of the
domain (i.e., $\Of(t)$). For this reason, in
equation~\eqref{eq:incompressible-navier-stokes} there are no initial
conditions for the fluid velocity, which is assumed to adjust
instantaneously to changes in boundary conditions and in domain shape.

Equations~\eqref{eq:conservation-momentum-fluid}
and~\eqref{eq:incompressibility-fluid} represent the conservation of
momentum and mass in Eulerian form for a low Reynolds number flow,
while equation~\eqref{eq:conservation-momentum-solid-material} is the
conservation of momentum for a solid body, written in Lagrangian form.

We will restrict our attention to problems for which $\Os(t)$ is a
thin shell and we will consider the Kirchhoff-Love shell theory, where
%transverse shear deformation is neglected and 
the director, i.e., a vector normal to the middle surface, is assumed
to remain normal to the middle surface in the deformed configuration
(i.e., parallel to $\vv g_3$). With this assumption, the configuration
of the shell is uniquely determined once we know the configuration of
its middle surface $\Gamma(t)$, making this an ideal candidate for a
coupled FSI problem which requires only a surface description.

To summarise, here are the list of all assumptions we make in our model:
\begin{itemize}
\item the inertial terms of the fluid are negligible when compared
  with both the fluid viscosity and the inertial terms of the solid;
\item the transversal dimension $h$ of the solid is much smaller than
  all other directions, and can be neglected when considering the
  geometry of the problem;
\item the coupling conditions between the solid and the fluid are
  applied \emph{at the middle surface} $\Gamma(t)$ of the solid.
\end{itemize}

With these assumptions, the fluid equations reduce to
Stokes equations on the domain $\Re^3\setminus \Gamma(t)$. For a
given prescribed velocity on $\Gamma(t)$, we can compute the force
per unit area that the fluid exerts on the middle surface of the
solid, by pulling the \emph{jump} of the fluid normal stress on $\Gamma(t)$ back to the solid reference configuration.

We define the operator that performs this pull back $\DNt$, i.e., a
\emph{Dirichlet to Neumann map} such that:
\nomenclature{$\DNt$}{Dirichlet to Neumann map for the fluid system}
\begin{equation}
  \label{eq:dn-map}
  \vv f^\fsi = J \DNt \vs_g.
\end{equation}

Given a Dirichlet datum $\vs_g$ on the middle surface $\Gamma(t)$, this
returns the pull back of the jump of the normal stress associated with
the solution of the fluid problem, i.e.,
\begin{subequations}
  \label{eq:fsi-stokes-shell}
  \begin{alignat}{3}
    \label{eq:conservation-momentum-fluid-shell}
     - \nabla\cdot\cauchyf := &-\eta \Delta \vf + \nabla p = 0 \qquad
    & \text{ in } \Re^3\setminus \Gamma(t)\\
    \label{eq:incompressibility-fluid-shell}
    & \nabla \cdot \vf = 0 & \text{ in } \Re^3\setminus \Gamma(t)\\[.3cm]
    \label{eq:non-slip}
     & \vf = \vs_g& \text{ on } \Gamma(t)\\
    \label{eq:conservation-momentum-solid-shell}
     & J^{-1}\vv f^\fsi = \jump{\cauchyf} \cdot \vv g_3& \text{ on } \Gamma(t),
  \end{alignat}
\end{subequations}
where the symbol $\jump{\cauchyf}$ represents the difference between
$\cauchyf$ across the middle surface $\Gamma(t)$.
\nomenclature{$\jump{.}$}{Jump of a quantity across the middle surface}

For the structural analysis, we consider a Kirchhoff-Love shell 
in large deformations and small strains, i.e., a St.-Venant-Kirchhoff
material model is applied. In the following, we present the weak form
of the problem, based on the principal of virtual work. 

We will indicate with $\vv x(\vv s, t)$ the current configuration of
the middle surface with respect to the curvilinear coordinates
$\vv s$, and with $\vv X(\vv s)$ the middle surface in the undeformed
configuration.  Analogously to equations
\eqref{eq:base_vec}-\eqref{b_ab}, we define the tangent vectors
$\vv G_\alpha$, the unit normal vector $\vv G_3$, the metric
coefficients $G_{\alpha\beta}$ and the curvature coefficients
$B_{\alpha\beta}$ for the undeformed configuration $\vv X(\vv s)$.
\nomenclature{$\vv G^{\alpha}$}{Contravariant base
  vectors of the reference configuration}
\nomenclature{$\vv G_{\alpha}$}{Covariant base vectors of the
  reference configuration}

As strain measure we use the the Green-Lagrange strain tensor, where
only in-plane strains are considered,
$\vv E=E_{\alpha\beta}\,\vv G^{\alpha} \otimes
\vv G^{\beta}$, with: \nomenclature{$\vv E$}{Green-Lagrange
  strain tensor}
\begin{align}
  E_{\alpha\beta} &= \varepsilon_{\alpha\beta}+\theta^3\kappa_{\alpha\beta}  \label{E_ab} \\
  \varepsilon_{\alpha\beta} &= \frac{1}{2}(g_{\alpha\beta}-G_{\alpha\beta})  \label{eps_ab} \\
  \kappa_{\alpha\beta} &= B_{\alpha\beta} - b_{\alpha\beta}  \label{kap_ab}
\end{align}
where $\varepsilon_{\alpha\beta}$ represents the membrane strain while
\nomenclature{$\varepsilon_{\alpha\beta}$}{Membrane strain} $\kappa_{\alpha\beta}$
describes the change in curvature or bending (pseudo-)strain. 
\nomenclature{$\kappa_{\alpha\beta}$}{Bending (pseudo-)strain} 
As stress measure, we use the energetically
conjugate second Piola-Kirchhoff stress tensor $\vv S$:
\begin{align}
  \vv S &= \mathbb{C}:\vv E  \label{S^ab}
\end{align}
where $\mathbb{C}$ is the fourth-order material tensor. 
\nomenclature{$\mathbb{C}$}{Fourth-order material tensor}
%For the shell formulation, both strains and stresses are separated in membrane and bending action, and the 
Stresses are represented by the stress resultants $\vv n$ and
$\vv m$, which are the normal forces and bending moments,
respectively. 
\nomenclature{$\vv n$}{Normal forces}
\nomenclature{$\vv m$}{Bending moments}
They are obtained by integrating the constant and the linear parts separately through the shell thickness $h$ as follows:
\nomenclature{$h$}{Shell thickness}
\begin{align}
  \vv n &= \int_{-h/2}^{h/2}\vv S(s^3=0) \mathrm d s^3 = h\; \hat{\mathbb{C}}:\boldsymbol{\varepsilon} \label{n^ab} \\
  \vv m &= \int_{-h/2}^{h/2}(\vv S-\vv S(s^3=0))\;s^3 \mathrm ds^3 = \frac{h^3}{12}\; \hat{\mathbb{C}}:\boldsymbol{\kappa}, \label{m^ab}
\end{align}
where $\hat{\mathbb{C}}$ is the plane stress material tensor \cite{bischoff_models_2004}. 
With membrane strains \eqref{eps_ab}, change in curvature
\eqref{kap_ab}, normal forces \eqref{n^ab}, and bending moments
\eqref{m^ab}, the internal virtual work of the shell can be written as:
\begin{align}
  \delta W_{int} = - \int_{\Gamma_0} \left(\rho \ddot \us \cdot\delta\us  +
    \vv n:\delta\boldsymbol{\varepsilon}+\vv m:\delta\boldsymbol{\kappa}
    \right)
  \mathrm dA, \label{W_intKL}
\end{align}
where $\delta$ denotes that these variables derive from a virtual
displacement $\delta \us$ and $\mathrm dA$ is the differential
area element of the middle surface.

The external virtual work is defined as:
\begin{align}
  \delta W_{ext}&=  \int_{\Gamma_0}  \left(\vv f^\fsi + \vv b\right) \cdot\delta\us \, \mathrm dA,
\end{align}
where $\vv f^\fsi$ is the term coming from the fluid-structure
interaction and $\vv b$ some additional body load (e.g., gravity)
acting on the shell.

The system is in equilibrium if the sum of internal and external
virtual work vanishes
\begin{align}
  \delta W_{int} + \delta W_{ext} = 0 \label{dW},
\end{align}
which must hold for an arbitrary variation of $\delta\us$. 

The final fluid-structure interaction system is given by
\begin{subequations}
  \label{eq:final-fsi-system}
  \begin{alignat}{3}
     \int_{\Gamma_0} \bigg(  &
     \rho \ddot \us \cdot\delta \us 
     -   J \left( \DNt \dot \us \right)\delta \us 
     \nonumber \\
     & \vv n:\delta\boldsymbol{\varepsilon} +
     \vv m:\delta\boldsymbol{\kappa} 
     - \vv b \cdot \delta \us 
     \bigg) \d A = 0,  \qquad && \forall \delta \us \in
    H^2(\Gamma_0)  \\[.3cm]
    \label{eq:intial-condition-u-final}
    & \us|_{t=0} = \us_0  && \text{ on } \Gamma_0\\
    \label{eq:initial-condition-ut-final}
    & \dot \us|_{t=0} = \vs_0  && \text{ on } \Gamma_0.
  \end{alignat}
\end{subequations}

An explicit construction of the Dirichlet to Neumann operator $\DNt$
is given in section~\ref{sec:hydr-inter}.
In the general case, we allow the solid to be either free, hinged or
clamped. In the first case, the functional space of virtual
displacement for which the variational formulation
\eqref{eq:final-fsi-system} makes sense is
\begin{equation}
  \label{eq:free-V}
  V^{\text{free}} := \{ \delta \us \in H^2(\Gamma_0) \},
\end{equation}
where we denote with $H^k(\Gamma_0)$ the Sobolev space of
three-dimensional vector functions on $\Gamma_0$ with square integrable weak derivatives up to order $k$. 
For the hinged case, we assume that the solid is fixed on the
portion $\partial \Gamma_{0,D} $, but it is free to rotate there, and
the appropriate functional space would be
\begin{equation}
  \label{eq:hinge-V}
  V^{\text{hinged}} := \{ \delta \us \in H^2(\Gamma_0) \text{ s. t. } \delta \us = 0 \text{ on
  }  \partial \Gamma_{0,D} \}.
\end{equation}

If the body is clamped on $\partial \Gamma_{0,D} $, then the correct
functional space is given by
\begin{equation}
  \label{eq:clamp-V}
  V^{\text{clamped} } := \{ \delta \us \in H^2(\Gamma_0) \text{ s. t. } \delta \us = 0, \delta
  \varphi_n  = 0 \text{ on
  }  \partial \Gamma_{0,D} \},
\end{equation}
where $\varphi_n$ describes the normal rotation on the boundary
(rotation around the edge) and it is defined by
$\varphi_n = \nabla_S \X \cdot \vv n_0$, with $\nabla_S$ indicating
the surface gradient and $\vv n_0$ as the outward normal vector on the
boundary.

\section{Algorithmic analysis of the FSI problem}
\label{sec:algor-analys-fsi}

We rewrite system \eqref{eq:final-fsi-system} in operator form, to
make some considerations on possible solution algorithms for the the
final fluid-structure interaction problem.  We indicate with the
functional space $V$ either one of \eqref{eq:free-V},
\eqref{eq:hinge-V} or \eqref{eq:clamp-V}, with $V^*$ its dual space
(the space of all linear operators on $V$) and with
$\duality{\cdot,\cdot}$ the duality product between $V^*$ and $V$,
i.e., 
\begin{equation}
  \label{eq:duality-definition}
  \duality{a,b} := \int_{\Gamma_0} a b \d A, \qquad \forall a \in
  V^*, \forall b \in V.
\end{equation}

With this notation, \eqref{eq:final-fsi-system}  can be rewritten as
\begin{align}
  &\OMm\ddot \us -\OC(\u)\dot \us+\OFi(\us) - \OFe = 0 \qquad \text{ in } V^*.
    \label{eq:operator-motion}
\end{align}

The operators in equation~\eqref{eq:operator-motion} are defined through their
action on arbitrary virtual displacements $\delta \us$:
\begin{align}
  & \duality{\OMm\ddot \us, \delta \us} &  :=  & \int_{\Gamma_0} \rho
\ddot \us \cdot\delta \us \d A& \forall \delta \us \in V \\
  & \duality{\OC(\u)\dot \us, \delta \us} & := & \int_{\Gamma_0} J
 \left( \DNt \dot \us \right)\delta \us \d A& \forall \delta \us \in V \\
  &  \duality{\OFi(\us), \delta\us} & := & \int_{\Gamma_0} \vv n:\delta\boldsymbol{\varepsilon} +
 \vv m:\delta\boldsymbol{\kappa} \d A& \forall \delta
 \us \in V \\
  & \duality{\OFe, \delta\us} & := & \int_{\Gamma_0}  \vv b \cdot
 \delta \us \d A& \forall \delta \us \in V.
    \label{operator-motion}
\end{align}

We observe that the fluid operator $\OC(\us)$ is nonlinear in the
displacement field $\us$, but it is linear in the \emph{velocity}
$\dot \us$, while the elastic operator $\OFi(\us)$ is nonlinear in
$\us$, and, for our choice of elastic constitutive model, it is rate
independent.

The presence of the fluid is felt by the structure solely through the
non-linear operator $\OC(\us)\dot \us$, which acts as a damping term
for the dynamics of the elastic structure. Classical visco-elastic
shells have a very similar structure, where $\OC \dot \us$ is usually
taken to be linear in $\dot \us$, and independent on $\us$.

Given the linearity of the problem in both $\ddot \us$ and $\dot \us$,
a possible solution strategy is to introduce a time discretization
$\mathcal T := \{ t_0, t_1, \dots, t_N = T\}$ and write
$\us_k := \us(t_k)$. At each time step $k$, we can approximate
$\dot \u_k $ and $\ddot \u_k $ as a linear combination of the previous
solution steps $\u_p$ with $p\leq k$, such that the problem reduces to
a nonlinear system in $\us_k$:

\begin{align}
  &\OMm\ddot \us_k -\OC(\u_k)\dot \us_k+\OFi(\us_k) - \OFe_k =:
    \mathcal R(\us_k) = 0 \qquad \text{ in } V^*,
    \label{eq:operator-motion-time-discrete}
\end{align}
whose solution can be formally computed by a Newton iteration method,
i.e., given a guess $\us_k^0 = \us_{k-1}$, we compute $\us_k^{m+1} =
\us_k^{m} +\Delta \us_k^m$ where formally\\
\begin{align}
  &\Delta \us_k^m = -(D_{\us}   \mathcal R(\us^m_k))^{-1} \mathcal
    R(\us_k^m),
    \label{eq:operator-motion-newton}
\end{align}
where $\mathcal R(\us_k^m)$ is the residual at step $m$, and the term
$D_{\us}   \mathcal R(\us^m_k)$ contains the Fr\'echet derivative of the residual
w.r.t. $\us$, evaluated at $\us_k^m$, i.e., 
\begin{equation}
  D_{\us}   \mathcal R(\us^m_k) = c_0 \OMm - c_1 \OC(\u_k^m) +
  D_{\us}\OC(\u_k^m) \dot \us^m_k + D_{\us}\OFi(\u_k^m) ,
  \label{eq:full-nonlinear-Jacobian}
\end{equation}
where $c_0$ and $c_1$  are the linear coefficients of the $\us_k$
term used to approximate $\ddot \us_k$ and $\dot \us_k$. 
The Fr\'echet derivative in
equation~\eqref{eq:full-nonlinear-Jacobian} translates to the Jacobian
of the residual in a finite dimensional setting. If computed directly using
equation~\eqref{eq:full-nonlinear-Jacobian}, such a Jacobian
can be quite complex to approximate, and several simplifications can
be proposed, leading to a Newton-Rapson iteration method in which the
Jacobian of the residual is not exact, but only approximate.

Among these methods, the most commonly used are \emph{segregated}
methods, where the solution of the fluid system is done separately
with respect to the solution of the solid system. These methods are
equivalent to a variation of the following systems of equations
\begin{align}
  &\OMm\ddot \us_k -\OC(\u_{k-1})\dot \us_{k-1}+\OFi(\us_{k\phantom{-1}}) - \OFe_k = 0 \qquad \text{ in } V^*
    \label{eq:operator-motion-time-discrete-leap-frog-solid-dominates} \\
  &\OMm\ddot \us_k -\OC(\u_{k\phantom{-1}})\dot \us_{k\phantom{-1}}+\OFi(\us_{k-1}) - \OFe_k =
  0 \qquad \text{ in } V^*,
  \label{eq:operator-motion-time-discrete-leap-frog-fluid-dominates}
\end{align}
where, in the first case
(equation~\eqref{eq:operator-motion-time-discrete-leap-frog-solid-dominates})
the fluid terms (i.e., $\OC(\u_{k-1})\dot \us_{k-1}$) are computed at
the previous time step, and a full nonlinear solution step is iterated
on the solid part (solid-dominated segregated FSI schemes), while in the
second case
(equation~\eqref{eq:operator-motion-time-discrete-leap-frog-fluid-dominates})
the opposite happens (fluid-dominated segregated FSI schemes).

Due to the nature of the fluid solver, the second family of segregated
solvers is in general difficult to achieve for Boundary Element
Methods, since it requires computing the Jacobian of the 
fluid-structure operator $D_{\us} \OC(\u_k^m)$, whose computational cost is
in the order of $O(n^3)$, where $n$ is the number of degrees of
freedom of the system. 

On the other hand, the structure of the problem suggests naturally a
semi-implicit solution scheme, in which the nonlinearity of the fluid
structure interaction is removed from the system, by evaluating the
fluid-structure operator at the previous time step, but retaining the
evaluation of the velocity field at the current time step, i.e., solving
\begin{align}
  &\OMm\ddot \us_k -\OC(\u_{k-1})\dot \us_{k}+\OFi(\us_k) - \OFe_k = 0 \qquad \text{ in } V^*.
    \label{eq:operator-motion-time-discrete-semi-implicit}
\end{align}

This solution strategy can be further refined by replacing the
computation of $\OC(\u_{k-1})$ with the current nonlinear iterate
$\OC(\u_{k}^m)$, resulting in a Newton-Rapson iteration scheme, in
which the fully implicit nonlinear
system~\eqref{eq:operator-motion-time-discrete} is resolved by
replacing the exact Jacobian in equation~\eqref{eq:full-nonlinear-Jacobian}
by an approximation in which $D_{\us} \OC(\us_{k}^{m})\dot \us_k$ is
neglected, i.e.
\begin{equation}
  D_{\us}   \mathcal R(\us^m_k) \sim c_0 \OMm - c_1 \OC(\u_k^m)
  + D_{\us}\OFi(\u_k^m).
  \label{eq:inexact-Jacobian}
\end{equation}

The details of the full discrete scheme are given in
section~\ref{sec:KLtheory}, where a generalized $\alpha$-scheme is
coupled with the fully implicit solver for the coupled system, with
inexact Jacobian given by the discrete version
of~\eqref{eq:inexact-Jacobian}.

\section{Isogeometric boundary integral representation of the hydrodynamic equations}
\label{sec:hydr-inter}
The fluid part of the equations of motion takes the form
\begin{subequations}
 \label{eq:main-fluid-system}
  \begin{alignat}{2}
    & - \eta \Delta \vf + \nabla p = -\nabla \cdot
    \vv{\cauchyf} = 0 \qquad && \mbox{in } \Of(t) \label{eq:momentum-cons} \\
    & \nabla\cdot \vf =0 && \mbox{in } \Of(t) \label{eq:mass-cons}\\
    & \vf = \vs_g && \mbox{on } \Gamma(t) \label{eq:dirichlet-bc}
  \end{alignat}
\end{subequations}
where $\vf$ and $p$ are the velocity and hydrodynamic pressure fields
in the domain $\Of(t):=\Re^3\setminus\Gamma(t)$, $\eta$ is the
viscosity of the fluid, $\vs_g$ is the (given) velocity of the middle
surface of the shell and $\vv \cauchyf$ is the Cauchy stress tensor
for an incompressible Newtonian fluid:
\begin{equation}
  \label{eq:cauchy-stress}
  \vv \cauchyf := - p \vv I +\eta (\nabla \vf+ (\nabla\vf) ^T).
\end{equation}

Equations~\eqref{eq:momentum-cons} and~\eqref{eq:mass-cons} describe
the conservation of linear momentum and volume in the Stokes fluid,
while \eqref{eq:dirichlet-bc} is a Dirichlet boundary condition. The
pressure $p$ can be regarded as the Lagrange multiplier associated
with the conservation of volume~\eqref{eq:mass-cons}, and it is
uniquely determined by $\vf$ up to an additive constant.

Following \cite{Pozrikidis-1992-a} or \cite{Steinbach2008},
we can write a boundary integral representation of the solution $\vf$
and $\cauchyf$ of system~\eqref{eq:main-fluid-system} using the
free-space Green's functions ${\vv {\mathcal
    S}}$ %, ${\vv {\mathcal P}}$
and ${\vv {\mathcal T}}$:
\begin{equation}
  \label{eq:free-space-green-functions}
  \begin{aligned}
    {\mathcal S}_{ab}(\vv r) & =   \frac{1}{8\pi \eta}\left(\frac{r_a r_b}{|\vv r|^3} +
      \frac{\delta_{ab}}{|\vv r|} \right)
  \\
%    \vv {\mathcal P}_i(\vv r) & = (d-1)\frac{r_i}{|\vv r|^d}\\
    {\mathcal T}_{abc}(\vv r) & = -  \frac{3}{4\pi }  \frac{r_ar_br_c}{|\vv r|^{5}},
  \end{aligned}
\end{equation}
where $\vv r$ is a shorthand notation for $(\vv x - \vv y)$, and 
\nomenclature{$\mathcal{S}$}{Free space Green function for the velocity}
\nomenclature{$\mathcal{T}$}{Free space Green function for the cauchy tensor}
$\mathcal{S} \vv b$ and $\mathcal{T} \vv b$ are the velocity and
stress fields in free space associated to a Dirac force with intensity
$\vv b$ centered in $\vv y$.

Given an arbitrary control volume $S$ such that $S\cap \Gamma(t) =
\emptyset$, it is possible to express the velocity $\v$ at arbitrary
points $\vv x \in S$, as \\
\begin{equation}
  \label{eq:bie-stokes-on-S}
  v_a(\vv x) +\int_{\partial S} {\mathcal T}_{abc}(\x-\y) \nu_b(\vv y)v_c(\vv y)
  \d\Gamma_y = \int_{\partial S} {\mathcal S}_{ab}(\x-\y)
  \cauchyf_{bc}(\vv y)\nu_c(\vv y) \d \Gamma_y.
\end{equation}

If we select $S = \Re^3\setminus S_\varepsilon$, defined as
\begin{equation}
  \label{eq:S-epsilon}
  S_\varepsilon := \{ \x(\vv s) + h \vv \nu(\vv s), \qquad \vv s \in B^2,
  h \in (-\varepsilon/2, \varepsilon/2) \},
\end{equation}
where $\vv \nu$ here and above is $\vv g_3$, i.e., the normal vector to
the middle surface, then taking the limit for $\varepsilon \to 0$,
the domain $S$ would coincide with $\Gamma(t)$, and
equation~\eqref{eq:bie-stokes-on-S} would collapse to (using compact notation)
\begin{equation}
  \label{eq:bie-stokes-on-Gamma}
  \vf(\x) +\int_{\Gamma(t)} \vv {\mathcal T}(\x-\y) \vv \nu(\y)
  \jump{\vf(\y)} \d\Gamma_y = \int_{\Gamma(t)} \vv {\mathcal S}(\x-\y) 
  \jump{\cauchyf(\y)} \vv \nu(\y) \d \Gamma_y.
\end{equation}

Such a limit may not be well posed if we considered
equation~\eqref{eq:bie-stokes-on-Gamma} \emph{as is}, since we went
from a closed surface with no boundaries to a surface with boundaries,
and the boundary element method may present singularities on the
curves representing the boundary of the two-dimensional middle
surface. However, since we impose a no-slip boundary condition on
$\Gamma(t)$, the velocity of the fluid on the middle surface coincides
with the velocity of the solid on both sides, making the second term
on the left hand side in equation~\eqref{eq:bie-stokes-on-Gamma}
identically zero. This term would be responsible for singularities on
the one dimensional boundary of the middle surface, which are not
there if one considers only the single layer. If we take the trace of
this equation on $\Gamma(t)$, the integral on the right hand side
becomes weakly singular but integrable, and we obtain a boundary
integral equation on the surface $\Gamma(t)$, which can be used to
explicitly compute the force per unit area applied by the fluid on the
solid $\vv f := \jump{\cauchyf}\vv \nu$, solving the following
integral equation of the first kind: \begin{equation}
  \label{eq:boundary-integral-equation-single-layer}
  \vf(\vv x) = \int_{\Gamma(t)} \vv {\mathcal S}(\vv x-\vv y)
  \vv f(\vv y)\d \Gamma_y \qquad \forall \vv x
  \text{ on } \Gamma(t).
\end{equation}

Such a boundary integral equation generates a fluid velocity field $\vf$
which is globally in $H^1(\Re^3)$ for any surface traction in
$H^{-1/2}(\Gamma)$. An effective way to numerically solve this
boundary integral equation is given by the boundary element method, in
which $\vf$ and $\vv f$ are sought for in a finite dimensional space
defined on $\Gamma(t)$, and the Dirichlet to Neumann map becomes an
invertible matrix.

Here we exploit the isogeometric NURBS spaces defined in
Section~\ref{sec:bound-integr-appr} to define the finite dimensional
spaces, as well as the discrete versions of the boundary integral
equations~\eqref{eq:boundary-integral-equation-single-layer}.

We \emph{collocate} the boundary integral equation at $n$ distinct
\emph{collocation points} $\{\vv x^{\vv i}\}_{\vv i \in \mathcal J}$,
and we restrict both $\vf(\vv x)$ and $\vv f(\vv x)$ to live in the
finite dimensional space $V_h^3(\Gamma)$:
\begin{align}
  \label{eq:collocated-bie}
  \vf^{\vv J} \vv \Phi^{\vv J}(\vv x^{\vv i}) = 
  \int_{\Gamma(t)} \vv {\mathcal S}(\vv x^{\vv I}-\vv y)
  \vv f^{\vv J} \vv \Phi ^{\vv J}(\vv y)\d \Gamma_y \qquad \vv i \in
  \mathcal{J}, \vv J \in \mathcal{J}^3.
\end{align}

For each collocation point $\vv x^{\vv I}$,
Equations~\eqref{eq:collocated-bie} are systems of $3$ equations in
$6n$ unknowns (the $3n$ coefficients of $\vf$ and the $3n$
coefficients of $\vv f$), which can be compactly rewritten as
\begin{alignat}{2}
  \label{eq:discrete-bie}
  \vv M_c \vf &= \vv D_c \f
\end{alignat}
\nomenclature{$\vv M_c$}{Collocation matrix}
\nomenclature{$\vv D_c$}{Single layer collocation matrix}
where the (square) matrices $\vv M_c$ and $\vv D_c$ are given by
\begin{alignat}{2}
  \label{eq:mass-matrix}
  \vv M_c^{(3\vv i +a)\, (\vv J)} &:= \delta_{ab}  \Phi^{\vv J}_b
  (\vv x^{\vv i})\\
  \label{eq:D-matrix}
  \vv D_c^{(3\vv i +a)\, (\vv J)} &:= \int_{\Gamma(t)}  {\mathcal
    S}_{ab}(\vv  x^{\vv I}-\vv y)  \Phi^{\vv J}_b(\vv y)
  \d \Gamma_y .
\end{alignat}

A common approach for the choice of the collocation points is given by
the Greville absciss\ae ~(see, for example,~\cite{Greville-1964-a},
or~\cite{AuricchioDaVeigaHughes2010aa}), which are defined as
\begin{equation}
  \label{eq:greville}
  \vv x^{\vv i} := \vv x(\vv s^{\vv i}), \qquad \vv
  s^{i_m} := \frac{\sum_{j=1}^{p} k^m_{i_m+j}}{p},
\end{equation}
where $k^m_{i_m+j}$ are the knots of the knot vector
$\Theta_i$. Care should be taken in order to avoid collapsing
collocation points, which would result in singular matrices.

A discrete version of the Dirichlet to Neumann operator $\DNt$ is
then given by the \emph{damping matrix} $\vv C(\u)$:
\begin{equation}
  \label{eq:discrete-DN}
  \vv C(\u) :=  \vv M_u \vv D_c^{-1} \vv M_c, 
\end{equation}
\nomenclature{$\vv M_u$}{Pseudo mass matrix}
where the matrix $\vv M_u$ is a pseudo mass matrix, defined as
\begin{equation}
  \label{eq:mass-matrix-deformed}
  \vv M^{\vv I\vv J}_u := \int_{\Gamma(t)} \vv \Phi^{\vv I}
  \cdot \vv \Phi^{\vv J} \d \Gamma =   \int_{\Gamma^0} \vv \Phi^{\vv I}
  \cdot \vv \Phi^{\vv J} J_u \d \Gamma.
\end{equation}

In general it is not necessary to explicitly assemble the matrix $\vv C(\u)$, as long as we can \emph{compute its action} on arbitrary vectors.  Such action requires the solution of the flow problem around $\Gamma(t)$, obtained through the inversion of the (dense) operator $\vv D_c$, pre and post multiplied by two (sparse) matrix multiplications.

% ========================================================================================
% ========================================================================================
% ========================================================================================
% ========================================================================================
\section{Isogeometric Galerkin approach for nonlinear shell dynamics} \label{sec:KLtheory}
We solve the structural dynamics problem by isogeometric Galerkin discretizations, and rewrite equation \eqref{eq:operator-motion-time-discrete} in the discrete form:
\begin{align}
 \Gm(\us_k) = \vv M_s\udd_k-\vv C(\us_k)\ud_k+\Fi(\us_k)-\Fe=\mathbf{0}. \label{motion2}
\end{align}
In equation \eqref{motion2}, $\us_k, \ud_k, \udd_k$ indicate the vectors of nodal displacements, velocities, and accelerations at a time step $k$, $\Fi$ and $\Fe$ are the vectors of internal and external nodal forces, respectively, $\Gm$ is the residual vector, $\vv C$ is the viscous damping matrix, representing the action of the surrounding fluid and obtained according to equation \eqref{eq:discrete-DN}, and $\vv M_s$ is the structural mass matrix, obtained in the reference configuration as:
\begin{align}
  \vv M^{\vv I\vv J}_s := \rho h \int_{\Gamma_0} \vv \Phi^{\vv I} \cdot \vv \Phi^{\vv J} \d \Gamma \label{M} 
\end{align}
\nomenclature{$\vv M_s$}{Structural mass matrix}
with $\rho$ as the density and $h$ as the shell thickness. \par
%equation \eqref{motion2} is linearized and solved for the displacements, using the approximated Jacobian according to equation \eqref{eq:inexact-Jacobian}. 
As time integration scheme, we use a generalized $\alpha$-method \cite{chung_time_1993,CottrellHughesBazilevs-2009-a}, where the displacements, velocities, and accelerations are interpolated at time instants between two discrete time steps $t_{k-1}$ and $t_{k}$ as follows: 
\begin{align}
  \ua &= \alpha_{f}\us_{k}+(1-\alpha_{f})\us_{k-1}  \text{ ,}\label{ualf} \\
  \uda &= \alpha_{f}\ud_{k}+(1-\alpha_{f})\ud_{k-1}  \text{ ,}\label{valf} \\
  \udda &= \alpha_{m}\udd_{k}+(1-\alpha_{m})\udd_{k-1}  \text{ ,}\label{aalf} 
\end{align}
where the velocity and displacement at time step $t_{k}$ are defined by a Newmark update:
\begin{align}
  \us_{k}&=\us_{k-1}+\Delta t\ud_{k-1}+\frac{1}{2}(\Delta t)^2 \left((1-2\beta)\udd_{k-1}+2\beta\udd_{k}\right)\text{ ,} \label{uofa} \\
  \ud_{k}&=\ud_{k-1}+\Delta t \left((1-\gamma)\udd_{k-1}+\gamma\udd_{k}\right)\text{ ,} \label{vofa}
\end{align}
with $\beta$ and $\gamma$ as the Newmark parameters and $\Delta t=t_{k}-t_{k-1}$ as the time step size. Solving for the displacements $\us_{k}$ first, the Newmark updates of velocities and accelerations are obtained as:
\begin{align}
  \ud_{k}&=\frac{\gamma}{\beta\Delta t}(\us_{k}-\us_{k-1}) + \left(1-\frac{\gamma}{\beta}\right)\ud_{k-1} + \left(1-\frac{\gamma}{2\beta}\right) \Delta t\udd_{k-1}\text{ ,} \label{vofu} \\
\udd_{k}&=\frac{1}{\beta(\Delta t)^2}(\us_{k}-\us_{k-1}) - \frac{1}{\beta\Delta t}\ud_{k-1} - \left(\frac{1}{2\beta}-1\right)\udd_{k-1} \text{ .}\label{aofu} 
\end{align}
The $\alpha$ and Newmark parameters are determined by the numerical dissipation parameter $\rho_\infty \in [0,1]$ as follows:
\begin{align}
  \alpha_{m}=\frac{2-\rho_\infty}{1+\rho_\infty} \text{ ,}\quad  \alpha_{f}=\frac{1}{1+\rho_\infty}\text{ ,}
  \quad   \beta=\frac{(1-\alpha_{f}+\alpha_{m})^2}{4}   \text{ ,} \quad   \gamma=\frac{1}{2}-\alpha_{f}+\alpha_{m}\text{ ,}
\end{align}
where $\rho_\infty=0.5$ is adopted in this paper.

% External forces, considered to be independent of $\us$, are evaluated as:
%\begin{align}
%  \Fe &= \alpha_{f}\Fe_{k}+(1-\alpha_{f})\Fe_n \text{ .} \label{Fealf2} 
%\end{align}
With the interpolated variables \eqref{ualf}-\eqref{aalf}, equation \eqref{motion2} is linearized and solved for the displacements, using the approximated Jacobian according to equation \eqref{eq:inexact-Jacobian}, which yields the following system of equations:
\begin{align}
  \left(\alpha_{m}\frac{1}{\beta(\Delta t)^2}\vv M_s-\alpha_{f}\frac{\gamma}{\beta\Delta t}\vv C(\ua^m)+\alpha_{f}\K(\ua^m)\right) &\Delta\us_{k}^m =& \notag \\
   -\vv M_s\udda^m+\vv C(\ua^m)\uda^m&-\Fi(\ua^m)+\Fe_\alpha \text{ ,} \label{genalf_fully}
\end{align}
with $\K$ being the structural stiffness matrix. Equation \eqref{genalf_fully} represents the fully implicit nonlinear system corresponding to \eqref{eq:inexact-Jacobian}. As outlined in Section \ref{sec:algor-analys-fsi}, we further consider a semi-implicit \eqref{eq:operator-motion-time-discrete-semi-implicit} and a segregated \eqref{eq:operator-motion-time-discrete-leap-frog-solid-dominates} approach for the fluid-structure coupling. For the semi-implicit approach, the damping matrix $\vv C(\ua^m)$ is approximated by $\vv C(\us_{k-1})$
\begin{align}
  \left(\alpha_{m}\frac{1}{\beta(\Delta t)^2}\vv M_s-\alpha_{f}\frac{\gamma}{\beta\Delta t}\vv C(\us_{k-1})+\alpha_{f}\K(\ua^m)\right) &\Delta\us_{k}^m =& \notag \\
   -\vv M_s\udda^m+\vv C(\us_{k-1})\uda^m&-\Fi(\ua^m)+\Fe_\alpha \text{ ,} \label{genalf_semi}
\end{align}
with the effect that the fluid equations have to be assembled and solved only once per time step. For the segregated approach the whole damping term $\vv C(\us)\ud$ is considered constant during one time step. In this case, the contribution from the fluid can be considered as an additional external force $\Fe_{\alpha}^{fsi}=\vv C(\us_{k-1})\ud_{k-1}$, and the damping term on the left hand side of \eqref{genalf_fully} vanishes:
\begin{align}
  \left(\alpha_{m}\frac{1}{\beta(\Delta t)^2}\vv M_s + \alpha_{f}\K(\ua^m)\right) \Delta\us_{k}^m =
   -\vv M_s\udda^m+\Fe_{\alpha}^{fsi}-\Fi(\ua^m)+\Fe_\alpha \text{ .} \label{genalf_semi_2}
\end{align}
% -------- stiffness matrix and internal forces -------------------
The internal force vector $\Fi(\us)$ and stiffness matrix $\K(\us)$ are obtained by linearization of the static terms of internal virtual work of the shell model \eqref{W_intKL} with respect to discrete displacement variables $u_{\vv I}$ and $u_{\vv J}$:
\begin{equation} \label{eq:Fint}
 \vv P^{\vv I} =-\int_A \left(
     \vv n:\frac{\partial\boldsymbol{\varepsilon}}{\partial u_{\vv I}}+
     \vv m:\frac{\partial\boldsymbol{\kappa}}{\partial u_{\vv I}}
     \right)\mathrm dA
\end{equation}
\begin{equation} \label{eq:Kint}
  \vv K^{\vv I\vv J} = \int_A \left(
      \frac{\partial\vv n}{\partial u_{\vv J}}:\frac{\partial\boldsymbol{\varepsilon}}{\partial u_{\vv I}} 
      +\vv n:\frac{\partial^2\boldsymbol{\varepsilon}}{\partial u_{\vv I} \partial u_{\vv J}} 
      +\frac{\partial\vv m}{\partial u_{\vv J}}:\frac{\partial\boldsymbol{\kappa}}{\partial u_{\vv I}} 
      +\vv m:\frac{\partial^2\boldsymbol{\kappa}}{\partial u_{\vv I} \partial u_{\vv J}}
      \right)\mathrm dA
\end{equation}
with
\begin{align} 
  \frac{\partial\vv n}{\partial u_{\vv J}} &=  h\; \hat{\mathbb{C}}:\frac{\partial\boldsymbol{\varepsilon}}{\partial u_{\vv J}} \label{dn^ab} \\
  \frac{\partial\vv m}{\partial u_{\vv J}} &= \frac{h^3}{12}\; \hat{\mathbb{C}}:\frac{\partial\boldsymbol{\kappa}}{\partial u_{\vv J}} \label{dm^ab} 
\end{align}
For the discrete model, $C^1$-continuity of the basis functions is required since second derivatives appear in the definition of the curvatures \eqref{b_ab}. NURBS-based isogeometric discretizations provide the necessary continuity and allow a straightforward implementation of this formulation. The control point displacements are identified as the displacement variables $u_{\vv I}$. The detailed linearization of the strain variables $\boldsymbol{\varepsilon}$ and $\boldsymbol{\kappa}$ with respect to $u_{\vv I}$ and $u_{\vv J}$ is given in \ref{App_lin}. 

% ======================================================================================================

% ======================================================================================================
% ======================================================================================================
% ======================================================================================================
% ======================================================================================================
\section{Numerical tests}
In this section, we apply the presented methods to different numerical tests. First, we consider the free vibration of a beam immersed in a fluid 
and use it to compare the different coupling strategies. Furthermore, we consider a structure which is deformed by 
 externally applied loads and the damping effect of the surrounding fluid and, finally, we consider a free fall problem. 
\label{sec:numtests}
\subsection{Vibration of a cantilever inside a viscous fluid}
\begin{figure}\centering
  \subfigure[$\eta=10\, Pa\cdot s$]{\includegraphics[width=0.49\textwidth]{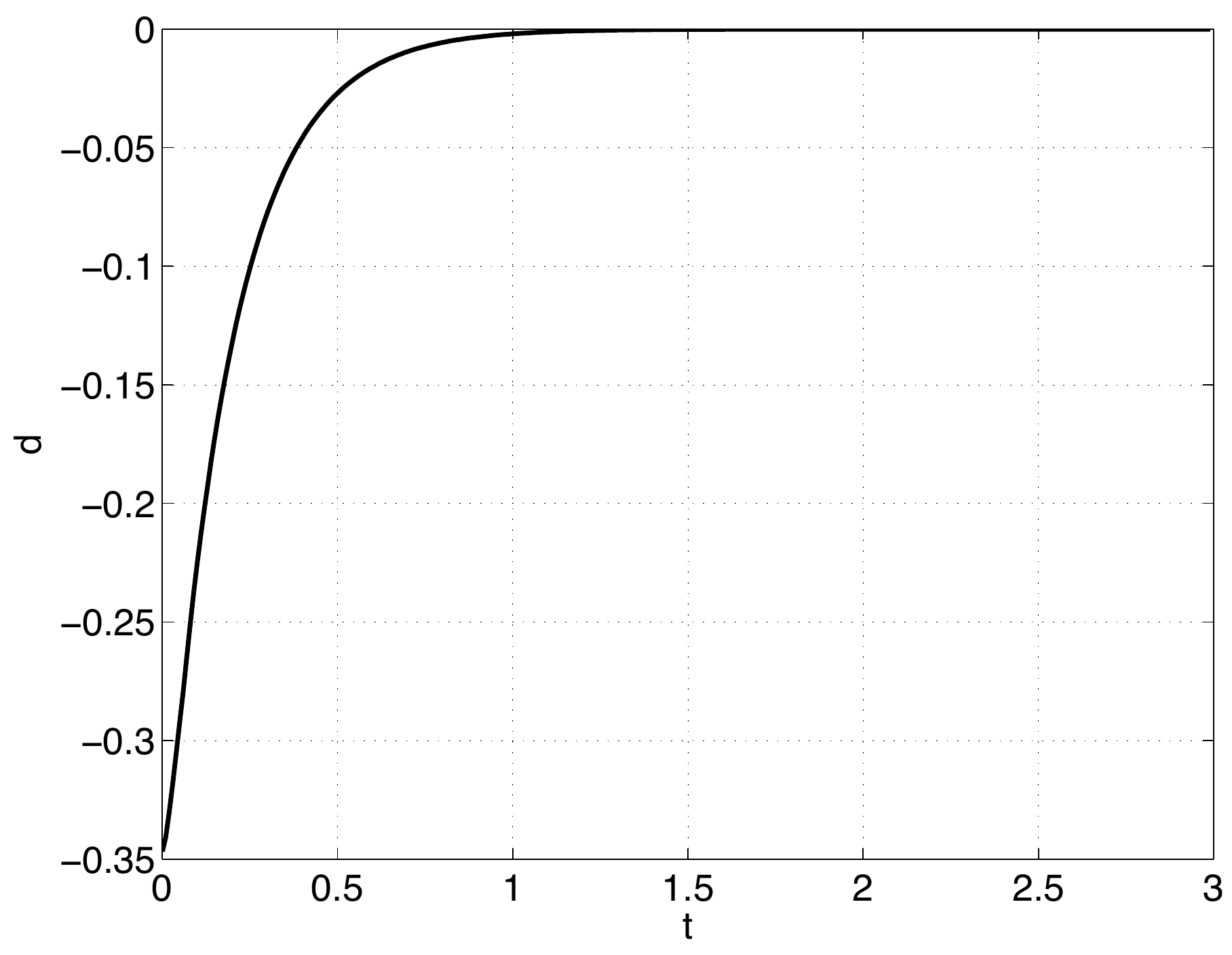}}
  \subfigure[$\eta=1\, Pa\cdot s$]{\includegraphics[width=0.49\textwidth]{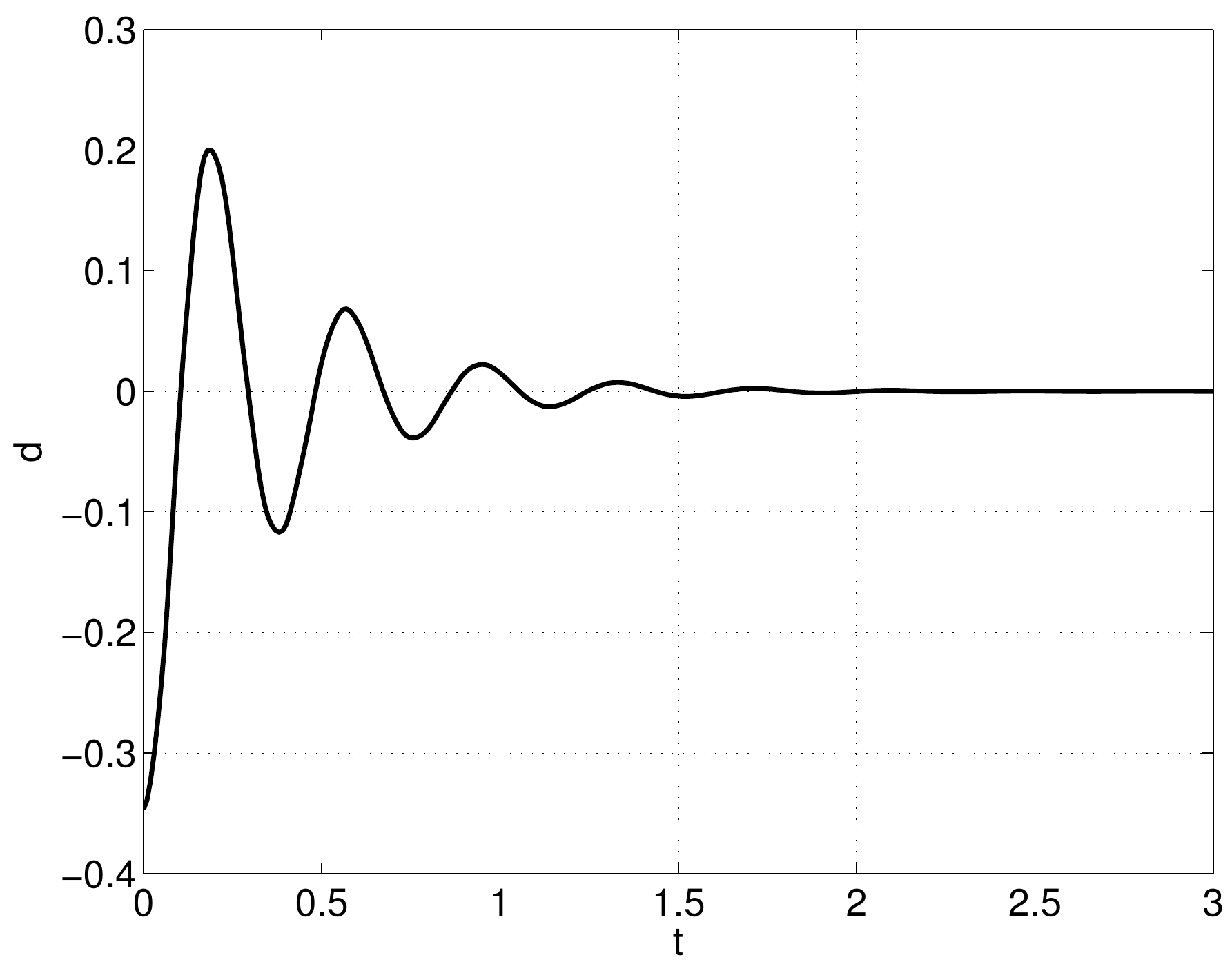}}
  \subfigure[$\eta=0.1\, Pa\cdot s$]{\includegraphics[width=0.49\textwidth]{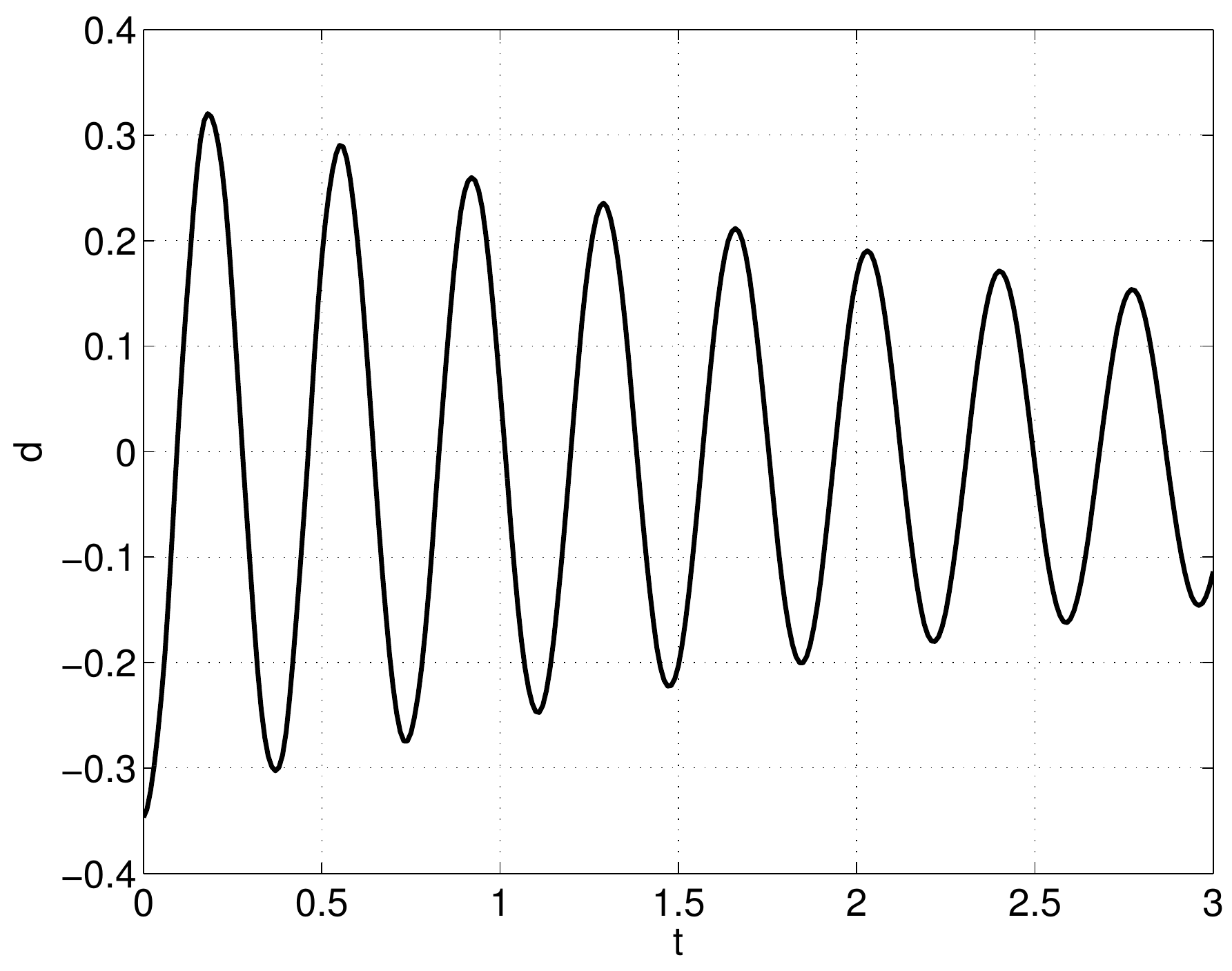}}
  \subfigure[$\eta=0.001\, Pa\cdot s$]{\includegraphics[width=0.49\textwidth]{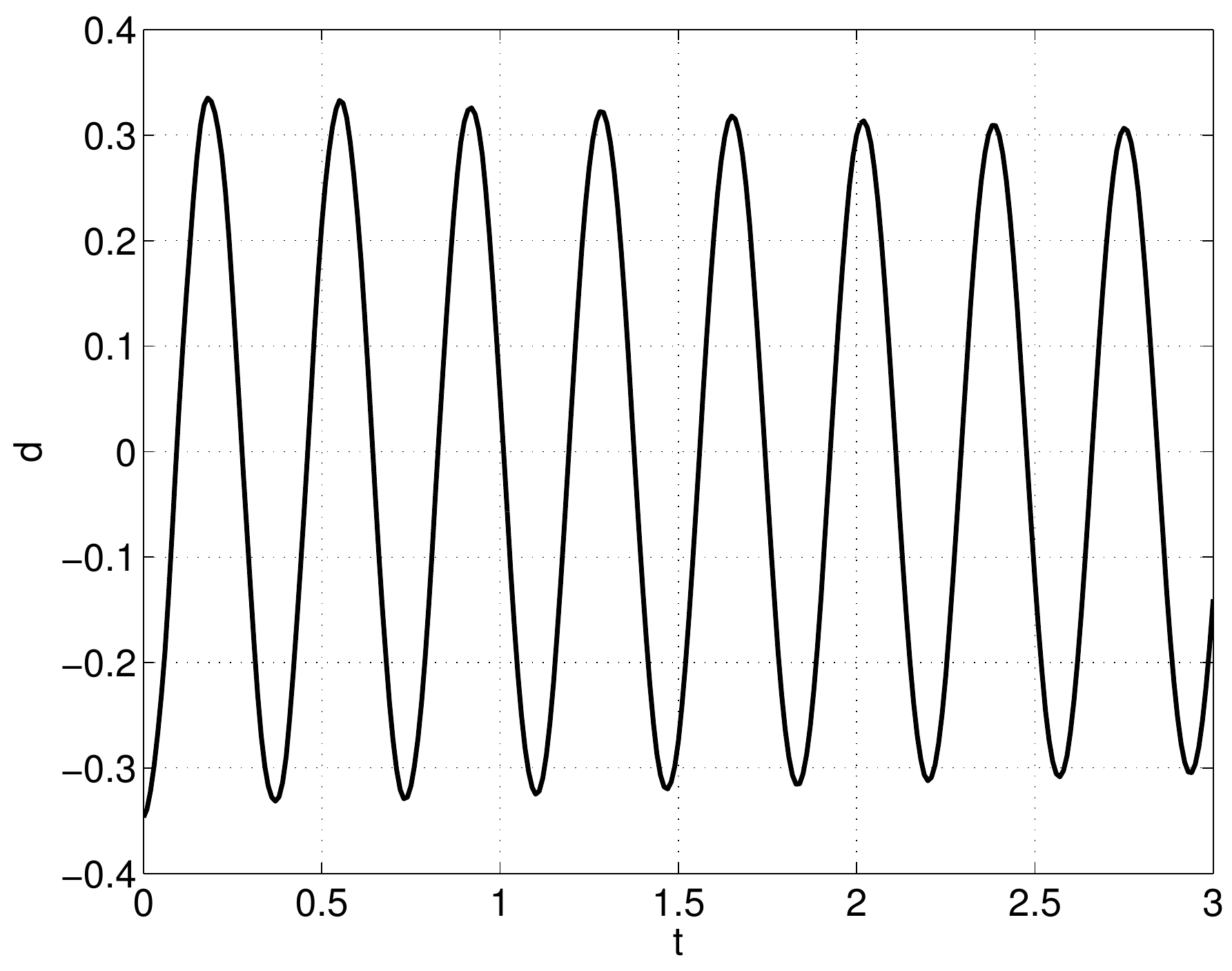}}
  \caption{Vibration of a cantilever immersed in a viscous fluid. Tip displacement for different viscosities and $\Delta t=0.01s$, obtained using the semi-implicit approach.} 
  \label{leaflet1}
\end{figure}

\begin{figure}\centering
  \includegraphics[width=0.49\textwidth]{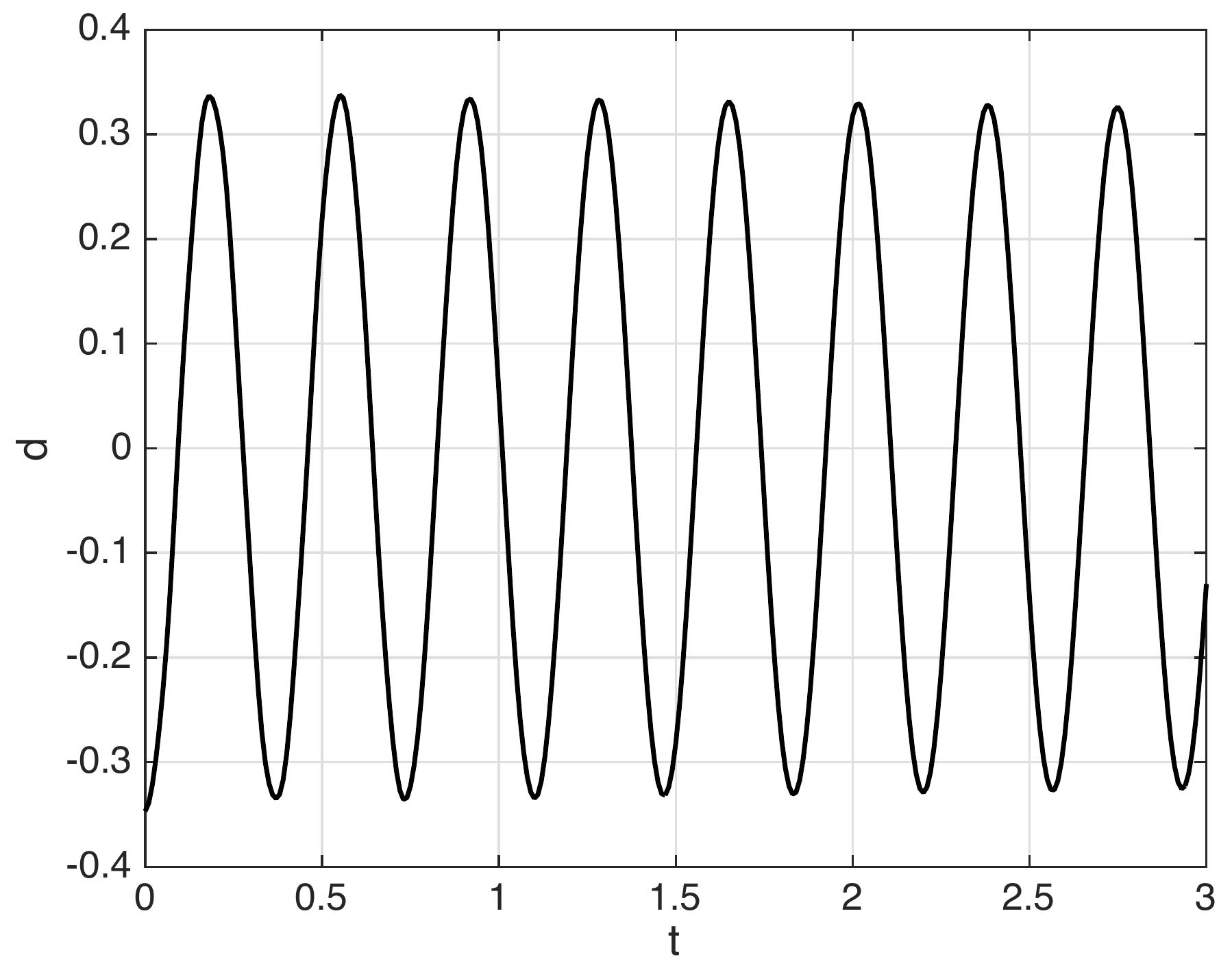}  \caption{Undamped vibration of a cantilever, obtained by a pure structural dynamics computation.} 
  \label{leaflet_free}
\end{figure}

We consider a cantilever plate surrounded by a viscous fluid. The plate dimensions are $1m\times0.1m\times1mm$ ($length\times width\times thickness$). The material parameters are $E=210.1\cdot 10^{10} Pa, \nu=0.3, \rho=7850 kg/m^3$. The plate is clamped at the left edge and initially deformed corresponding to a static load of $225 N/m$ at the right edge (tip). At time $t=0$, the load is removed and the vibrations of the plate are observed by plotting the tip displacement. In Figure \ref{leaflet1}, we plot the tip displacement for different viscosities $\eta=\{10,1,10^{-1},10^{-3}\} \,Pa\cdot s$, using a time step of $\Delta t=0.01s$. 
%These results have been obtained by the semi-implicit approach \eqref{eq:operator-motion-time-discrete-semi-implicit}. 
%As can be seen, good results without artificial oscillations are obtained for all cases. 
% New:
%It can be seen that for the very viscous case with $\eta=10$, the system is over-damped and no oscillations occur, while with decreasing $\eta$, the solution converges to the free vibration solution, shown in Figure \ref{leaflet_free}. This solution was obtained by a purely structural dynamics analysis without damping. The frequency obtained in this case is $f_{num}=2.73 Hz$. To the authors knowledge, there is no analytical solution available for the free vibration of a plate with nonlinear kinematics. However, as an approximate reference, we can compare with the analytical solution from linear beam theory. In this case, the first eigenfrequency of a cantilever beam is given as $f_{ref}=\cfrac{1.875^2}{2\pi L^2}\sqrt{\cfrac{E h^2}{12\rho}}=2.64 Hz$, with $L$ as the beam's length, which is in good agreement with the numerical result, considering the differences in the underlying kinematics. \par
The typical patterns of damped vibration can be observed, with an over-damping in the very viscous case ($\eta=10$), and damped oscillations for the other cases which converge towards the undamped vibration solution for very small viscosities. For comparison, we also perform a purely structural dynamics simulation with no damping, see Figure \ref{leaflet_free}. We measure the natural frequencies (averaged over the first 7 periods) obtained in these simulations and report them in Table \ref{tab:frequencies}, which shows that the frequencies of the oscillations in the fluid-structure interaction model are converging to that of the undamped vibration case as the viscosity $\eta$ converges to zero. We can also observe that, as to be expected, the frequency diminishes as the dissipation increases. 
Furthermore, we compute the analytical solution of an undamped vibrating cantilever beam, which is given by $f_{ref}=\cfrac{1.875^2}{2\pi L^2}\sqrt{\cfrac{E h^2}{12\rho}}$, where $L$ is the beam length. It should be noted that the analytical solution is based on linear beam theory, while the numerical solutions are obtained with a nonlinear shell formulation.
%, and, when $\eta \neq 0$, they include the coupling with the surrounding fluid. 
Since, to the best of our knowledge, there is no analytical solution available for the free vibration of a shell with nonlinear kinematics, we use this beam solution as an approximate reference solution showing that our results are in a physically sound range. The analytical frequency obtained from linear beam theory is $f_{ref}=2.6418 Hz$, showing good agreement with the numerical results, considering the differences in the underlying models. 
\begin{table}[!h]
  \centering
  \begin{tabular}{rl||c|c|c|c||c||}
     $\eta$&$[Pa\cdot s]$: & $1$ & $10^{-1}$ & $10^{-3}$ & Undamped Vibration \\ \hline
     $f$ &$[1/s]$:& 2.6254 & 2.7055 & 2.7254 & 2.7284 
  \end{tabular}
  \caption{Vibration frequencies (averaged over the first 7 periods) for different values of fluid viscosity and for undamped vibration.}
  \label{tab:frequencies}
\end{table}

 \par
\begin{figure}\centering
  \subfigure[fully impl., $\Delta t=0.01s$]{\includegraphics[width=0.32\textwidth]{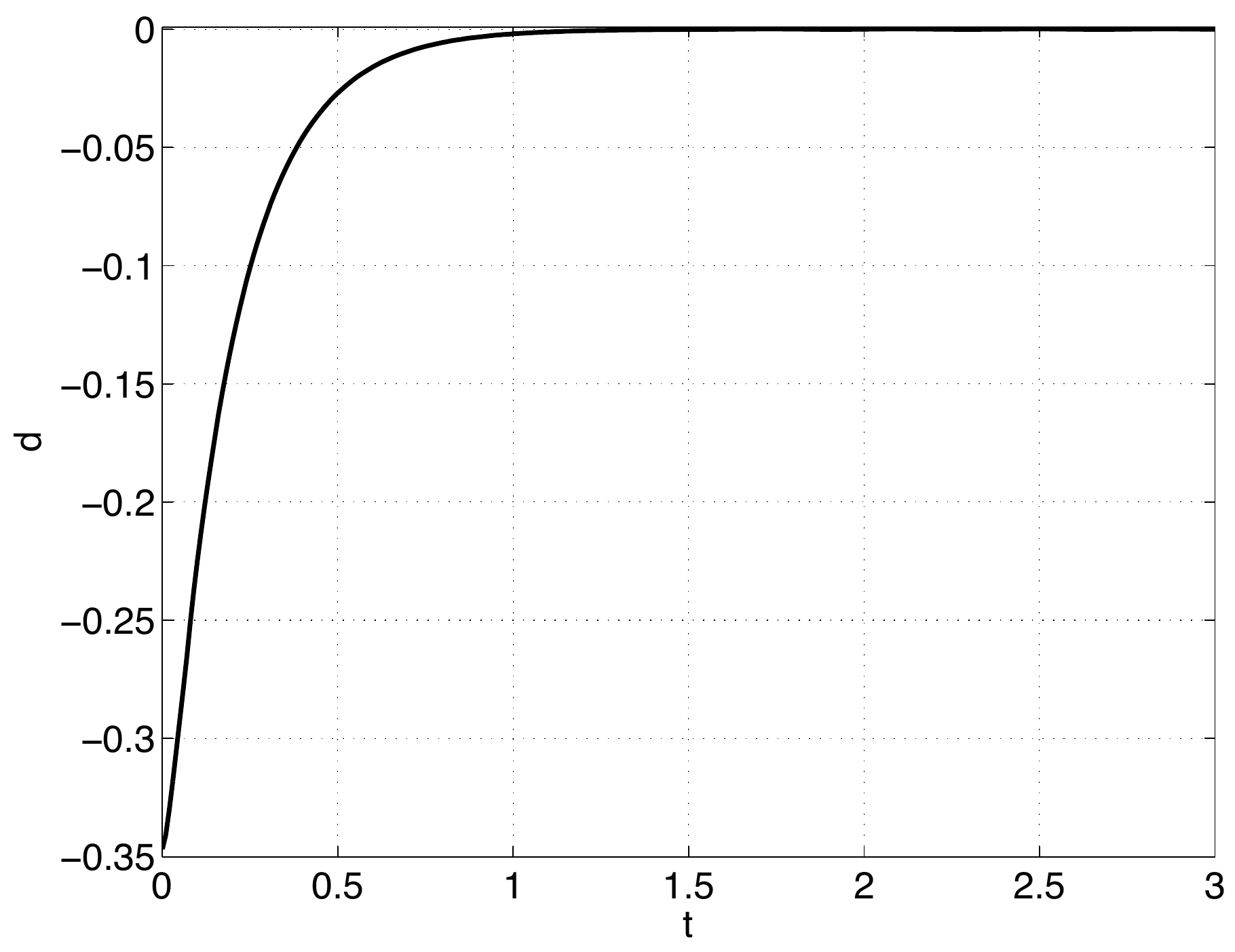}}
  \subfigure[semi-impl., $\Delta t=0.01s$]{\includegraphics[width=0.32\textwidth]{figures/leaflet_new_eta1ep1_fsi1_damp1_dt100_d.pdf}}
  \subfigure[segregated, $\Delta t=0.01s$]{\includegraphics[width=0.32\textwidth]{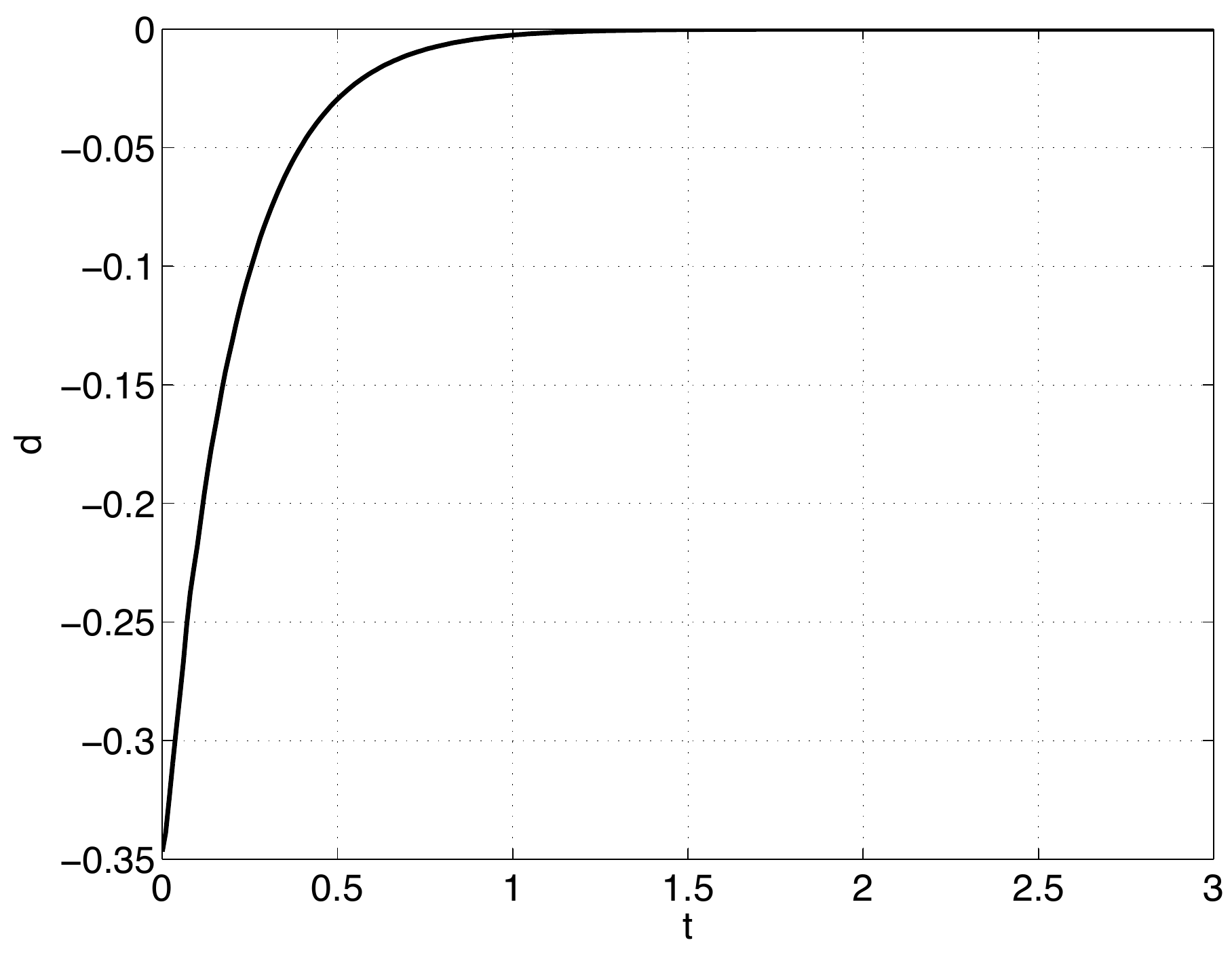}}
  \subfigure[fully impl., $\Delta t=0.05s$]{\includegraphics[width=0.32\textwidth]{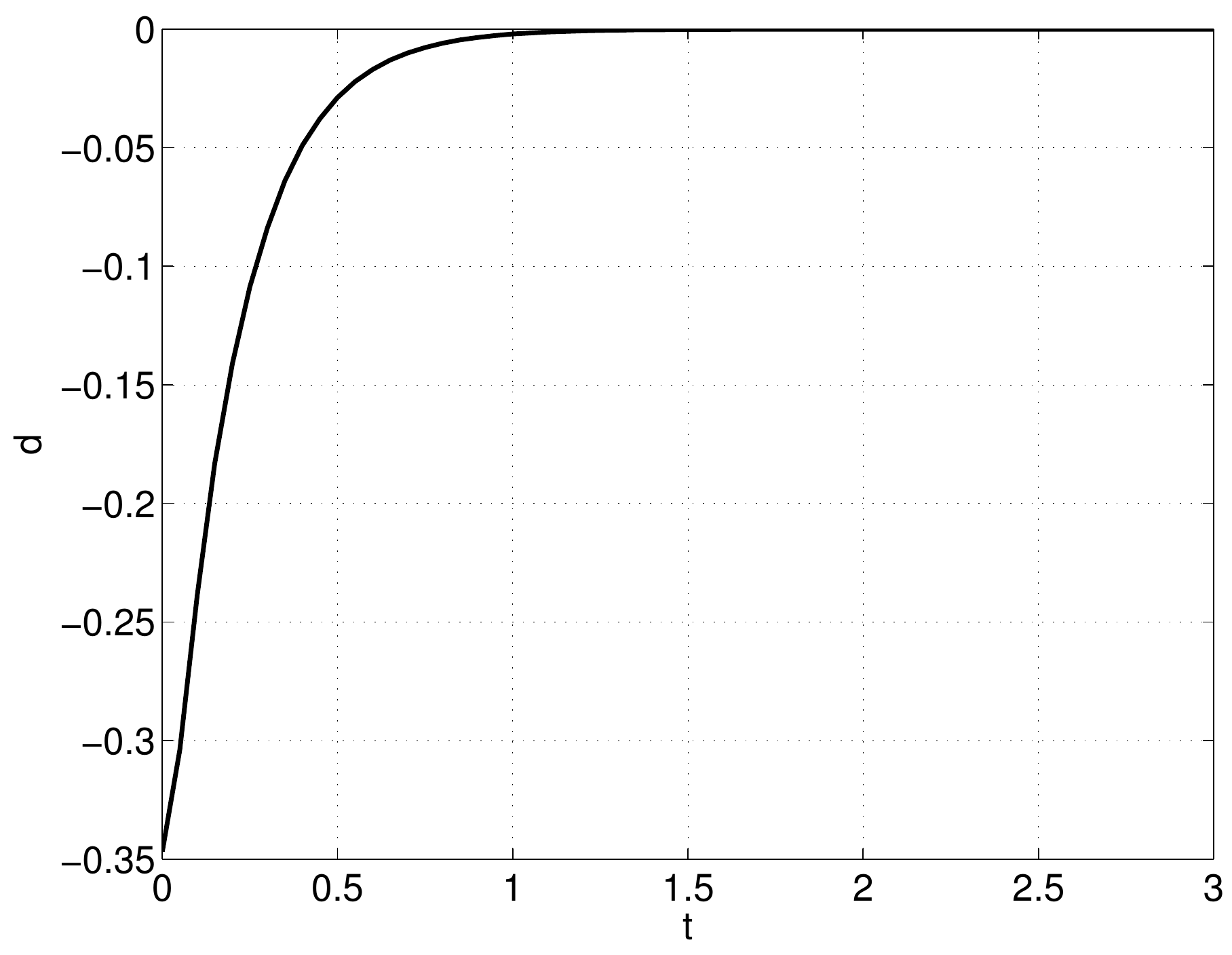}}
  \subfigure[semi-impl., $\Delta t=0.05s$]{\includegraphics[width=0.32\textwidth]{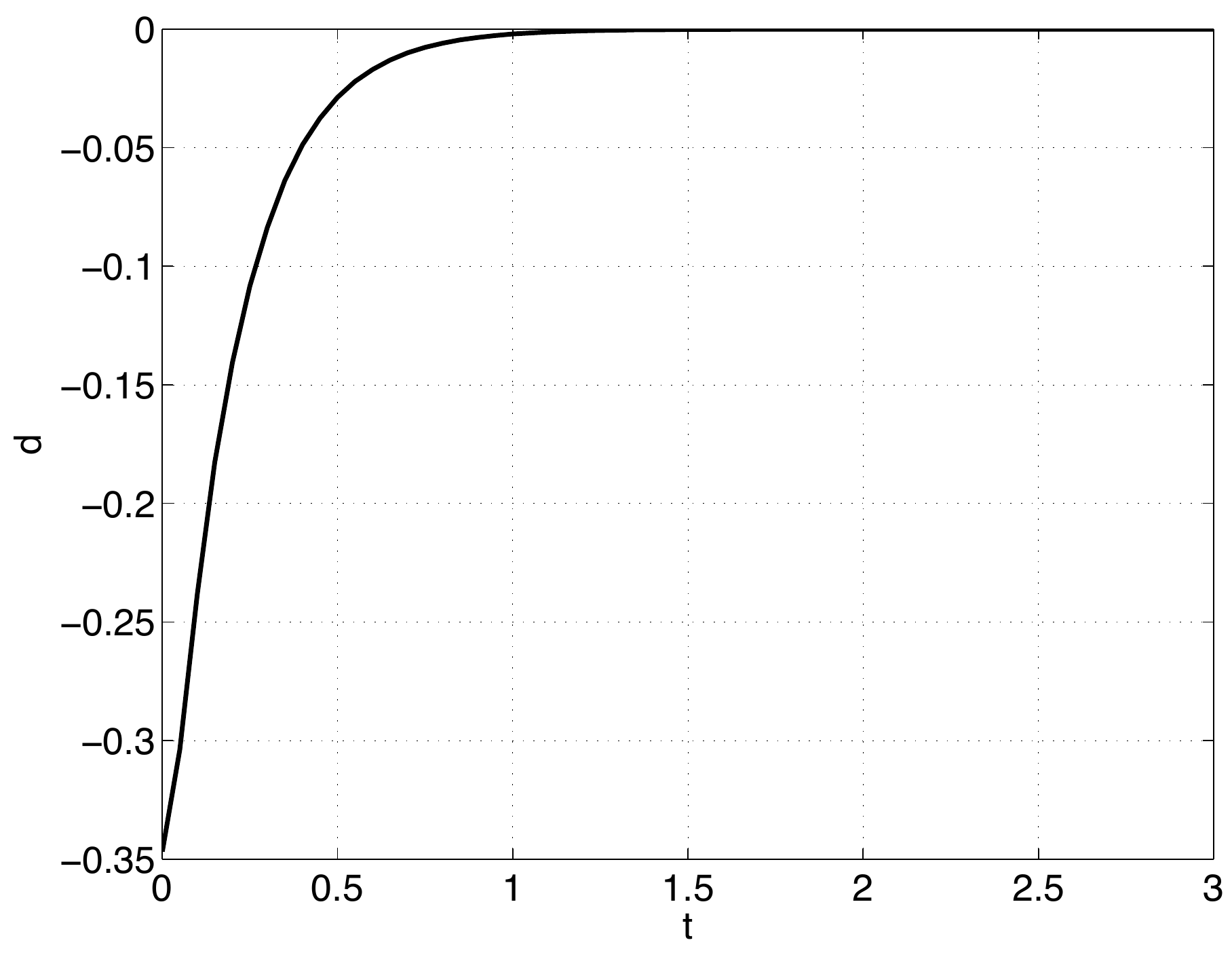}}
  \subfigure[segregated, $\Delta t=0.05s$]{\includegraphics[width=0.32\textwidth]{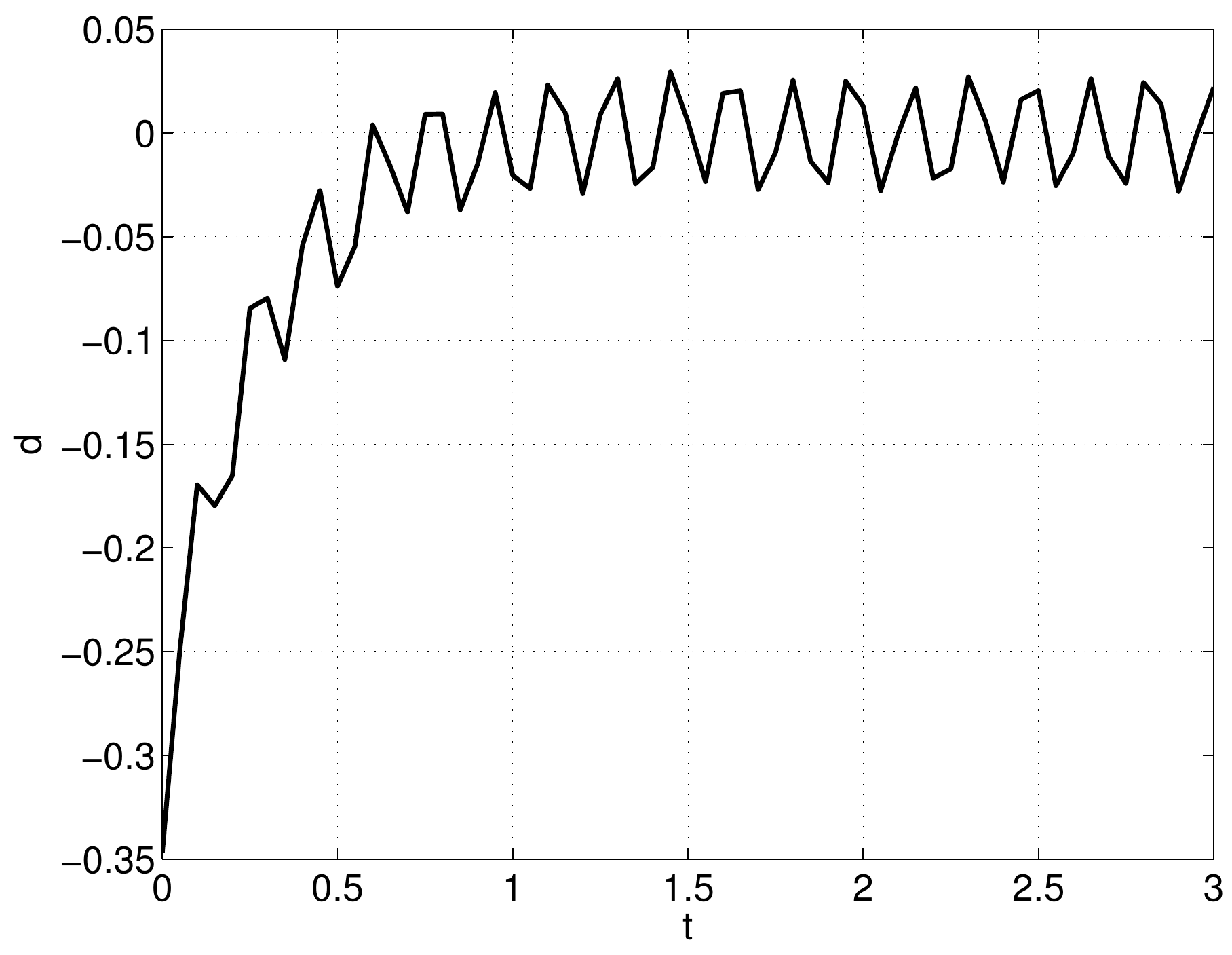}}
  \subfigure[fully impl., $\Delta t=0.1s$]{\includegraphics[width=0.32\textwidth]{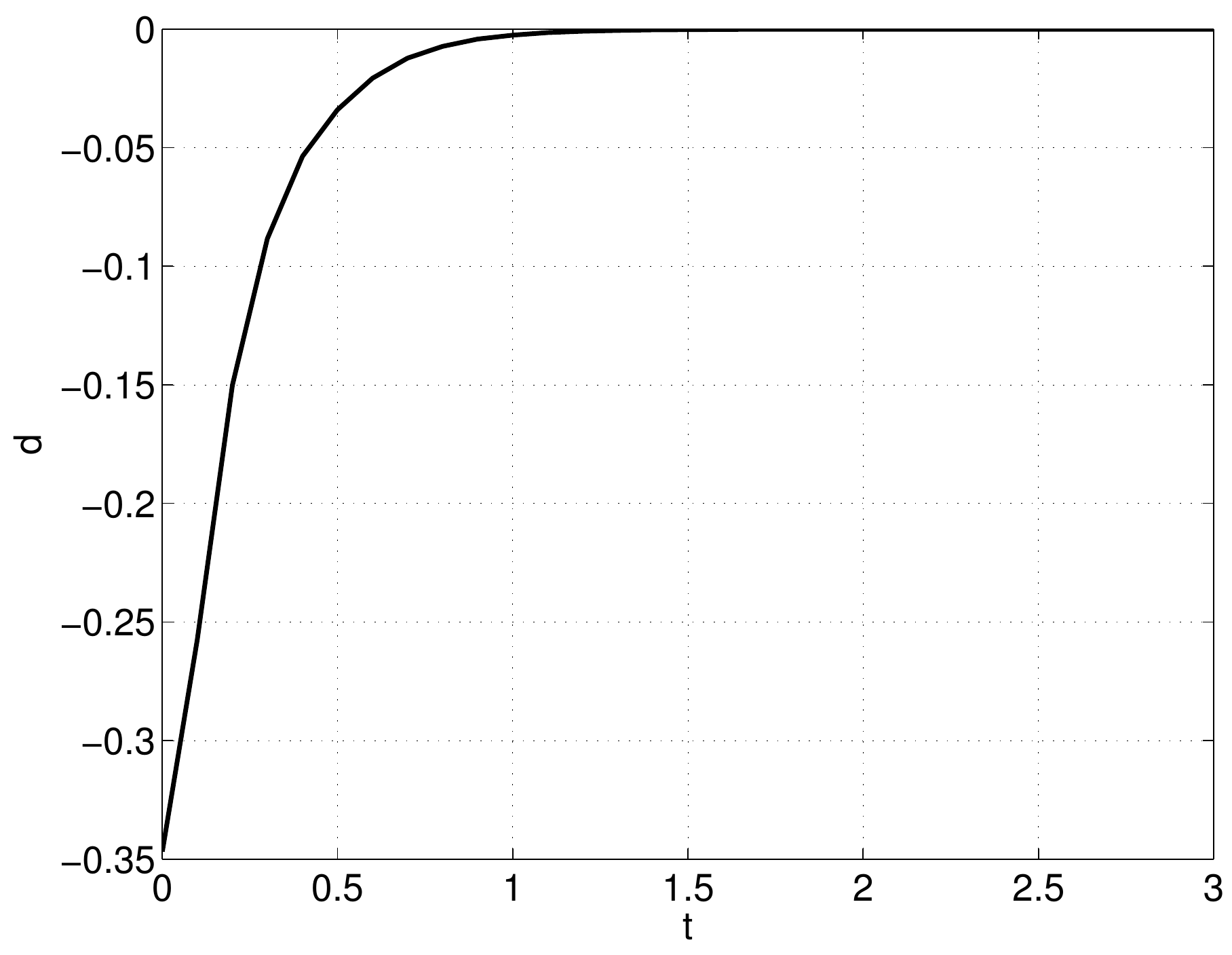}}
  \subfigure[semi-impl., $\Delta t=0.1s$]{\includegraphics[width=0.32\textwidth]{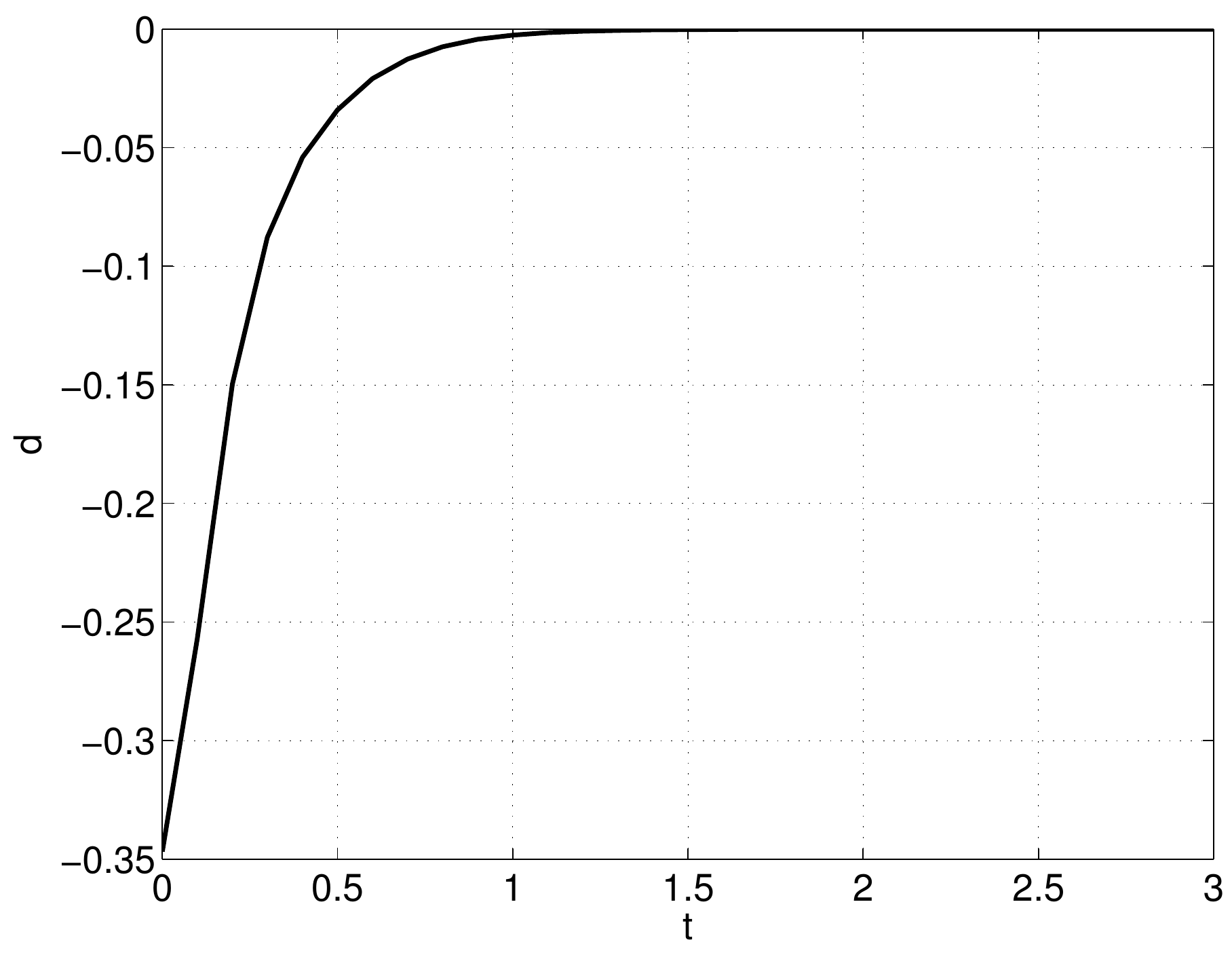}}
  \subfigure[segregated, $\Delta t=0.1s$]{\includegraphics[width=0.32\textwidth]{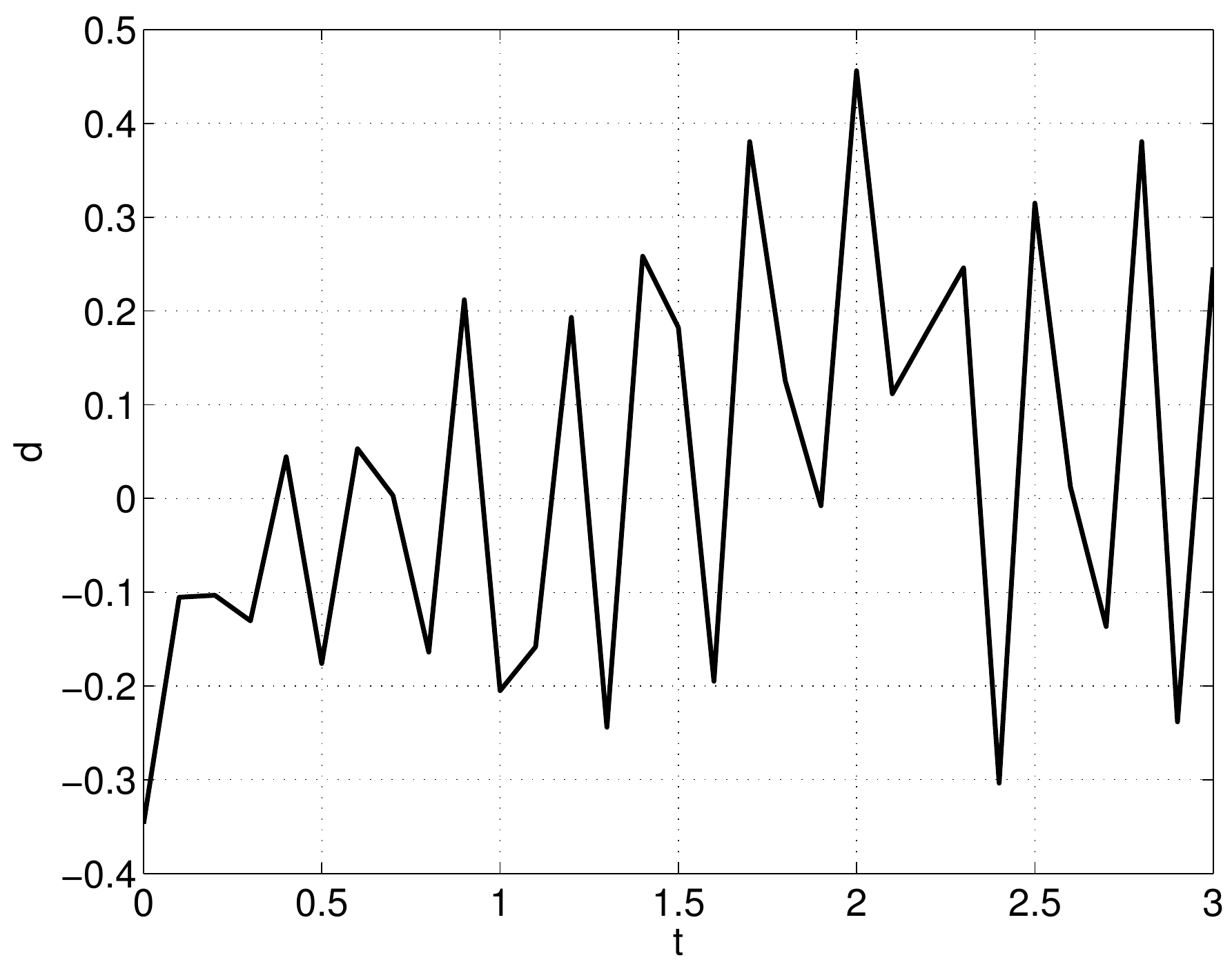}}
  \caption{Vibration of a cantilever immersed in a viscous fluid. Tip displacement for $\eta=10$ using different coupling approaches and different time steps: fully implicit (left), semi-implicit (middle), segregated (right), with time steps $\Delta t=\{0.01, 0.05, 0.1\} s$ (top, middle, bottom).} 
  \label{leaflet2}
\end{figure}

All results in Figure \ref{leaflet1} have been obtained by the semi-implicit approach \eqref{eq:operator-motion-time-discrete-semi-implicit}. 
% end new
%These results have been obtained by the semi-implicit approach \eqref{eq:operator-motion-time-discrete-semi-implicit}. 
%As can be seen, good results without artificial oscillations are obtained for all cases. 
%\begin{figure}\centering
%  \subfigure[$\eta=0.001\, Pa\cdot s$]{\includegraphics[width=0.49\textwidth]{figures/leaflet_new_eta1em3_fsi1_damp1_dt100_d.pdf}}
%  \subfigure[$\eta=0.1\, Pa\cdot s$]{\includegraphics[width=0.49\textwidth]{figures/leaflet_new_eta1em1_fsi1_damp1_dt100_d.pdf}}
%  \subfigure[$\eta=1\, Pa\cdot s$]{\includegraphics[width=0.49\textwidth]{figures/leaflet_new_eta1_fsi1_damp1_dt100_d.pdf}}
%  \subfigure[$\eta=10\, Pa\cdot s$]{\includegraphics[width=0.49\textwidth]{figures/leaflet_new_eta1ep1_fsi1_damp1_dt100_d.pdf}}
%  \caption{Vibration of a cantilever. Tip displacement for different viscosities and $\Delta t=0.01s$, obtained using the semi-implicit approach.} 
%  \label{leaflet1}
%\end{figure}
%Furthermore, 
In a next step, we use this example in the very viscous case ($\eta=10\, Pa\cdot s$) in order to compare the three different coupling approaches, i.e., the fully implicit \eqref{eq:inexact-Jacobian}, the semi-implicit \eqref{eq:operator-motion-time-discrete-semi-implicit}, and the segregated approach \eqref{eq:operator-motion-time-discrete-leap-frog-solid-dominates}. For all cases, we consider different time steps $\Delta t=\{0.01, 0.05, 0.1\} s$. The results are gathered in Figure \ref{leaflet2}. As can be seen, the results for the fully implicit and the semi-implicit approach are identical for all cases. The segregated approach yields stable results only for the smallest time step, while strong spurious oscillations appear when the time step is increased.
We highlight that the computational cost of the semi-implicit approach is the same as in the segregated approach, assembling the matrices of the fluid problem only once per time step, and, therefore, significantly less than in the fully implicit approach where the fluid problem matrices are assembled in each Newton iteration. Accordingly, the semi-implicit approach appears to be a very efficient alternative combining the cost-effectiveness of the segregated approach with the accuracy and stability of the fully implicit approach. For the remainder of this paper, we use the semi-implicit formulation in all computations.
% ----------------------------------------------------------------------------------------------------------------------------------------
%  HONEYSPOON
\begin{figure}\centering
  \includegraphics[width=0.37625\textwidth]{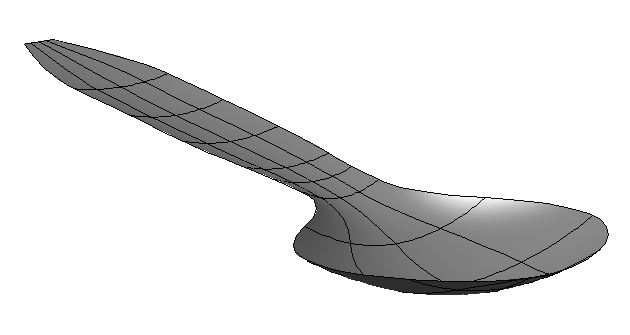} \hspace{1cm}
  \includegraphics[width=0.35\textwidth]{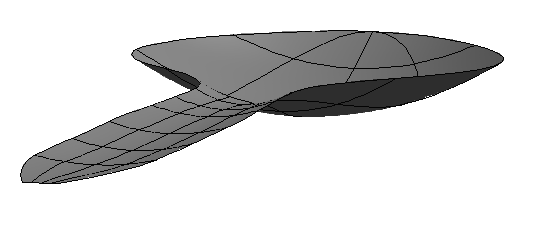}
  \includegraphics[width=0.525\textwidth]{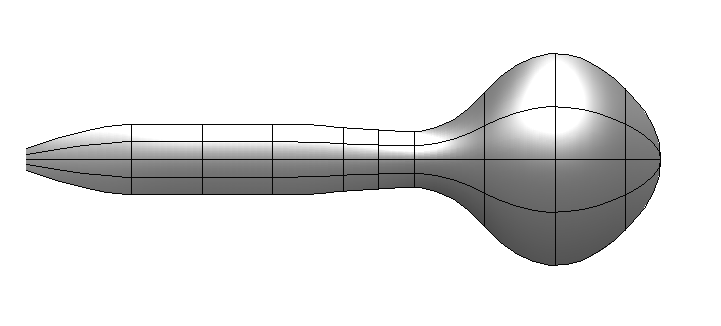}
  \caption{Honey-spoon geometry.} 
  \label{spoon_setup_geometry}
\end{figure}
\begin{figure}\centering
  \includegraphics[height=0.5\textwidth]{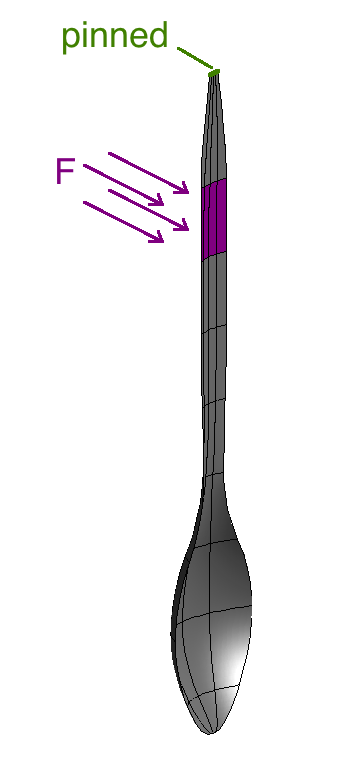} \hspace{1cm}
  \includegraphics[height=0.5\textwidth]{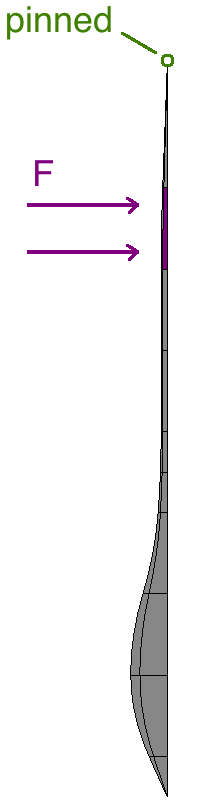}
  \caption{Honey-spoon problem setup.} 
  \label{spoon_setup}
\end{figure}
% spoon deformation
\begin{figure}\centering
  \subfigure[$t=0s$]{\includegraphics[width=0.32\textwidth]{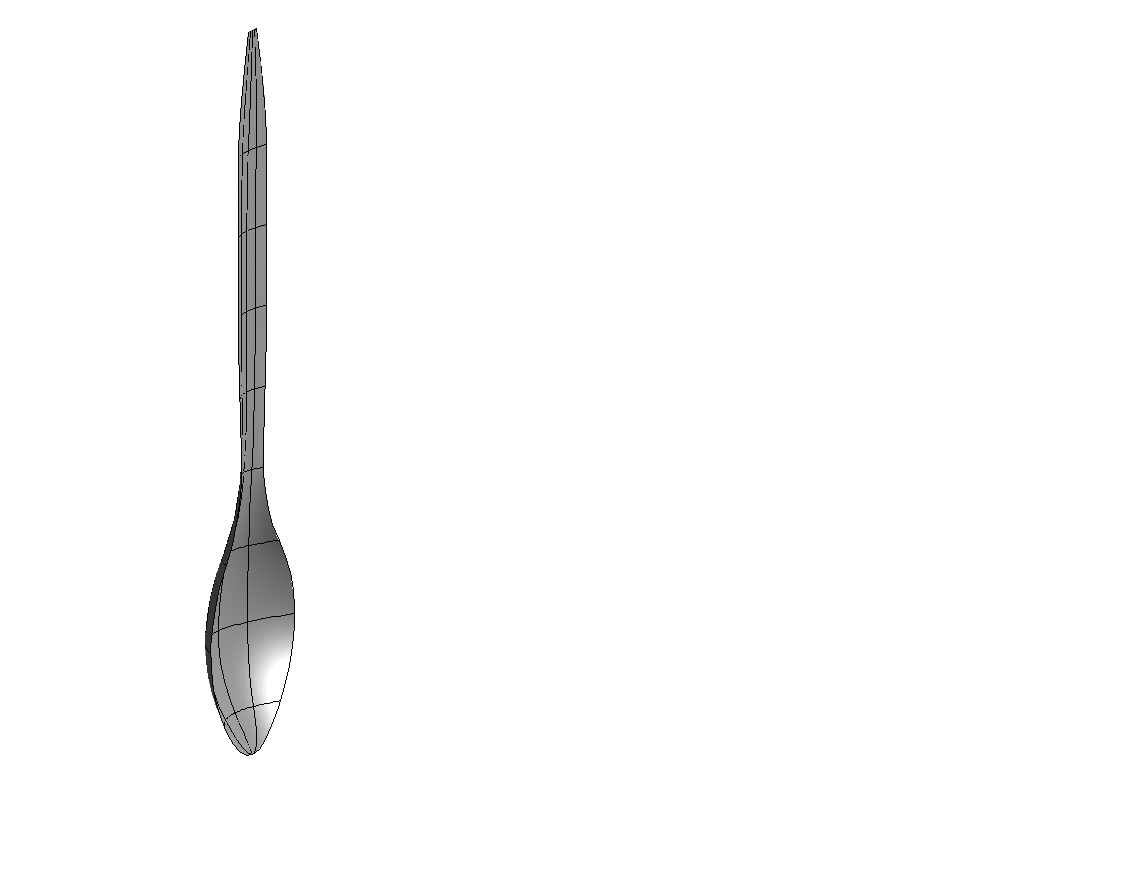}}
  \subfigure[$t=0.1s$]{\includegraphics[width=0.32\textwidth]{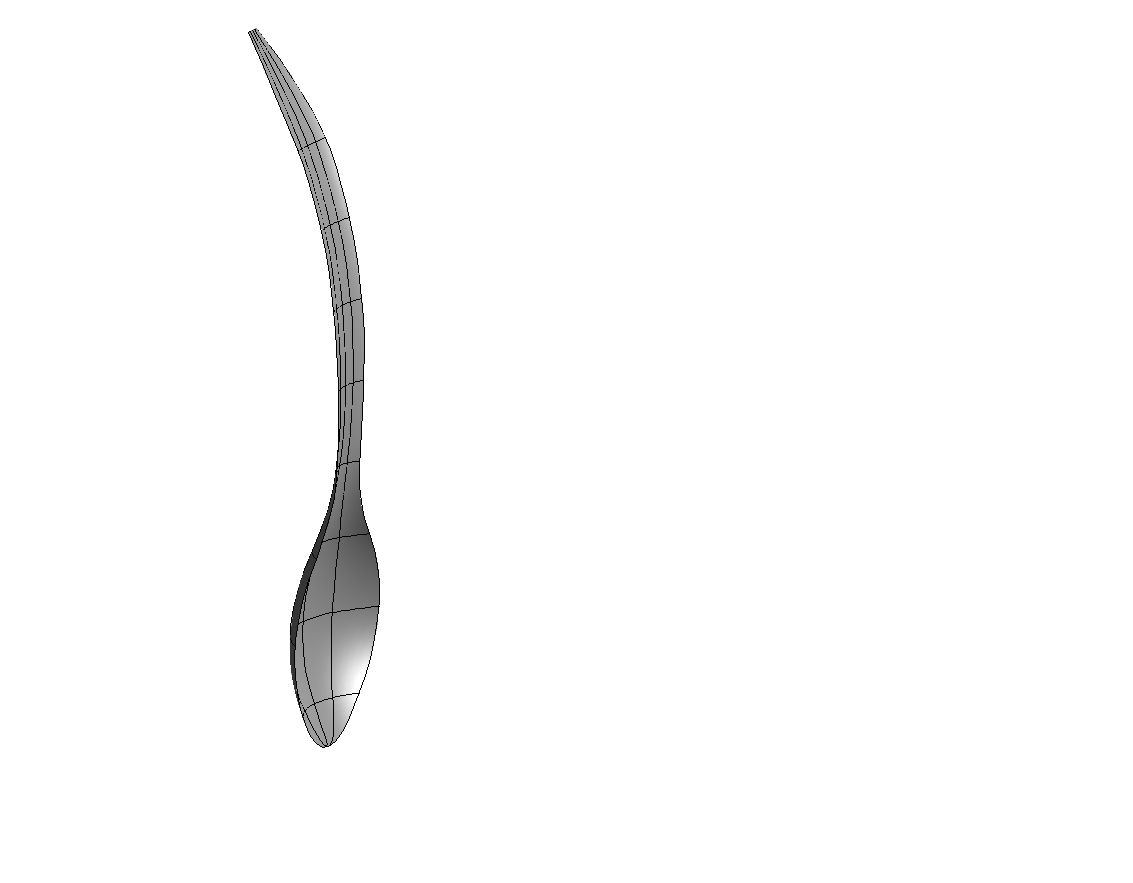}}
  \subfigure[$t=0.2s$]{\includegraphics[width=0.32\textwidth]{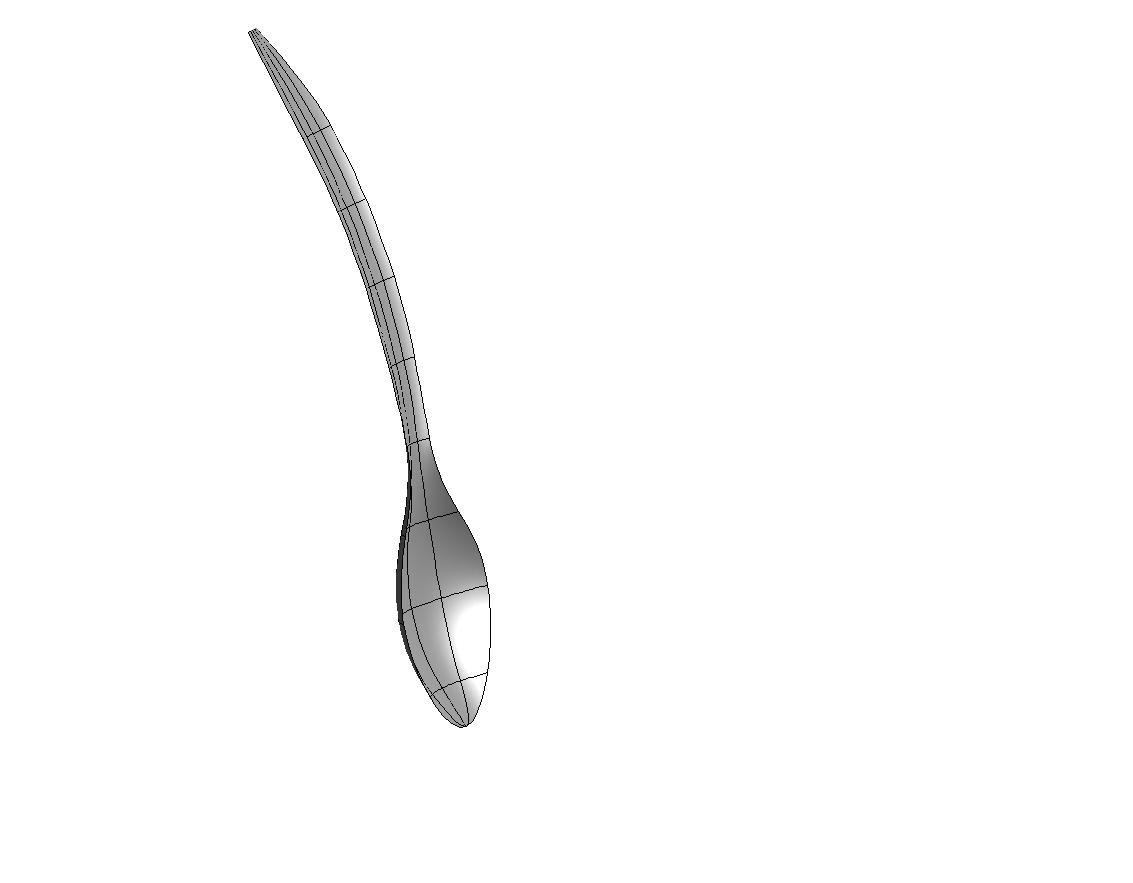}}
  \subfigure[$t=0.3s$]{\includegraphics[width=0.32\textwidth]{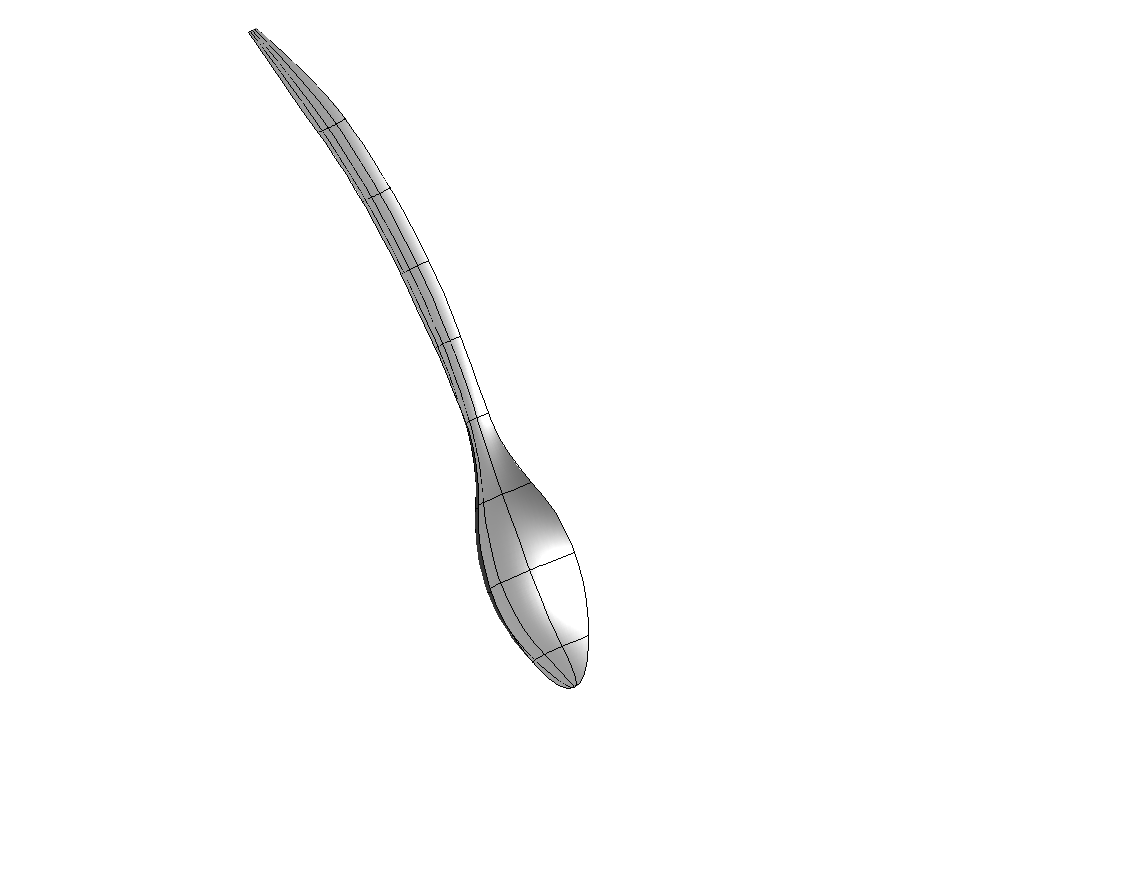}}
  \subfigure[$t=0.4s$]{\includegraphics[width=0.32\textwidth]{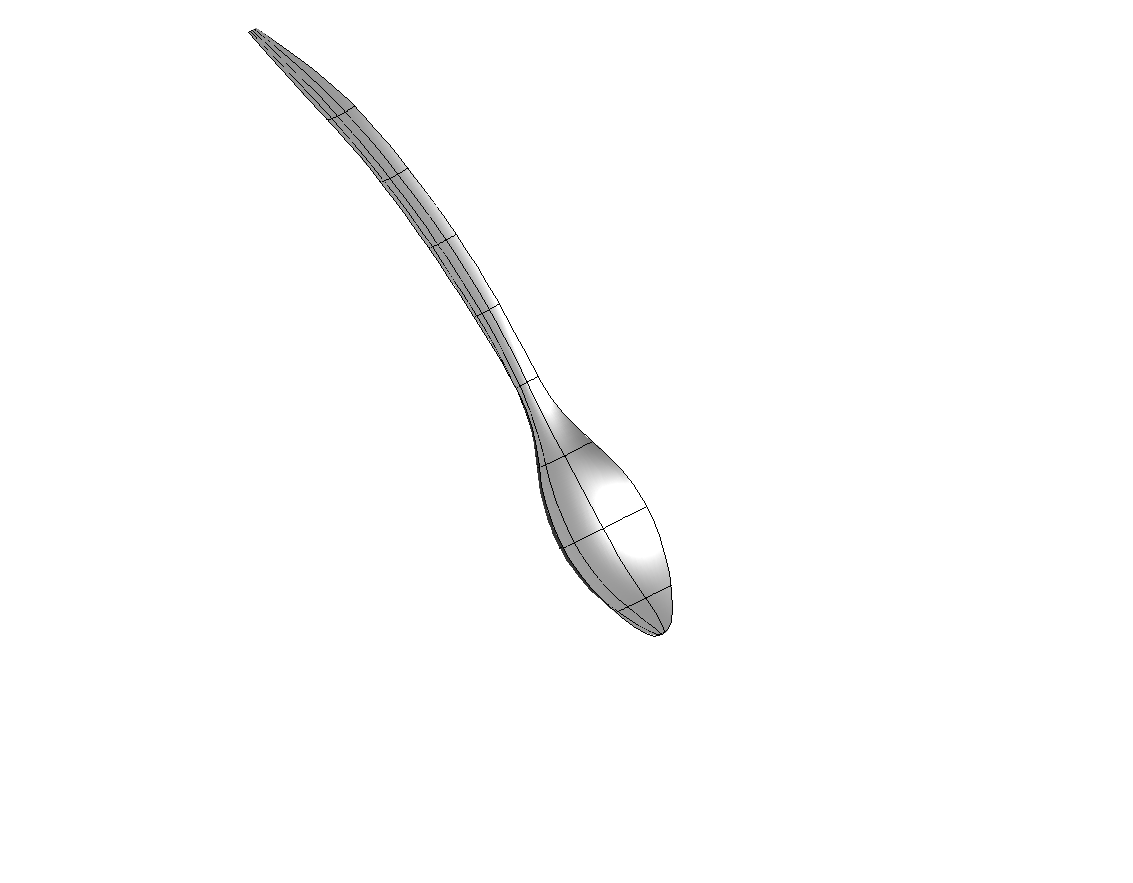}}
  \subfigure[$t=0.5s$]{\includegraphics[width=0.32\textwidth]{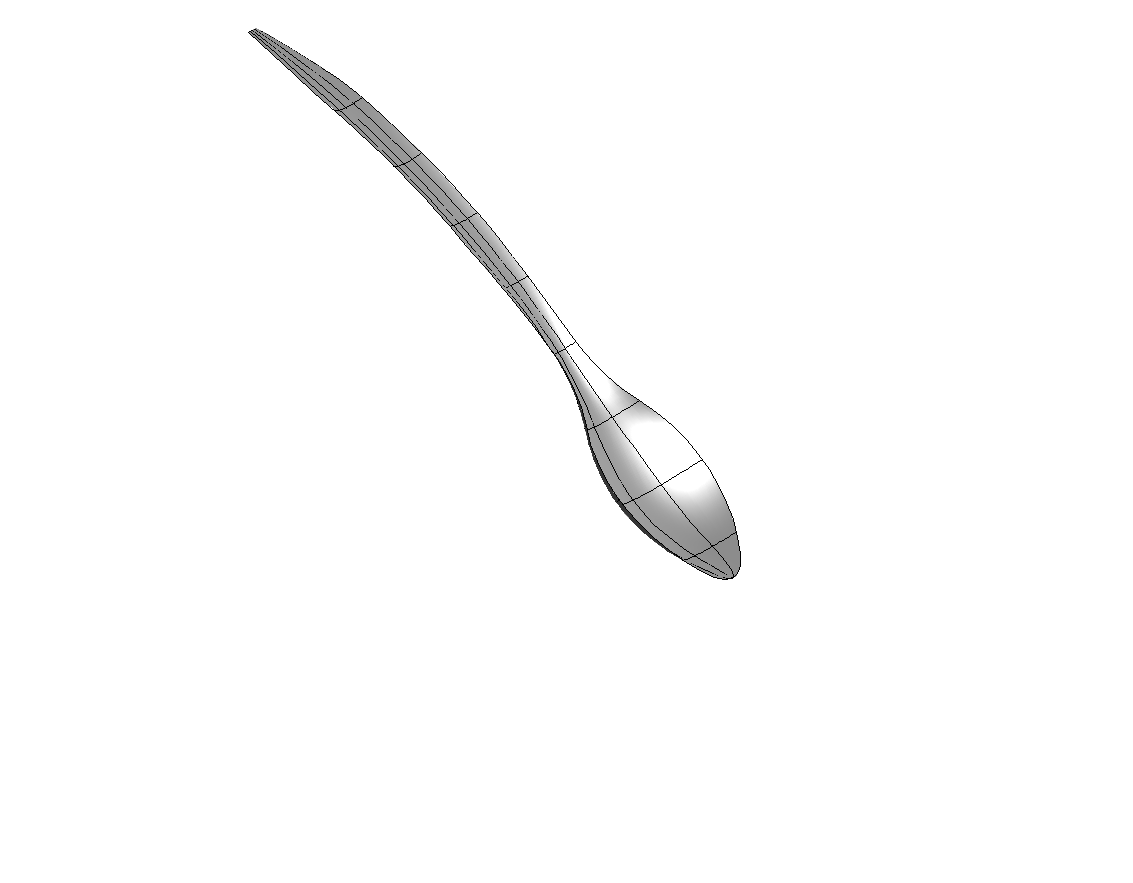}}
  \subfigure[$t=0.6s$]{\includegraphics[width=0.32\textwidth]{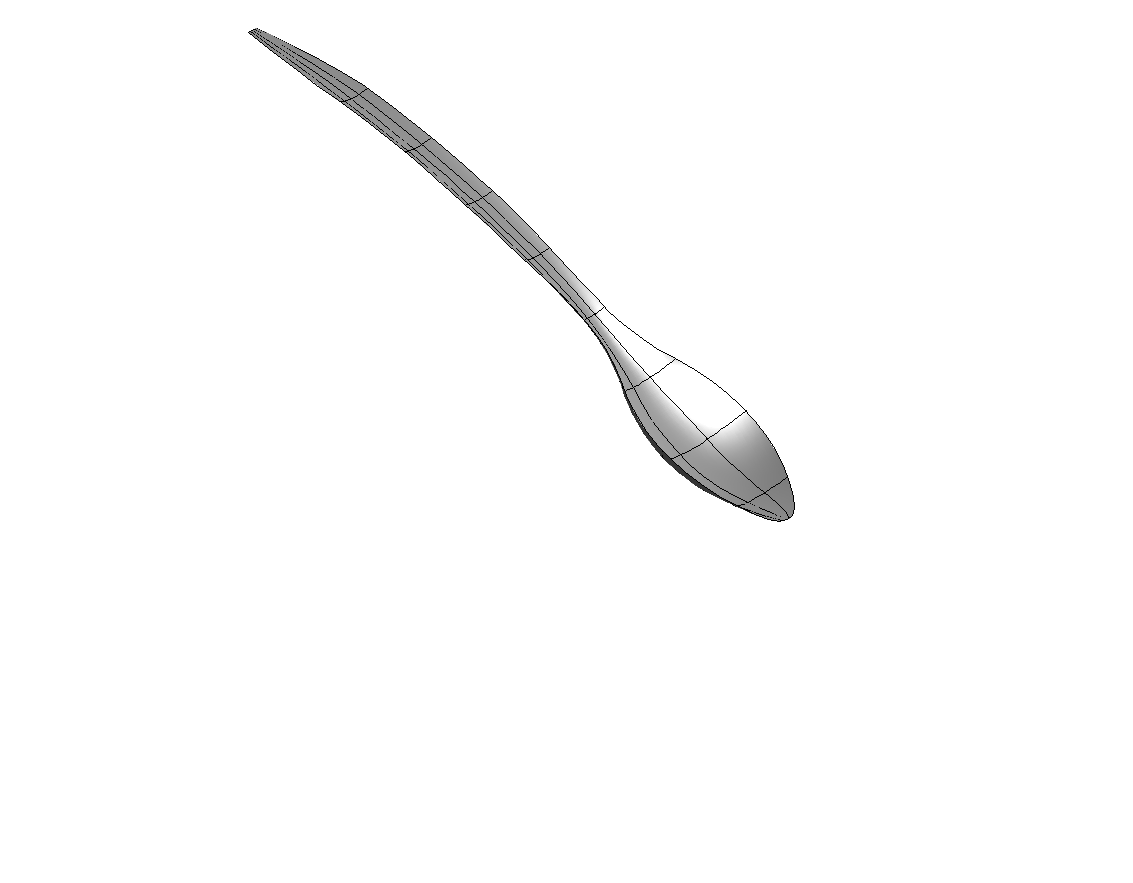}}
  \subfigure[$t=0.7s$]{\includegraphics[width=0.32\textwidth]{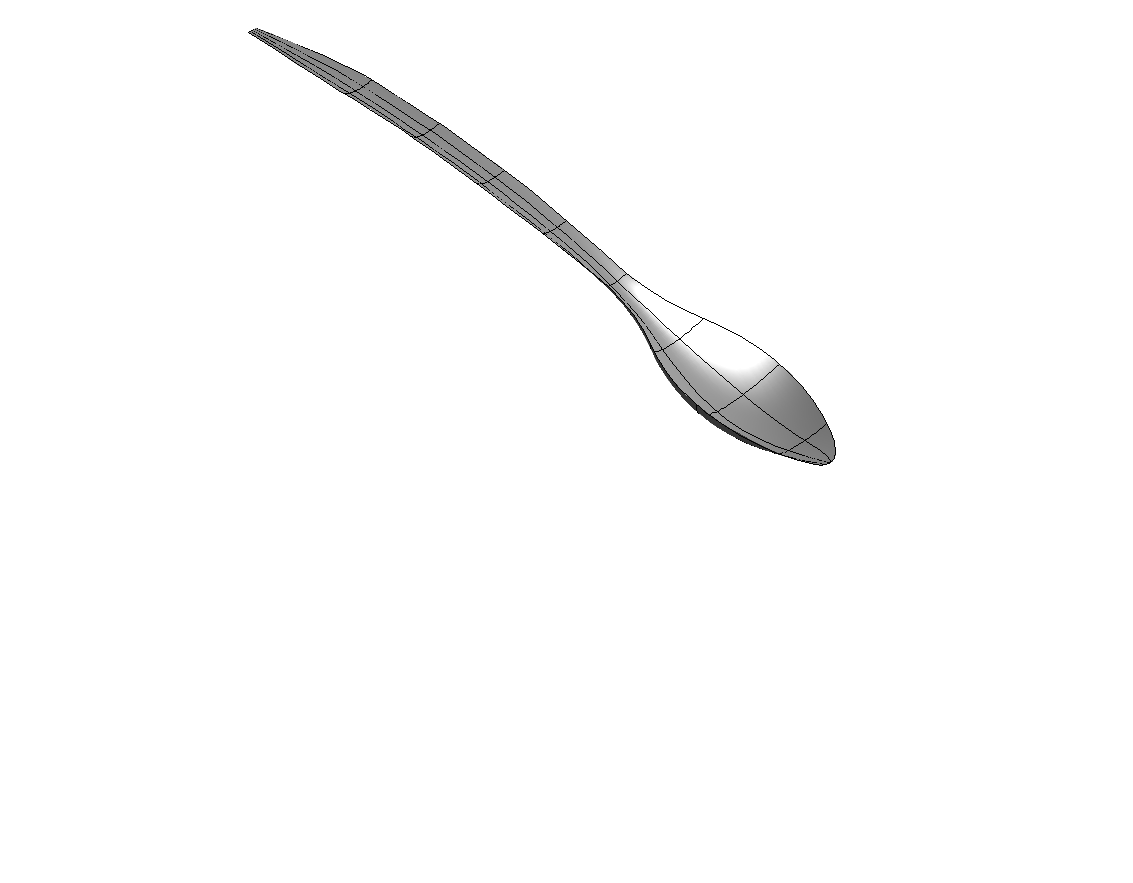}}
  \subfigure[$t=0.8s$]{\includegraphics[width=0.32\textwidth]{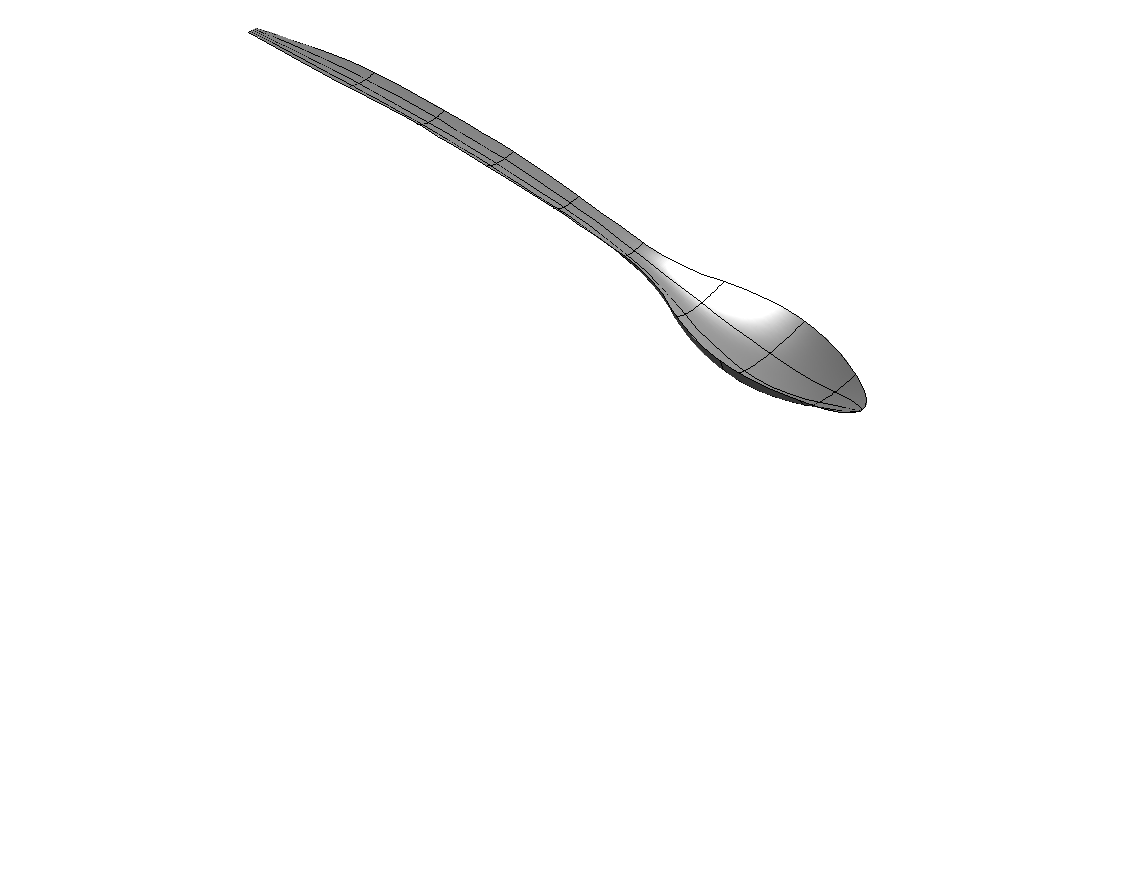}}
  \subfigure[$t=0.9s$]{\includegraphics[width=0.32\textwidth]{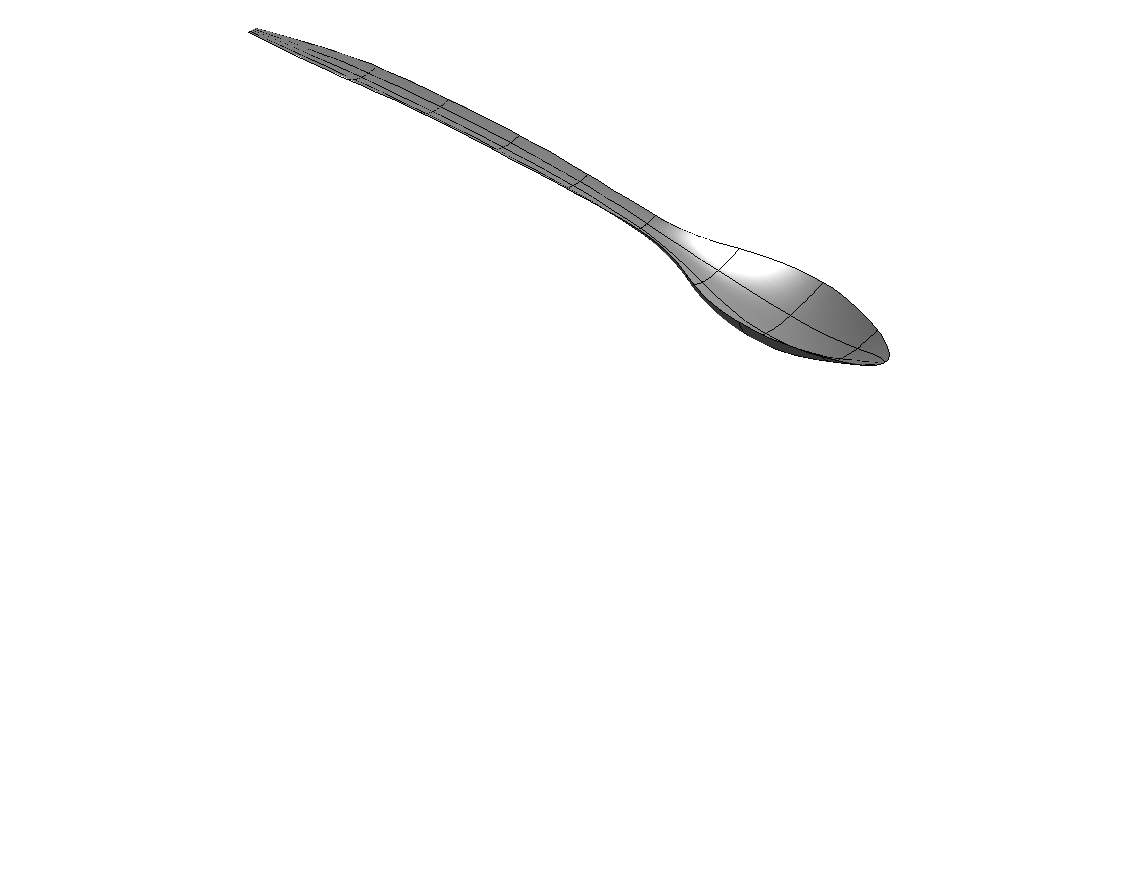}}
  \subfigure[$t=1.0s$]{\includegraphics[width=0.32\textwidth]{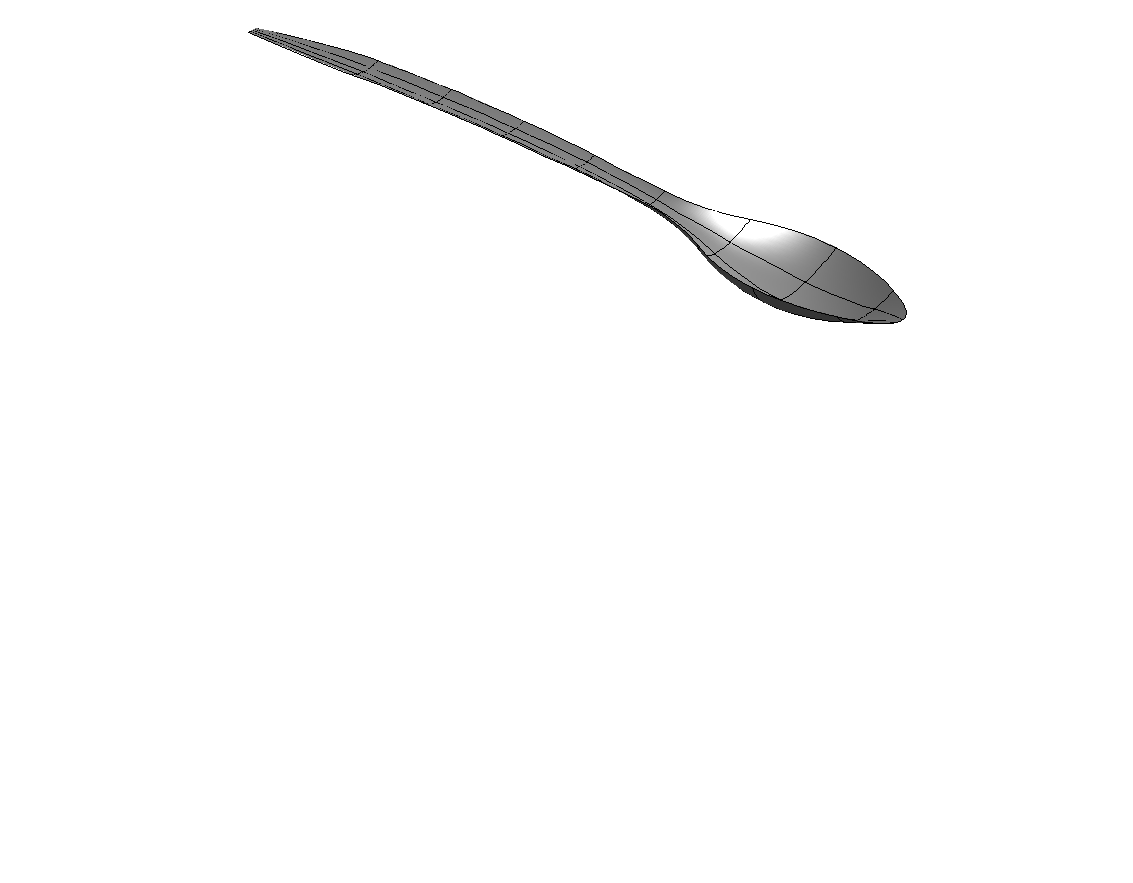}}
  \subfigure[$t=1.1s$]{\includegraphics[width=0.32\textwidth]{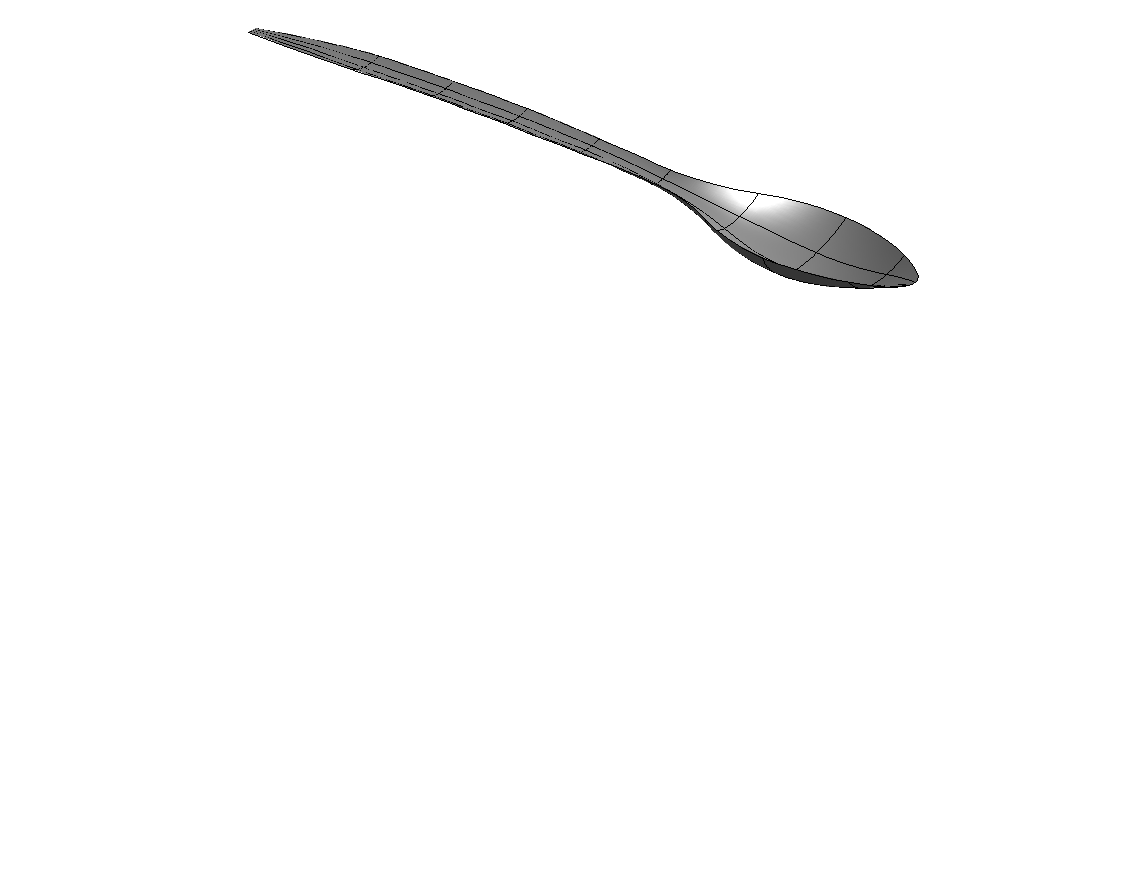}}
  \caption{Honey-spoon deformation plots.} 
  \label{spoon_def}
\end{figure}
% spoon external flow
\begin{figure}\centering
  \subfigure[$t=0.1s$]{\includegraphics[scale=.8]{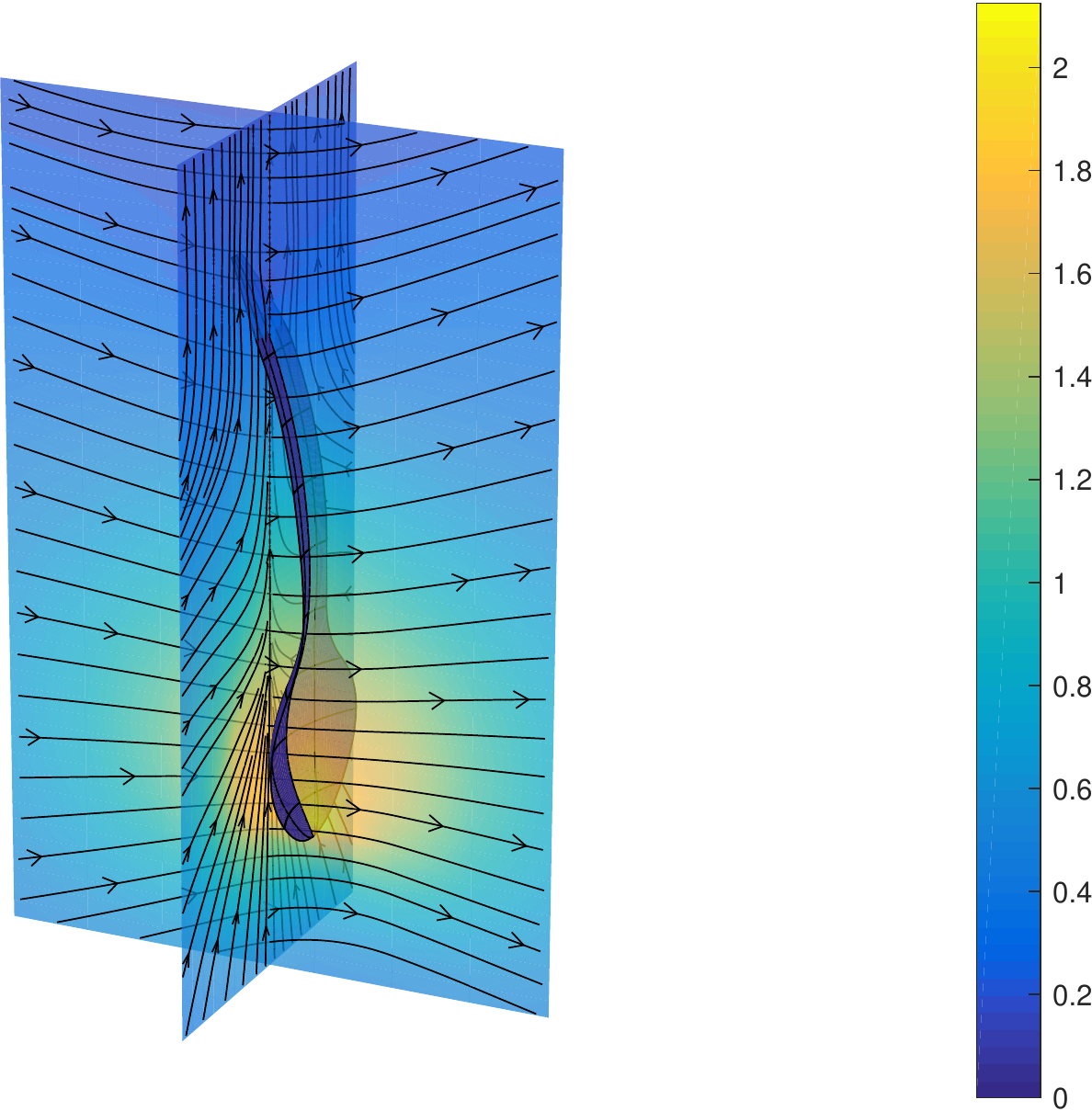}}

  \subfigure[$t=0.7s$]{\includegraphics[scale=.8]{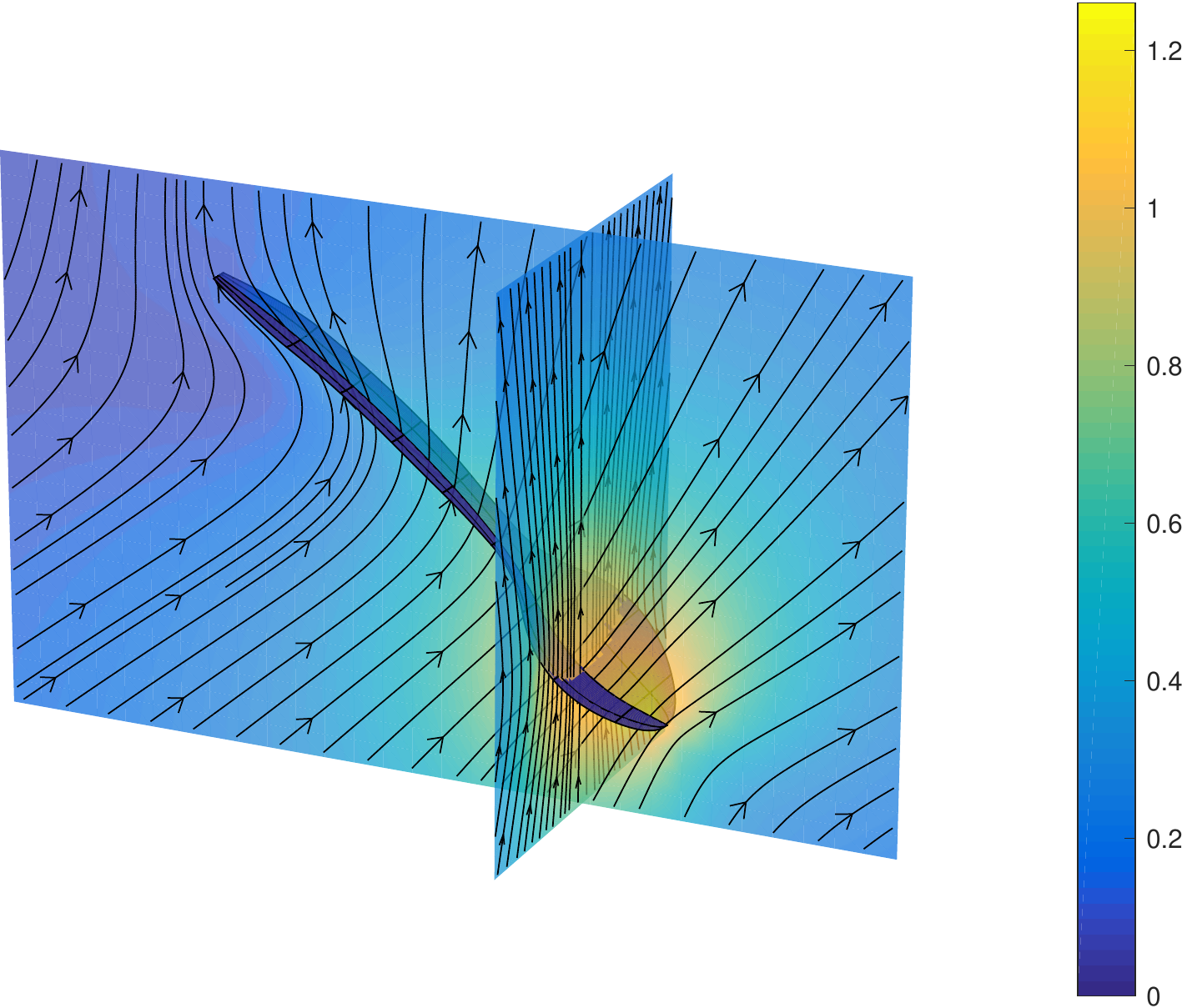}}
  \caption{Honey-spoon streamline plots. Streamlines and velocity magnitude plots, projected on two perpendicular planes. Velocity is measured in $m/s$.}
  \label{spoon_streamlines}
\end{figure}
\subsection{Honey-spoon}
In this example, we model the fluid-structure interaction of a spoon
moving through honey. The spoon geometry is shown in Figures
\ref{spoon_setup_geometry} and \ref{spoon_setup}, it is $4.5\, cm$
long with a thickness $h=0.2\, mm$ and material parameters
$E=2.8\cdot10^9\, N/mm^2, \nu=0.39, \rho=1.13\, kg/m^3$. The viscosity
of the honey is taken as $\eta=5.0 \,Pa\cdot s$. The problem setup and
boundary conditions are shown in Figure \ref{spoon_setup}, with
$F=7000\,N/m^2$, corresponding to a total load of $0.177\,N$. In
Figure \ref{spoon_def} the deformation at different time steps is
displayed, where a time step of $\Delta t=0.1\, s$ has been used.
Figure~\ref{spoon_streamlines} presents the details of the flow
streamlines and velocity magnitude for two significant time steps.
% ----------------------------------------------------------------------------------------------------------------------------------------
%  FALLING CAP 
\begin{figure}\centering
  \subfigure[$t=1s$]{\includegraphics[width=0.32\textwidth]{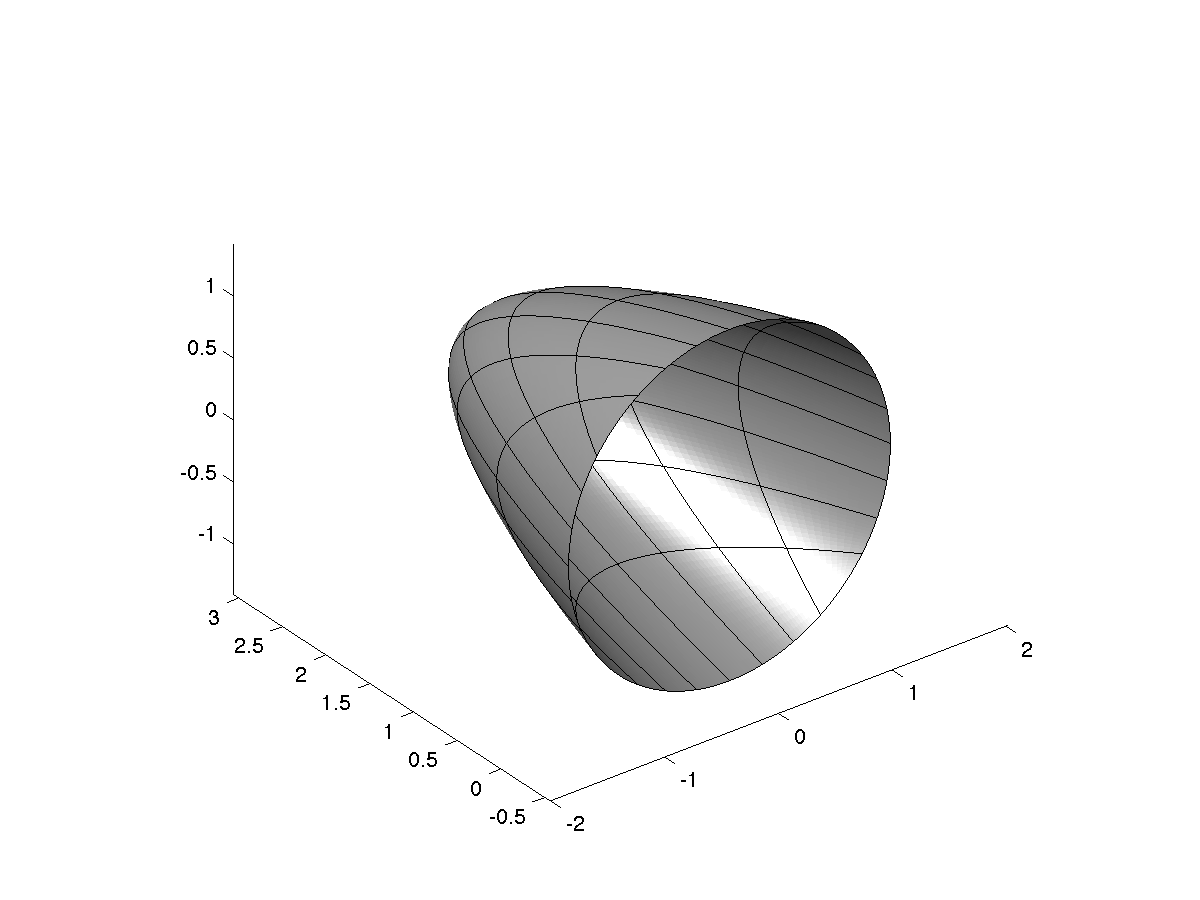}}
  \subfigure[$t=51s$]{\includegraphics[width=0.32\textwidth]{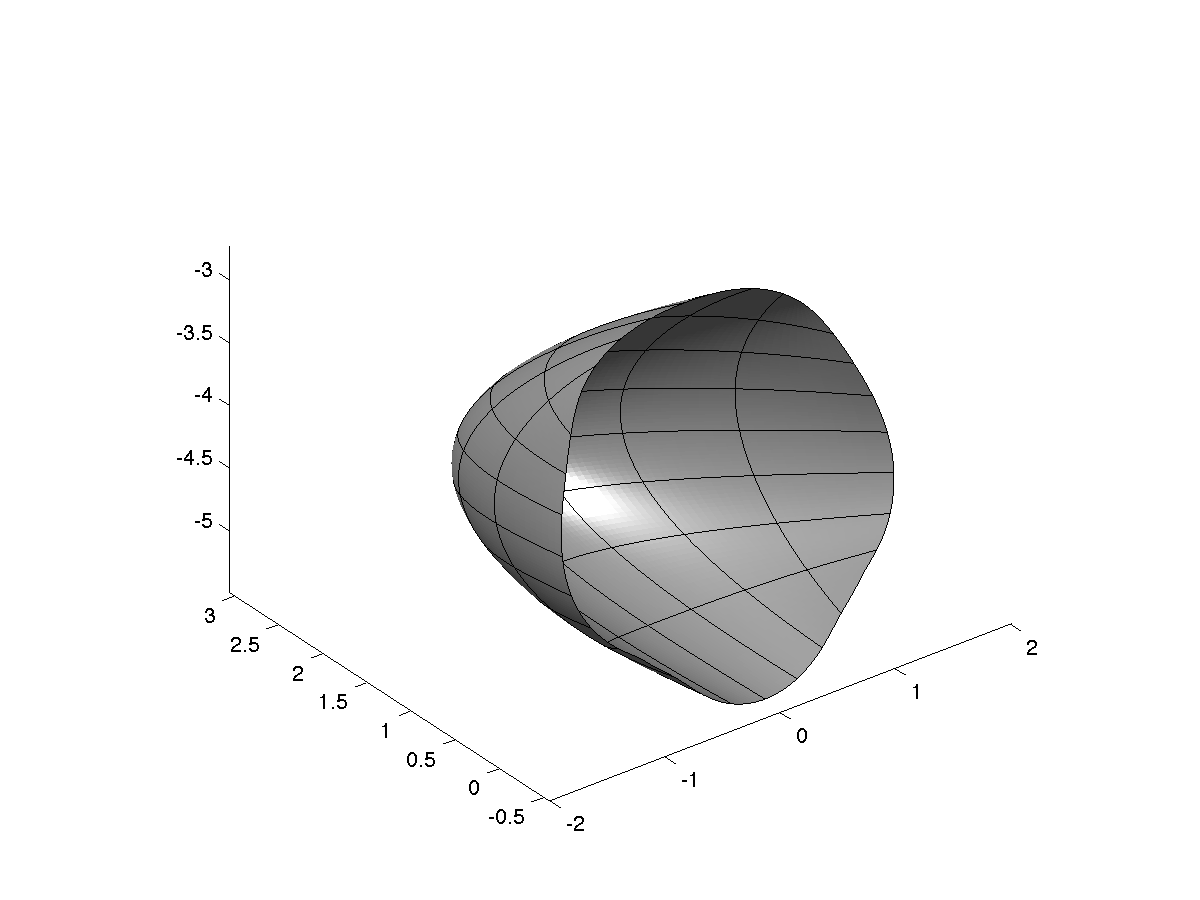}}
  \subfigure[$t=101s$]{\includegraphics[width=0.32\textwidth]{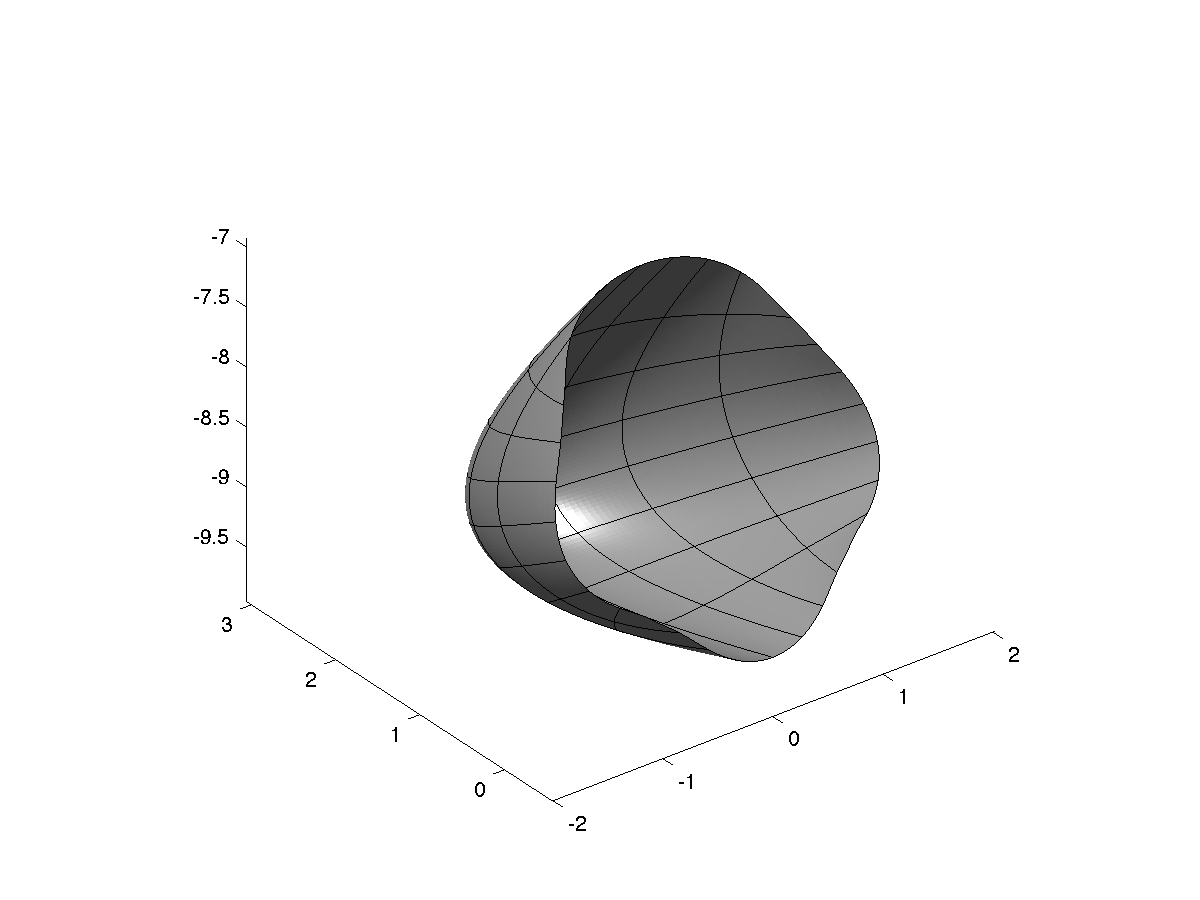}}
  \subfigure[$t=151s$]{\includegraphics[width=0.32\textwidth]{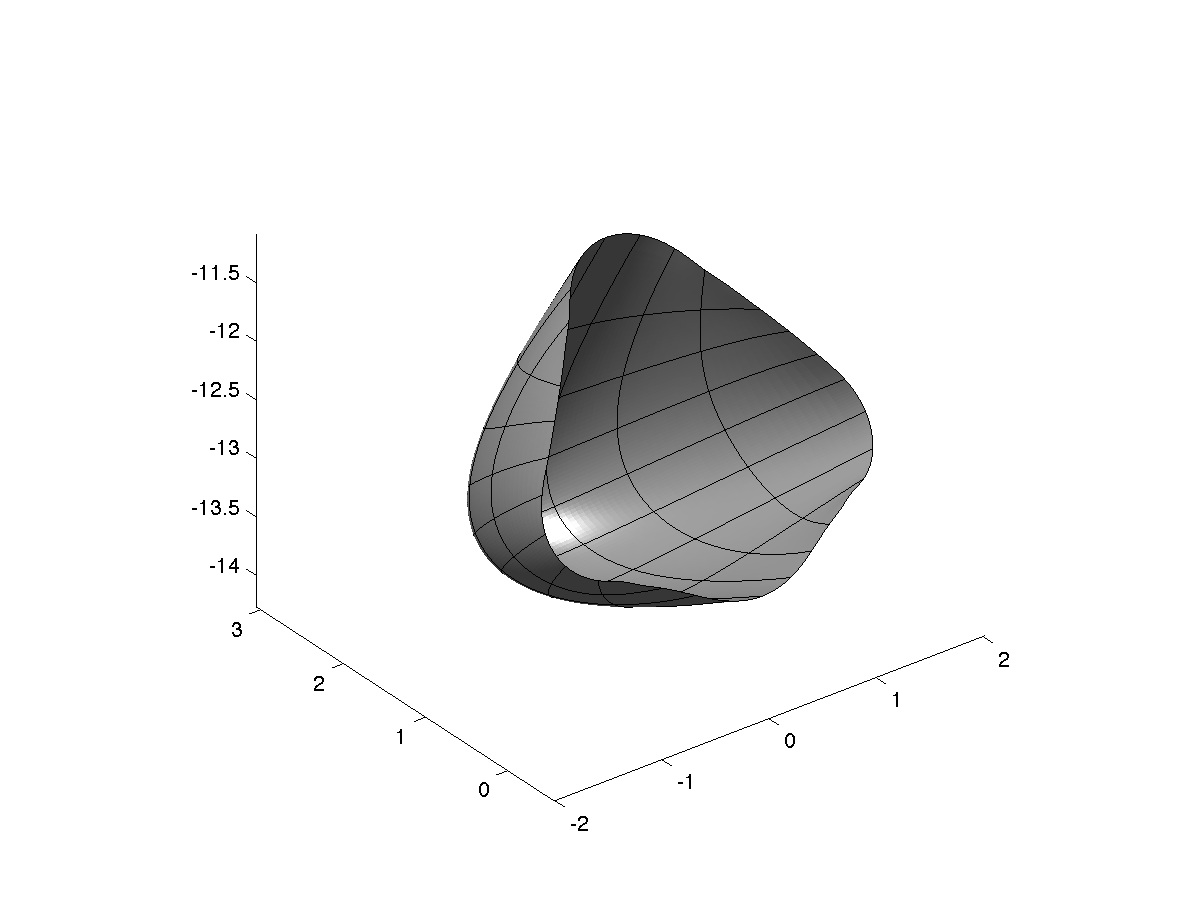}}
  \subfigure[$t=201s$]{\includegraphics[width=0.32\textwidth]{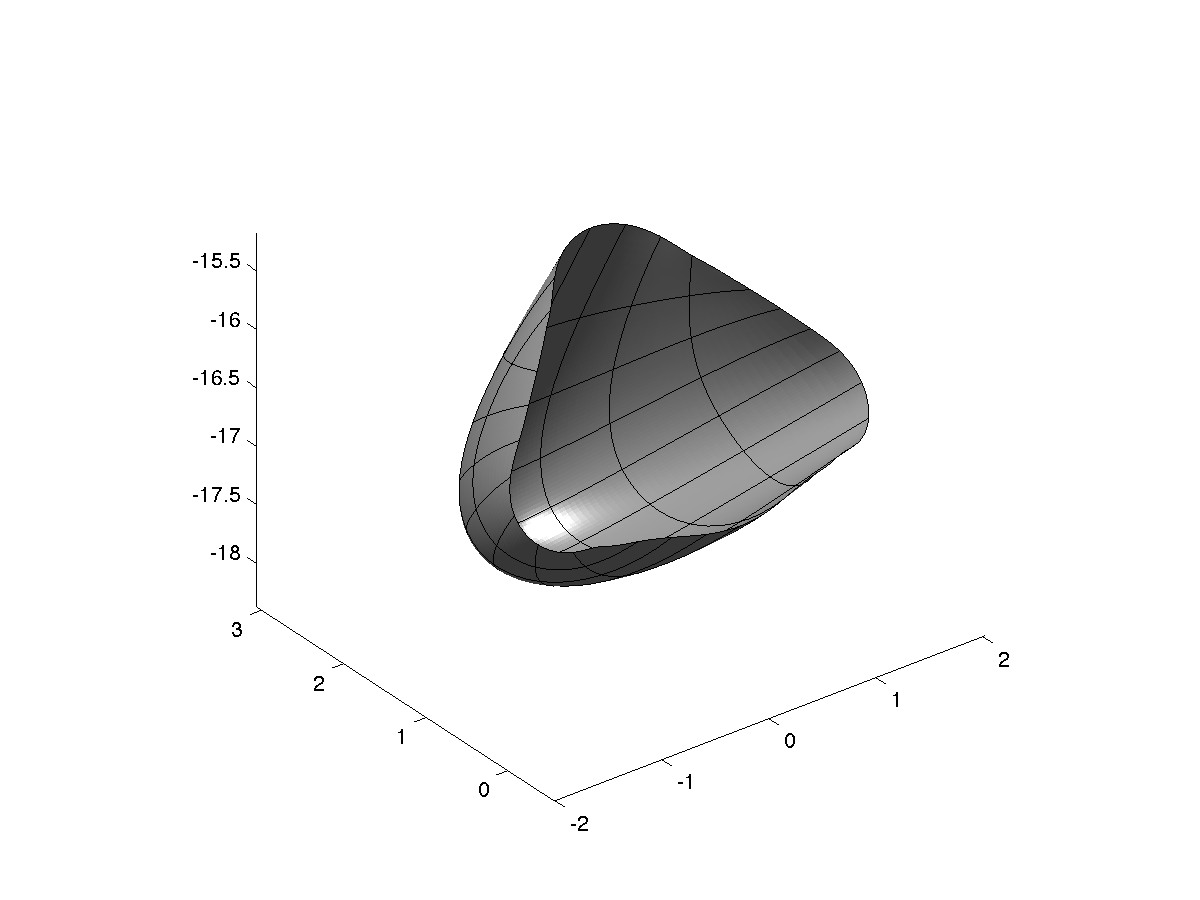}}
  \subfigure[$t=251s$]{\includegraphics[width=0.32\textwidth]{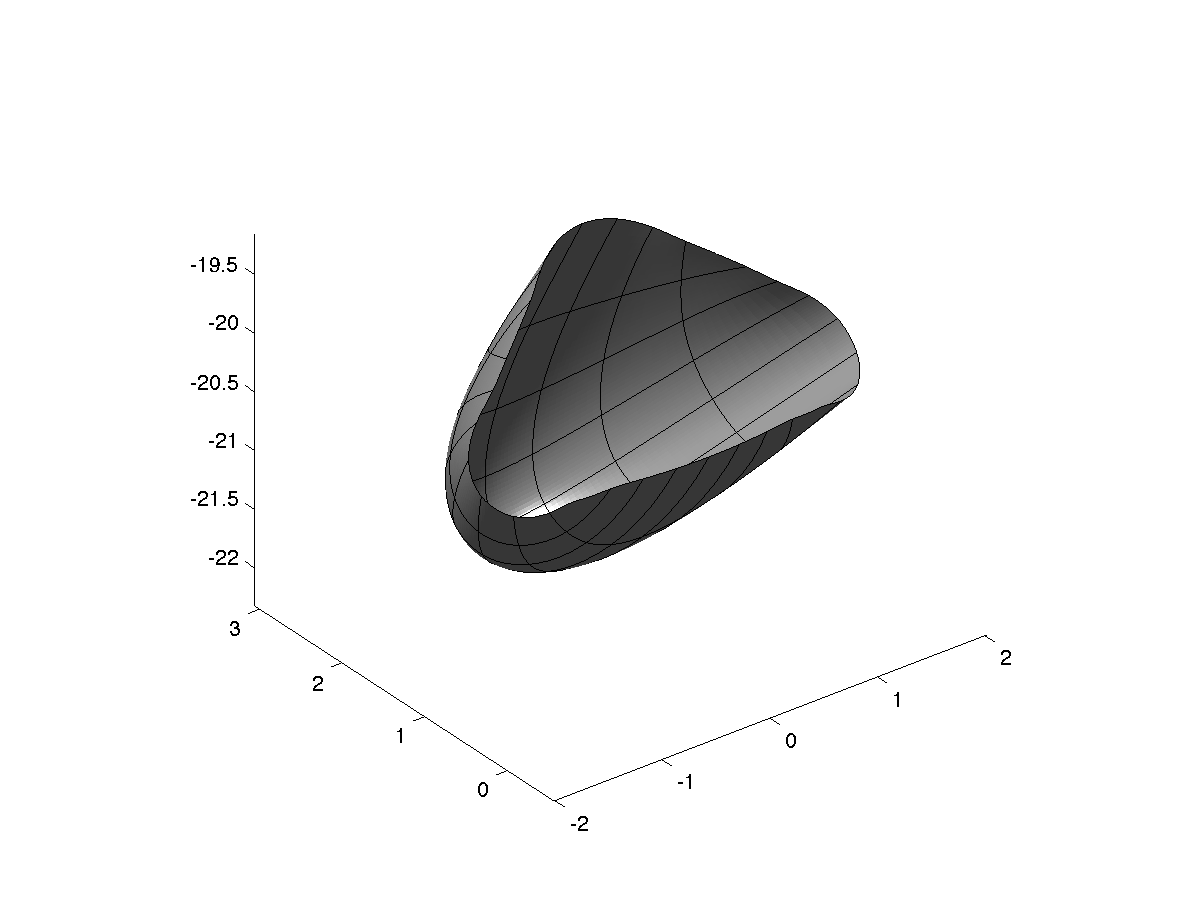}}
  \subfigure[$t=301s$]{\includegraphics[width=0.32\textwidth]{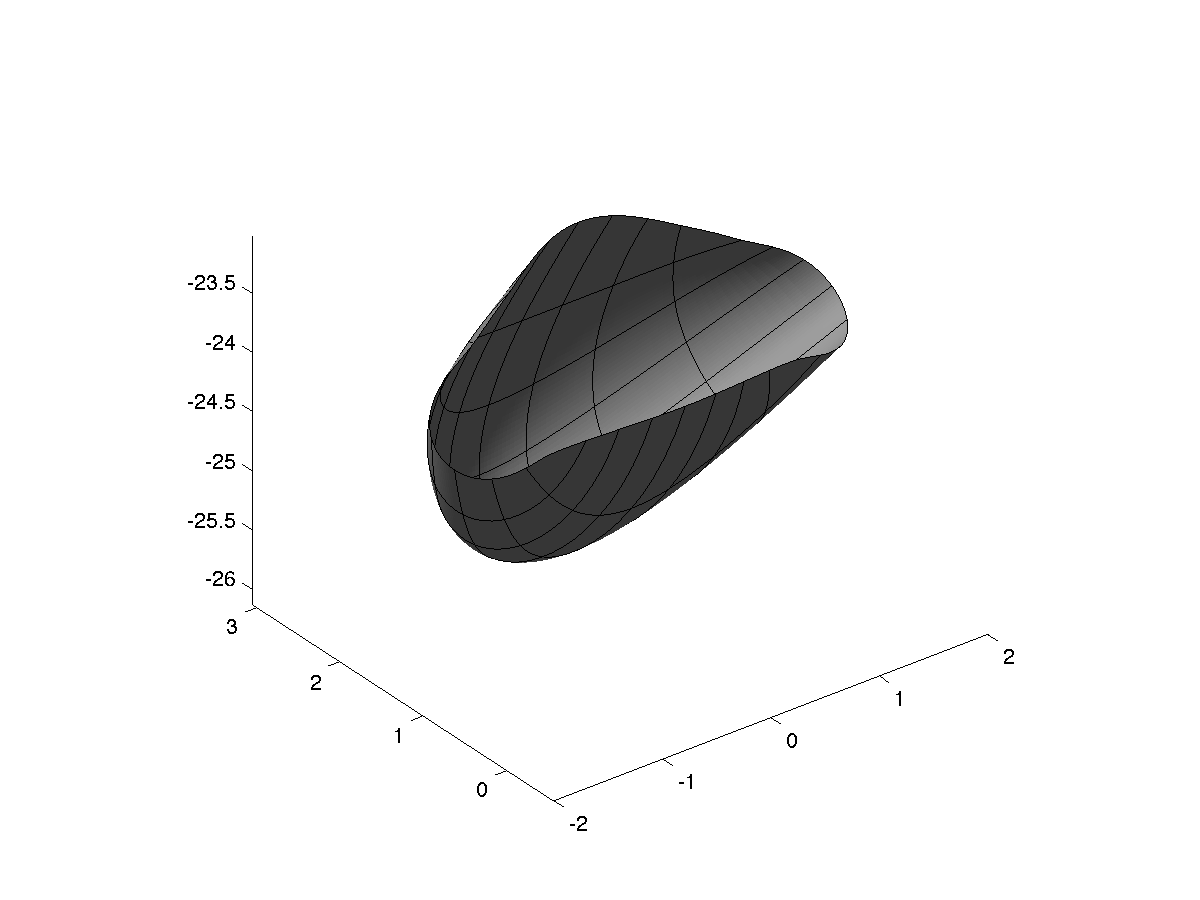}}
  \subfigure[$t=351s$]{\includegraphics[width=0.32\textwidth]{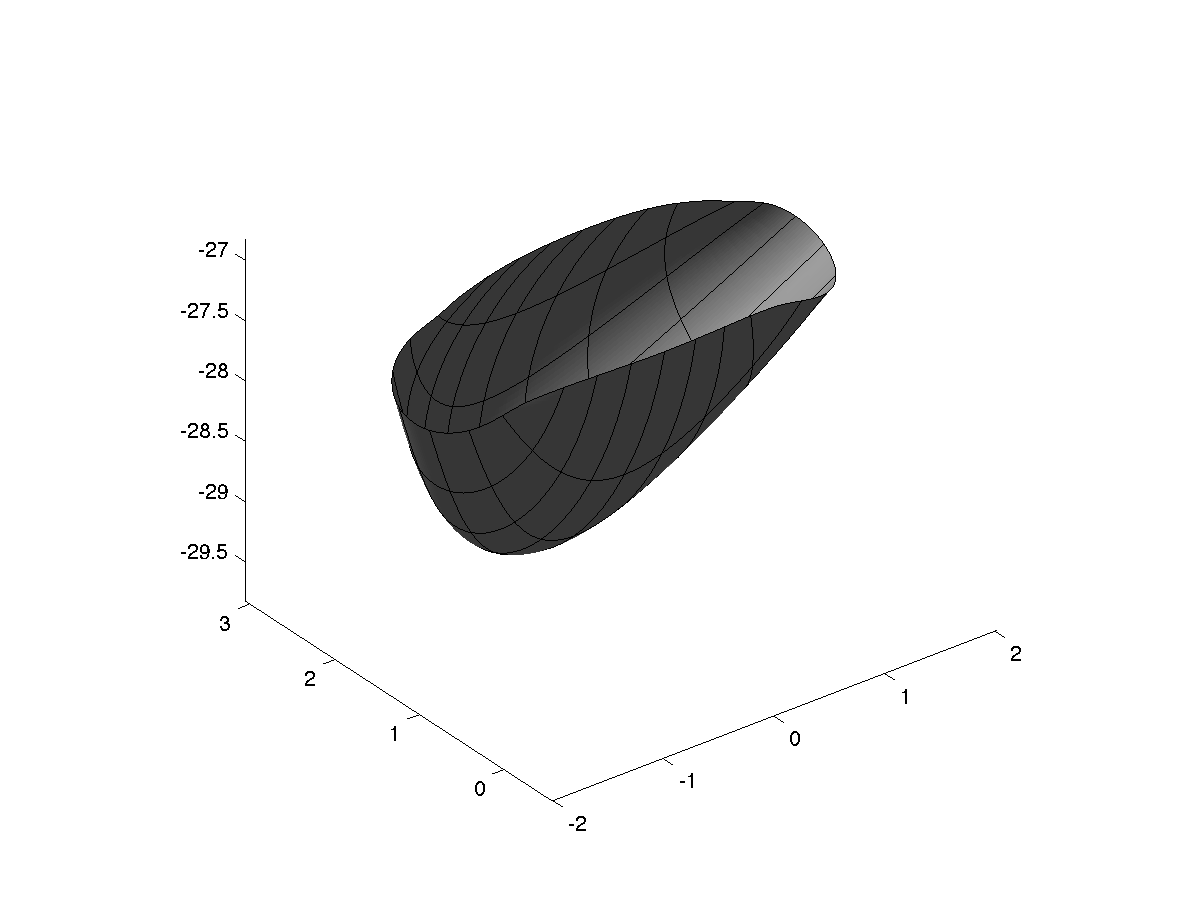}}
  \subfigure[$t=401s$]{\includegraphics[width=0.32\textwidth]{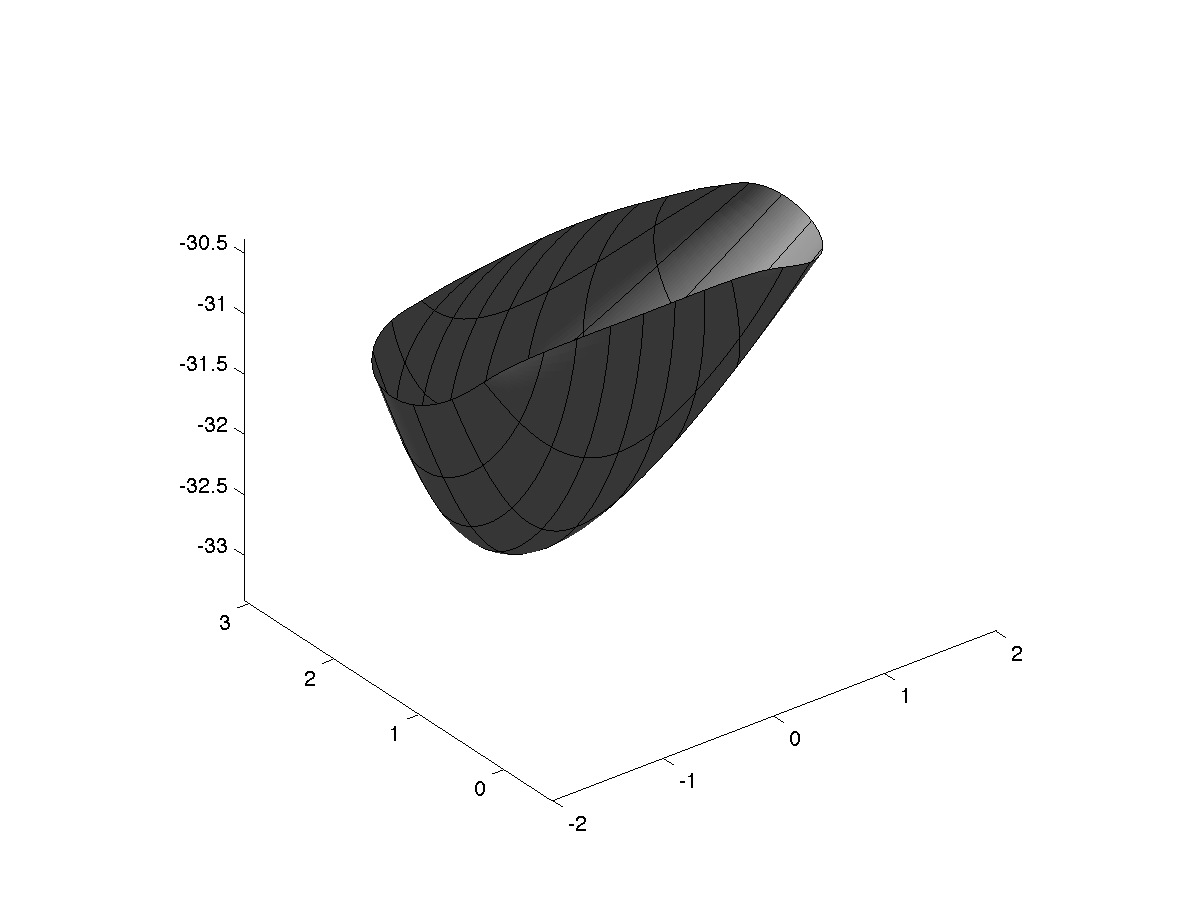}}
  \subfigure[$t=451s$]{\includegraphics[width=0.32\textwidth]{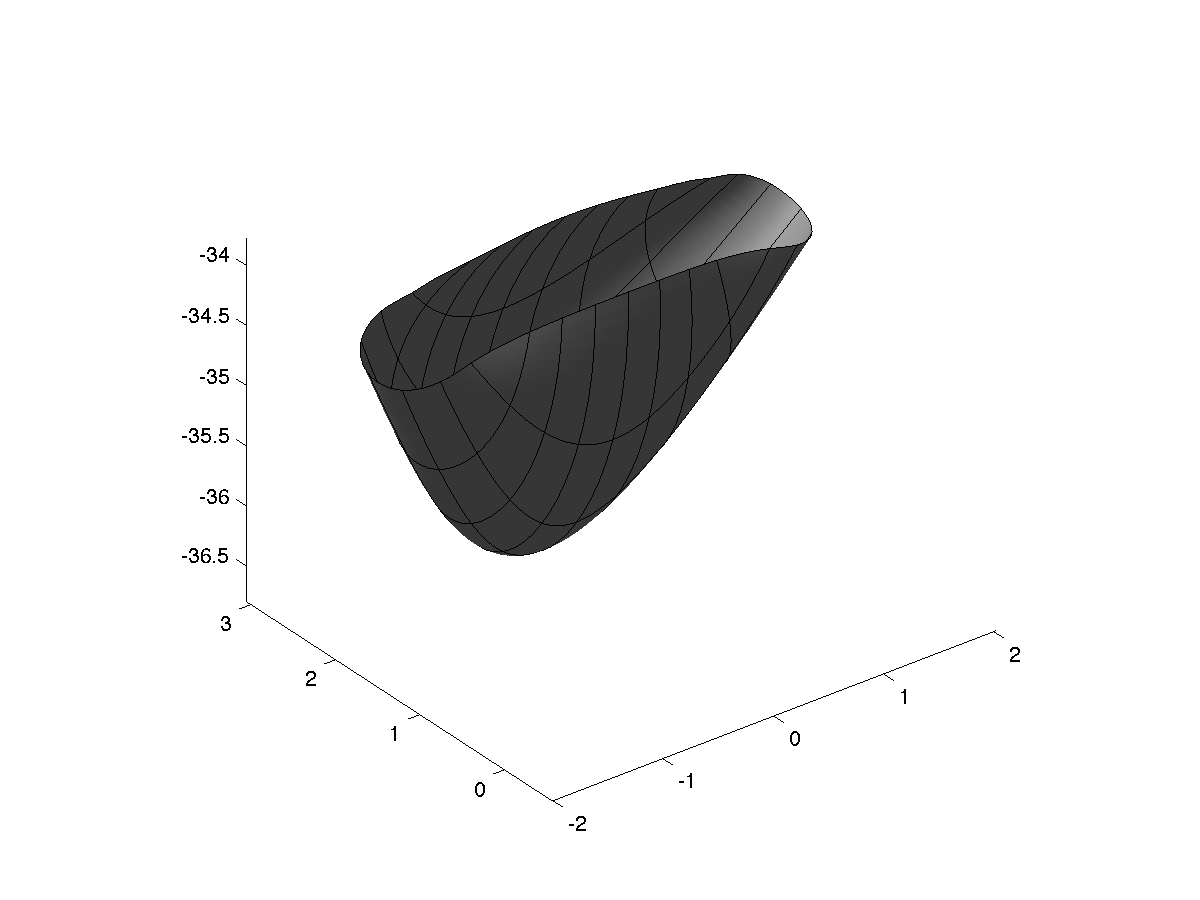}}
  \subfigure[$t=501s$]{\includegraphics[width=0.32\textwidth]{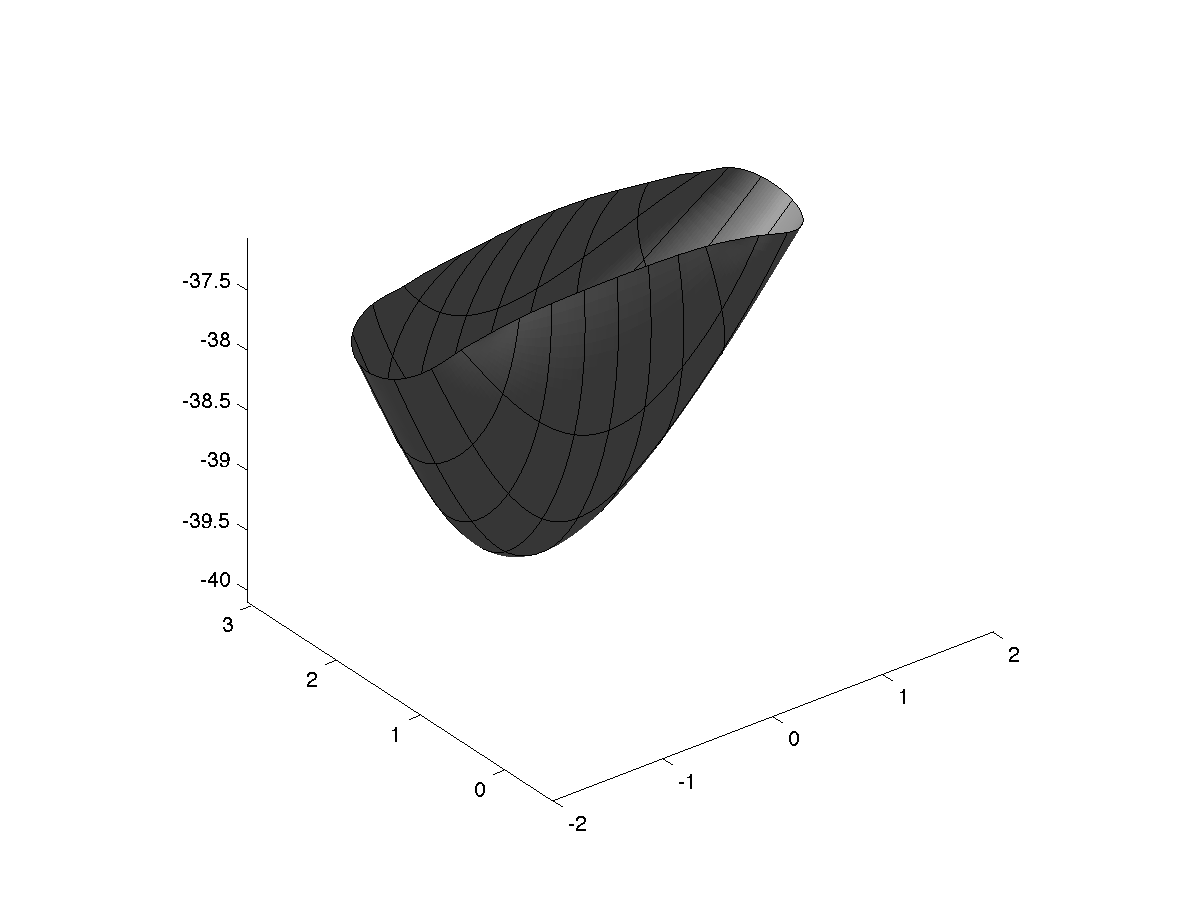}}
  \subfigure[$t=551s$]{\includegraphics[width=0.32\textwidth]{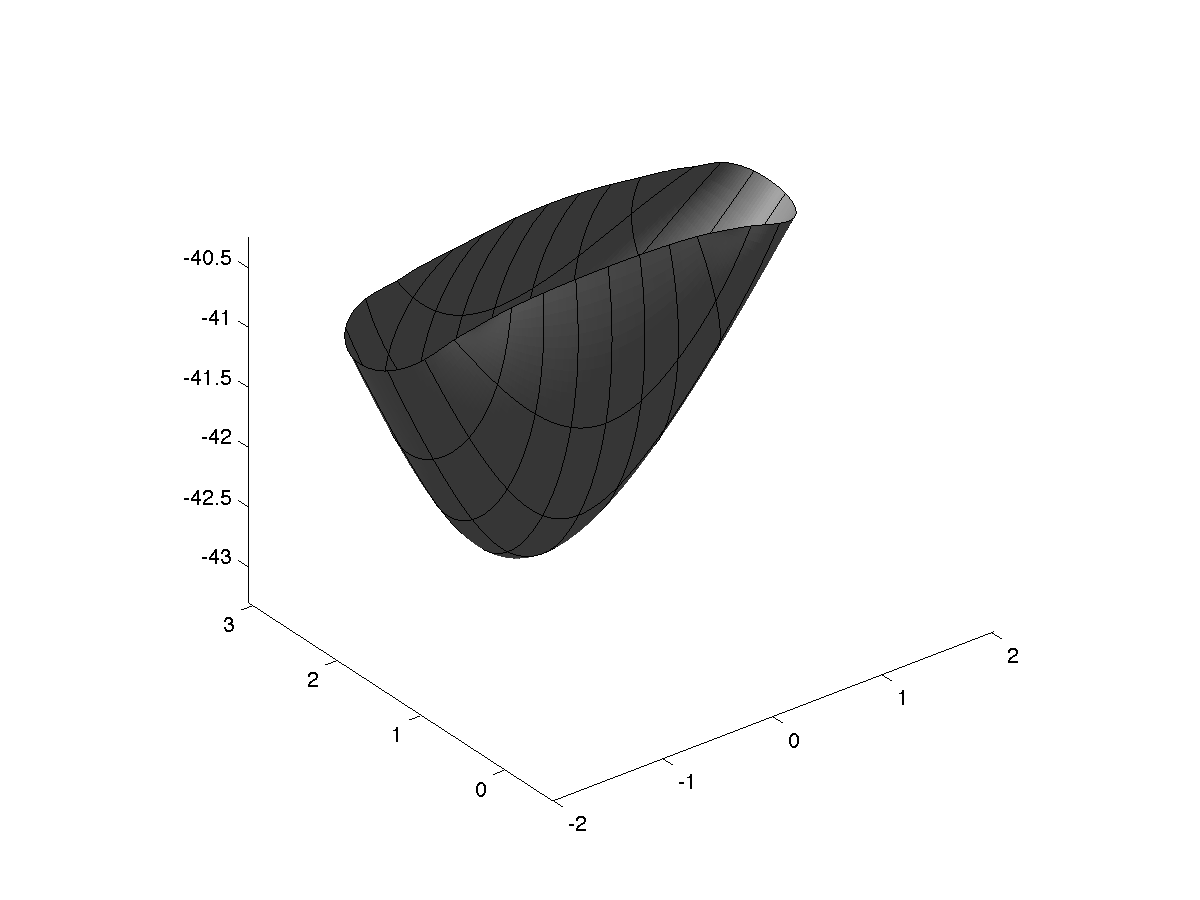}}
  \caption{Falling cap in water.} 
  \label{dome}
\end{figure}
% falling cap external flow
\begin{figure}\centering
  \subfigure[$t=101s$]{\includegraphics[width=.6\textwidth]{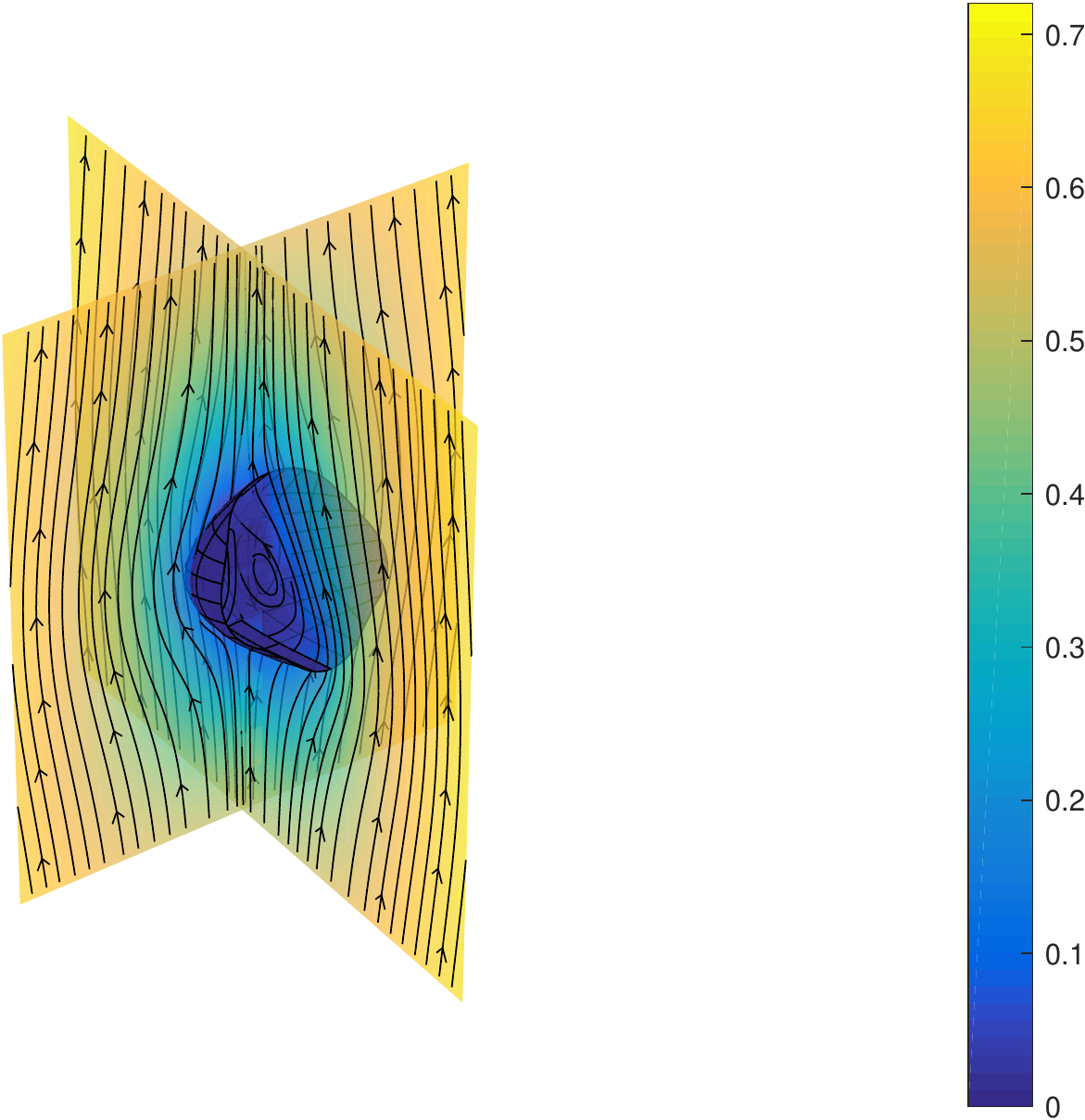}}

  \subfigure[$t=301s$]{\includegraphics[width=.6\textwidth]{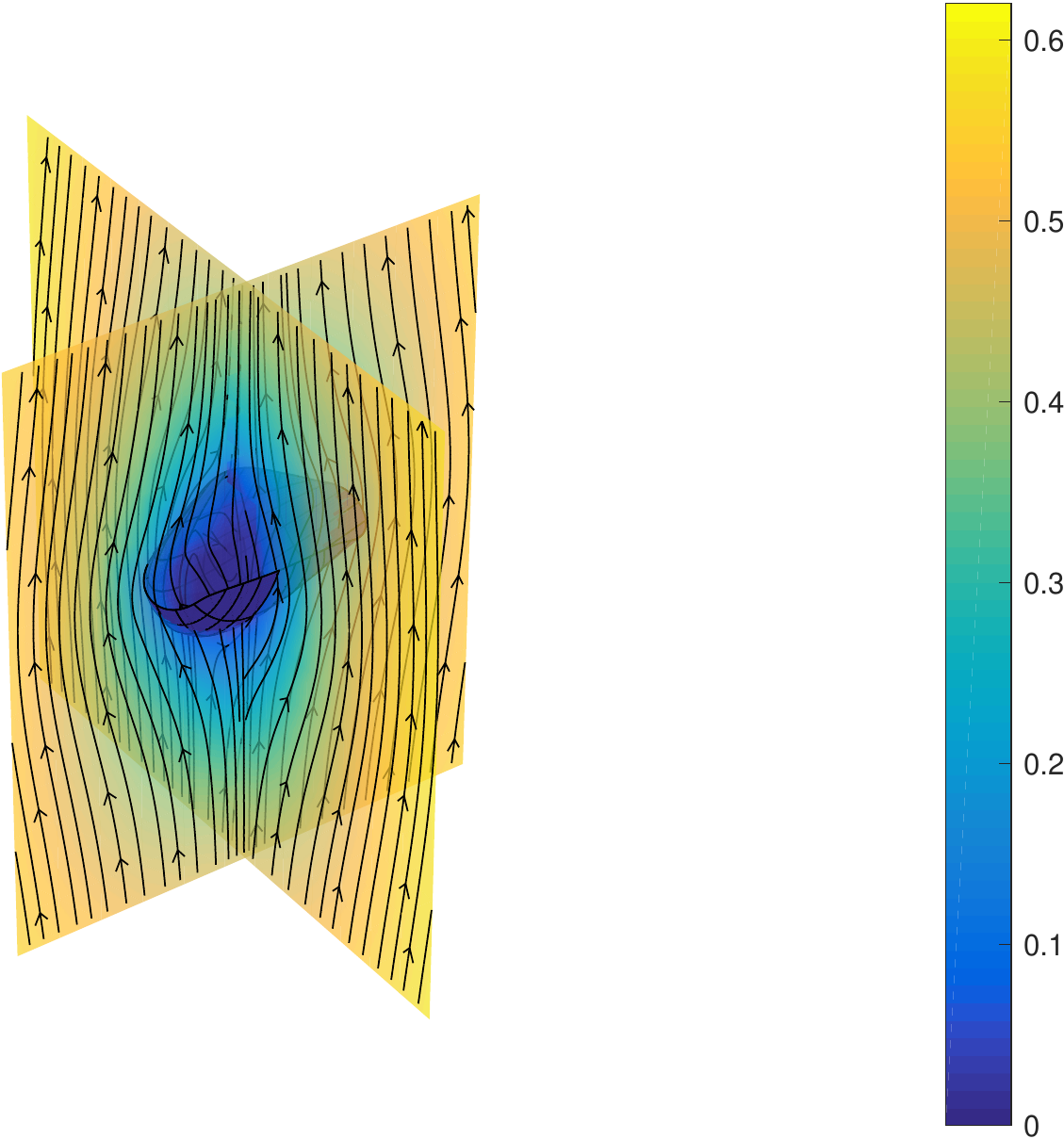}}
  \caption{Falling cap streamline plots. Streamlines and velocity magnitude plots for the relative velocity, projected on two perpendicular planes. Velocity is measured in $m/s$.} 
  \label{dome_streamlines}
\end{figure}
\subsection{Falling cap}
In this final example, we apply our method to the simulation of falling objects. Such problems are especially challenging in terms of fluid mesh generation and mesh update when using ALE approaches \cite{Franci2016520,Takizawa2012,Tezduyar2001}. 
In \cite{casquero_nurbs_2015,Hesch201251,Zhang20042051}, such problems were treated by immersed approaches, which avoid the difficulties of mesh update or remeshing but still require the discretization of a large fluid domain, whose size depends on the time interval to be observed. Using the approach presented in this paper, the computational domain is always confined to the surface model of the structure and is independent on the ``falling time".
The object under consideration is a cap as shown in Figure \ref{dome}(a). The structure's material parameters are $E=2.8\cdot10^5\, N/mm^2, \nu=0.39, \rho=1.13\, kg/m^3$ and the thickness is $h=1.0\,mm$. The structure is immersed into water with $\eta=9.0\cdot10^{-3}Pa\cdot s$ and subjected to gravity. A time step of $\Delta t=1.0\,s$ is used for the analysis and the results are depicted in Figure \ref{dome}. We can observe a combination of deformation modes, i.e., a rotation into the upright position and the deformation due to the flexibility of the structure.
Figure~\ref{dome_streamlines} presents the details of the relative flow
streamlines and relative velocity magnitude for two significant time steps.

\section{Conclusions}
\label{sec:conclusions}

In this paper, we have presented an isogeometric analysis framework to deal with a special class of FSI problems, where the fluid is represented by a Stokes flow and the structure by a shell.
The proposed framework can be defined ``truly isogeometric'' because it is entirely based on bivariate geometries, as those immediately given by CAD B-rep descriptions, given the shell nature of the considered structure and the fact that the fluid equations are solved by isogeometric boundary elements. This allows to completely circumvent the mesh generation process in many situations. In addition, the use of boundary elements may significantly limit the dimension of the discrete problem, in particular when large fluid domains are studied as in the case of falling objects.

For the solution of the coupled problem, we have chosen to adopt a semi-implicit algorithm where the effect of the surrounding flow is incorporated in the non-linear terms of the solid solver through its damping characteristics. This strategy seems capable of favourably combining the accuracy and stability of a fully implicit method with the cost-efficiency of a segregated approach.

Several numerical tests have been presented, showing the potential of the proposed analysis framework.
The extension to more complex situations, like, e.g., cell motion, and the study of the convergence properties of the model will be the subject of forthcoming research.

\section*{Acknowledgements}
LH acknowledges support by the project OpenViewSHIP, ``Sviluppo
di un ecosistema computazionale per la progettazione idrodinamica del
sistema elica-carena'', supported by Regione FVG - PAR FSC 2007-2013,
Fondo per lo Sviluppo e la Coesione and by the project ``TRIM -
Tecnologia e Ricerca Industriale per la Mobilit\`a Marina'',
CTN01-00176-163601, supported by MIUR, the italian Ministry of
Instruction, University and Research. Work by AR was supported by the
European Research Council through the FP7 Ideas Starting Grant
(project no. 259229) \emph{ISOBIO - Isogeometric Methods for
  Biomechanics}, and work by ADS was supported by the European
Research Council through AdG-340685– \emph{MicroMotility}.

\begin{appendix}

\section{Linearization of strain variables} \label{App_lin}

We begin with the linearization of the displacement vector $\us{}$, which, in discrete form, is defined as:
\begin{align}
\us{} = \sum_a^{n_{cp}} N^a \hat{\us}^a 
\end{align}
where $n_{cp}$ is the number of control points, $N^a$ are the NURBS basis functions, and $\hat{\us}^a$ are the nodal displacement vectors with components $\hat{u}^a_i \,(i=1,2,3)$ referring to the global $x-,y-,z-$components. The global degree of freedom number $\vv I$ of a nodal displacement is defined by $\vv I=3(a-1)+i$, such that $u_{\vv I}=\hat{u}^a_i$.
The variation with respect to $u_{\vv I}$ is denoted by $(\cdot),_{\vv I}$ for a compact notation and we obtain:
\begin{align}
  \frac{\partial \us{}}{\partial u_{\vv I}}=\us{},_{\vv I} = N^a \vv{e}_i
\end{align}
with $\vv{e}_i$ representing the global cartesian base vectors. For the second derivatives we obtain:
\begin{align}
  \frac{\partial^2 \us{}}{\partial u_{\vv I}\partial u_{\vv J}}=\us{},_{\vv I \vv J} = \vv{0}
\end{align}
Since variations with respect to $u_{\vv I}$ vanish for all quantities of the undeformed configuration, we obtain for the variation of $\vv{x}$:
\begin{align}
  \vv{x},_{\vv I} &= \us{},_{\vv I} = N^a \vv{e}_i \\
  \vv{x},_{\vv I \vv J}&=\us{},_{\vv I \vv J} = \vv{0}
\end{align}
Accordingly, we get the variations of the base vectors $\vv{g}_{\alpha}$ as:
% dg
\begin{align}
 \vv{g}_{\alpha},_{\vv I} &= N,_{\alpha}^a \vv{e}_i \label{aal_r}\\
 \vv{g}_{\alpha},_{\vv I \vv J}&=\vv{0} \label{aal_rs}
\end{align}
and for $\vv{g}_{\alpha,\beta}$:
\begin{align}
 \vv{g}_{\alpha,\beta},_{\vv I} &= N,_{\alpha\beta}^a \vv{e}_i \label{aalbe_r}\\
 \vv{g}_{\alpha,\beta},_{\vv I \vv J}&=\vv{0} \label{aalbe_rs}
\end{align}
With \eqref{aal_r}-\eqref{aal_rs} and $u_{\vv J}=\hat{u}^b_j$ we can express the variations of the metric coefficients $g_{\alpha\beta}=\vv{g}_{\alpha}\cdot\vv{g}_{\beta}$:
% dgab
\begin{align} \label{eq:dgab}
 g_{\alpha\beta},_{\vv I} &%= \vv{g}_{\alpha},_{\vv I}\cdot\vv{g}_{\beta} + \vv{g}_{\alpha}\cdot\vv{g}_{\beta},_{\vv I}
 = N,_{\alpha}^a \vv{e}_i\cdot\vv{g}_{\beta} + N,_{\beta}^a \vv{e}_i\cdot\vv{g}_{\alpha}  \\
 g_{\alpha\beta},_{\vv I \vv J} &= (N,_{\alpha}^a N,_{\beta}^b+N,_{\beta}^a N,_{\alpha}^b)\delta_{ij}
\end{align}
The variations of the unit normal vector $\vv{g}_{3}$ are more involved and, therefore, we introduce the auxiliary variables $\tilde{\vv{g}_{}}_3$ and $\bar{g_3}$
% g3
\begin{equation} \label{eq:g3}
\vv{\tilde g}_3 = \vv{g}_1 \times \vv{g}_2
\end{equation}
% lg3
\begin{equation} \label{eq:lg3}
\bar g_3 = \sqrt{\vv{\tilde g}_3\cdot\vv{\tilde g}_3}
\end{equation}
such that $\vv{g}_3$ can be written as:
% n
\begin{equation}
\vv{g}_3 = \frac{\vv{\tilde g}_3}{\bar g_3}
\end{equation}
In the following, we first compute the variations of the auxiliary variables which are then used for further derivations. It is convenient to follow this approach also in the implementation since these intermediate results are needed several times.
We first derive the variations of $\tilde{\vv{g}}_3$:
% dg3
\begin{align} \label{eq:dg3}
\vv{\tilde g}_3,_{\vv I} &= \vv{g}_1,_{\vv I} \times \, \vv{g}_2 + \vv{g}_1 \times \vv{g}_2,_{\vv I} \\
\vv{\tilde g}_3,_{\vv I \vv J} &= \vv{g}_1,_{\vv I} \times \, \vv{g}_2,_{\vv J} + \vv{g}_1,_{\vv J} \times \vv{g}_2,_{\vv I}
\end{align}
which are used for the variations of $\bar g_3$:
% dlg3
\begin{align} \label{eq:dlg3}
\bar g_3,_{\vv I} &= \vv{g}_{3}\cdot\vv{\tilde g}_3,_{\vv I}\\
\bar g_3,_{\vv I \vv J}  %&= \frac{\vv{\tilde g}_3,_{\vv I \vv J} \vv{\tilde g}_3 + \vv{\tilde g}_3,_{\vv I} \vv{\tilde g}_3,_{\vv J}} {\bar g_3}
%           - \frac{(\vv{\tilde g}_3,_{\vv I} \vv{\tilde g}_3) \ (\vv{\tilde g}_3,_{\vv J} \vv{\tilde g}_3)} {\bar g_3^3} \\
           &=  \bar g_3^{-1} (\vv{\tilde g}_3,_{\vv I \vv J} \cdot\vv{\tilde g}_3 + \vv{\tilde g}_3,_{\vv I} \cdot\vv{\tilde g}_3,_{\vv J} 
           - (\vv{\tilde g}_3,_{\vv I}\cdot\vv{g}_{3}) (\vv{\tilde g}_3,_{\vv J} \cdot\vv{g}_{3})
\end{align}
and finally for the variations of $\vv{g}_{3}$:
% dn
\begin{align} 
\vv{g}_3,_{\vv I} &= \bar g_3^{-1} (\vv{\tilde g}_3,_{\vv I} - \bar g_3,_{\vv I}\vv{g}_{3}) \label{a3_r}\\
\vv{g}_3,_{\vv I \vv J} & = \bar g_3^{-1} (\vv{\tilde g}_3,_{\vv I \vv J} - \bar g_3,_{\vv I \vv J}\vv{g}_{3}) +
                         \bar g_3^{-2} (2 \bar g_3,_{\vv I} \bar g_3,_{\vv J} \vv{g}_{3} - \bar g_3,_{\vv I} \vv{\tilde g}_3,_{\vv J} - \bar g_3,_{\vv J} \vv{\tilde g}_3,_{\vv I}) \label{a3_rs}
\end{align}
%The detailed steps of these derivations can be found in \cite{kiendl_isogeometric_2011}. 
With \eqref{aalbe_r}-\eqref{aalbe_rs} and \eqref{a3_r}-\eqref{a3_rs}, we can compute the variations of the curvatures $b_{\alpha\beta}=\vv{g}_{\alpha,\beta}\cdot\vv{g}_{3}$: 
% db
\begin{align}
b_{\alpha\beta},_{\vv I} &= \vv{g}_{\alpha},_{\beta},_{\vv I} \, \vv{g}_3 + \vv{g}_{\alpha},_{\beta} \, \vv{g}_3,_{\vv I} \label{balbe_r}\\
% ddb
b_{\alpha\beta},_{\vv I \vv J} &= \vv{g}_{\alpha},_{\beta},_{\vv I} \vv{g}_3,_{\vv J}
                    + \vv{g}_{\alpha},_{\beta},_{\vv J} \vv{g}_3,_{\vv I}
                    + \vv{g}_{\alpha},_{\beta} \vv{g}_3,_{\vv I \vv J} \label{balbe_rs}
\end{align}
With equations \eqref{aalbe_r}-\eqref{aalbe_rs} and \eqref{balbe_r}-\eqref{balbe_rs} we finally obtain the variations the strain variables:
% deps
\begin{align}
 \varepsilon_{\alpha\beta},_{\vv I} &= \frac{1}{2}(g_{\alpha\beta}-G_{\alpha\beta}),_{\vv I} = \frac{1}{2}g_{\alpha\beta},_{\vv I} \\
 \varepsilon_{\alpha\beta},_{\vv I \vv J} &= \frac{1}{2}g_{\alpha\beta},_{\vv I \vv J} \\
 \kappa_{\alpha\beta},_{\vv I} &= (B_{\alpha\beta}-b_{\alpha\beta}),_{\vv I} = -b_{\alpha\beta},_{\vv I}  \\
 \kappa_{\alpha\beta},_{\vv I \vv J} &= -b_{\alpha\beta},_{\vv I \vv J}
\end{align}

\section{}

\printnomenclature[2cm]

\end{appendix}

\bibliographystyle{abbrv}
\bibliography{iga_fsi_bem_shell,lit_jk_2014,lit_mch_2014,lit_fsi}

\end{document}